\documentclass[a4paper, 12pt]{article}

\usepackage[french,english]{babel} 
\usepackage[utf8]{inputenc} 

\usepackage[T1]{fontenc} 
\usepackage{csquotes} 

\usepackage{mathtools} 
\usepackage{amsthm} 
\usepackage{amssymb} 

\usepackage{mathrsfs} 
\usepackage{bbm} 
\usepackage{latexsym}
\usepackage{array} 
\usepackage{esint} 
\usepackage{stmaryrd} 
\usepackage{color}
\usepackage{graphicx}
\usepackage{verbatim} 
\usepackage{upgreek} 

\usepackage{appendix}
\usepackage{chngcntr}

\usepackage[natbib=true, style=numeric]{biblatex}
\newtheorem{dref}{Definition}[section] \newtheorem{lemma}[dref]{Lemma}
\newtheorem{theo}[dref]{Theorem} \newtheorem{prop}[dref]{Proposition}
\newtheorem{remark}[dref]{Remark}

\usepackage{anyfontsize} 
\newcommand{\e}{\mathrm{e}} 
\newcommand\iintJ{\int\hskip -2.5mm\int}

\usepackage[dvipsnames]{xcolor} 

\usepackage[pdfpagelabels,colorlinks=true,citecolor=Sepia,linkcolor=Red,urlcolor=Blue, draft=false]{hyperref}

\begin{document}


\title{\vskip-3cm {\Large\textbf{Stationary phase with Cauchy singularity.\\
}} \large{A critical point of signature {$(+,-)$}}}

 \author{Christian Klein 
 \\
 \small Université Bourgogne Europe, CNRS\\ 
   \small  IMB UMR 5584,\\
   \small 21000 Dijon, France,\\ 
  \small Institut Universitaire de France,\\
   \small  Christian.Klein@u-bourgogne.fr\\
   \and 
 Johannes Sjöstrand 
 \\ 
 \small Université Bourgogne Europe, CNRS\\ 
   \small  IMB UMR 5584,\\
   \small 21000 Dijon, France,\\ 
   \small johannes.sjostrand@u-bourgogne.fr
   \and 
   Maher Zerzeri 
   \\ 
   \small LAGA {\footnotesize - UMR7539 CNRS},\\
   \small Université Paris 13\\
   \small 99, avenue J.-B. Clément\\
   \small 93430 Villetaneuse, France\\ 
   \small zerzeri@math.univ-paris13.fr
 }
 
 
 \date{}

\maketitle

\noindent
\textbf{Abstract}: {\small Asymptotic expressions for an integral appearing in the 
solution of a d-bar problem are presented. The  integral is a solid
Cauchy transform of a function with a rapidly oscillating phase with 
a small parameter $h$, $0<h\ll 1$. 
Whereas standard steepest descent approaches can be applied to the 
case where the stationary points of the phase $\omega_{k}$, 
$k=1,\ldots, N$ are far from the 
singularity $\zeta$ of the integrand, a polarization approach is proposed for 
the case  that $|\zeta-\omega_{k}|<\mathcal{O}(\sqrt{h})$ for some $k$. In this 
case the problem is studied in $\mathbb{C}^{2}$ 
($\widetilde{\omega}:=\overline{\omega}$ is treated as an independent variable) 
on steepest descent contours. An application of Stokes' theorem allows 
for a decomposition of the integral into three terms for which 
asymptotic expressions in terms of special functions are given.}



\smallskip\noindent  
\thanks{\it This work is partially supported by 
the ANR project ISAAC-ANR-23-CE40-0015-01 and  the EIPHI Graduate School (contract ANR-17-EURE-0002)}

\tableofcontents

\section{Introduction}\label{int}
\setcounter{equation}{0}   
This paper is devoted to the asymptotic behaviour for small $h$  of the integral:
\begin{equation}\label{int.1}
I(a,\phi,\zeta;h)=-\frac{1}{\pi}\int_{\mathbb{C}}\frac{1}{\omega-\zeta} \, 
a(\omega;h) 
\e^{\frac{\mathrm{i}}{h}\phi(\omega)}\,\, {\ensuremath{L(\mathrm{d}\omega)}}\, ,\quad h\rightarrow 0\, ,
\end{equation}
where $\displaystyle {\ensuremath{L(\mathrm{d}\omega)}}=\mathrm{d}\mathrm{Re}(\omega) \wedge\mathrm{d}\mathrm{Im}(\omega)=\frac{\mathrm{d}\overline{\omega}\wedge\mathrm{d}\omega}{2\mathrm{i}}$ is the Lebesgue measure on the complex plane. Here
$\zeta\in\mathbb{C}$, $a,\phi$ are real-analytic in 
neigh$(\zeta;\mathbb{C})$\footnote{Let $M$ be a topological space. Let $N$ be a subset of $M$. The set neigh($N, M$) denotes
some neighborhood of $N$ in $M$.}. The function $a$ 
is assumed to be a symbol, the function $\phi$ to have a finite 
number of stationary points $\omega_{k}$, $k=1,\ldots,N$ with 
non-degenerate Hessian. 

Such integrals appear in the solution of d-bar problems since the solid Cauchy transform 
(\ref{int.1}) provides the solution of a d-bar equation,
$$
\frac{\partial I}{\partial\overline{\zeta}}(a,\phi,\zeta;h)= 
a(\zeta,\overline{\zeta};h) \,
\e^{\frac{\mathrm{i}}{h}\phi(\zeta,\overline{\zeta})}\,.
$$
The d-bar problems appear in the theory of integrable equations in  
2D as the Davey-Stewartson (DS) equation, see e.g. \cite{KlSa21_01},  and in 
\emph{electrical impedance tomography}  (EIT), see 
\cite{BrUh97_01,  MuSi12_01, Na96_01, Uh09_01}, and for Normal 
Matrix Models in Random Matrix Theory, see e.g.\ \cite{KlMcL17_01}. 
Since the integrals of the 
form (\ref{int.1}) are considered for all values of 
$\zeta\in\mathbb{C}$, they can be seen as a nonlinear variant of the 
Fourier transform. 

For general functions $a$, $\phi$, the integral \eqref{int.1} can 
only be computed numerically. Because of the importance of d-bar 
problems in applications there are  many 
numerical approaches to solve them. The standard method is to use 
Fourier techniques, i.e. to compute the integral after a Fourier 
transform in $\zeta$ by approximating the transform via a discrete 
Fourier transform in standard manner, see \cite{KnMuSi04_01, MuSi12_01} for a 
first order method.  
The first Fourier approach with an exponential decrease of the 
numerical error with the number of Fourier modes, or \emph{spectral 
convergence}, was presented in \cite{KlMcLSt20_01, KlMcLSt19_01, KlMcL17_01} for Schwartz class potentials via an analytical regularization of the integrand in the solid Cauchy 
transform. A spectral approach for potentials with compact support on 
a disk was presented in \cite{KlSt19_01}, a problem important in EIT applications. 
However this approach cannot be 
applied for all values of $h$ since the solution becomes too 
oscillatory for small $h$ to be computed numerically. Thus 
asymptotic formulae for small $h$ were presented in 
\cite{KlSjSt23_02,KlSjSt23_01} for functions $a$ with compact support on 
domains with a strictly convex smooth boundary. This allowed for a 
hybrid approach, the numerical computation for $h>h_{\min}$ and 
the asymptotic expressions \cite{KlSjSt23_02,KlSjSt23_01} for $h<h_{\min}$ 
with some suitably chosen $h_{\min}$ such that in both domains the 
numerical respectively asymptotic expressions are of similar accuracy. 

 In the present paper we aim at the first step in extending the 
 approach of \cite{KlSjSt23_02,KlSjSt23_01} to more general functions $a$, 
 for instance those appearing in 
 the asymptotic description of the \emph{reflection coefficient} in 
 the DS case, see Subsection \ref{sec1.1}. The basic idea is to apply 
 a standard stationary phase approximation in case $\zeta$ is not 
 close to any of the stationary points, 
 $|\zeta-\omega_{k}|>\mathcal{O}(\sqrt{h})$ for all $k=1,\ldots,N$. 
 If $|\zeta-\omega_{k}|<\mathcal{O}(\sqrt{h})$ for some 
 $k\in\{1,\ldots,N\}$, a \emph{polarization approach} is applied. This 
 means that $\overline{\omega}$ in (\ref{int.1}) is treated as an 
 independent variable $\widetilde{\omega}$, i.e. the integral 
 (\ref{int.1}) is considered in $\mathbb{C}^{2}$. After an 
 application of Stokes' theorem, suitable contours for the 
 integration are introduced in $\mathbb{C}^{2}$ via steepest descent 
 arguments. This allows to express the integrals in this case via 
 transcendental functions in leading order as in \cite{KlSjSt21_01}. 

\subsection{Background and motivation}\label{sec1.1}

One application of d-bar equations is in the inverse scattering 
approach to completely integrable equations in 2D as 
the Davey-Stewartson II (DS II) equation
\begin{equation}
\begin{split}
\mathrm{i} q_t + (q_{xx}-q_{yy}) + 2\sigma(\Phi + |q|^2)q&=0\\
\Phi_{xx}+\Phi_{yy} +2(|q|^2)_{xx}&=0\,,
\end{split}
    \label{eq:DSII}
\end{equation}
where  indices denote partial derivatives, and $q=q(x,y)$ is a 
complex-valued field. Thus equation (\ref{eq:DSII}) is  a two-dimensional nonlinear Schrödinger equation; the equation is \emph{defocusing} 
for  $\sigma=1$, and \emph{focusing} for $\sigma=-1$. Here $\Phi$ denotes a mean field. For further examples as the Kadomtsev-Petviashvili and the Novikov-Veselov 
equation, see \cite{KlSa21_01} for many references. 

\noindent
Both the scattering and the inverse scattering problem for the DS
II equation are given by the Dirac system 
\begin{equation}\label{dbarphi}
  \begin{cases}
    \overline{\partial}\upphi_{1}=\frac{1}{2}q\mathrm{e}^{\overline{k}\overline{z}-k z}\upphi_{2}\\
    \partial\upphi_{2}=\sigma\frac{1}{2}\overline{q}\mathrm{e}^{k z-\overline{k}\overline{z}}\upphi_{1}\,,\qquad \sigma=\pm1\,,
\end{cases}   
\end{equation}
subject to the asymptotic conditions 
\begin{equation}
    \lim_{|z|\to +\infty}\upphi_{1}(z)=1\,,\qquad \lim_{|z|\to +\infty}\upphi_{2}(z)=0\,,
    \label{Phisasym}
\end{equation}
where the \emph{spectral parameter} $k\in\mathbb{C}$ is independent of $z=x+\mathrm{i} y$, 
$(x,y)\in \mathbb{R}^{2}$, and 
where 
\begin{equation*}
\partial:=\frac{1}{2}\left(\frac{\partial}{\partial x}-\mathrm{i}\frac{\partial}{\partial y}\right)\quad\text{and}\quad
\bar{\partial}:=\frac{1}{2}\left(\frac{\partial}{\partial 
x}+\mathrm{i}\frac{\partial}{\partial y}\right)\,.
\end{equation*}
The system \eqref{dbarphi} appears also in the context of Calder\'on's problem, see \cite{Ca80_01}. 
The scattering data for DS, the so-called  \emph{reflection 
coefficient} $R$ are defined by 
\begin{equation}\label{reflc}
    \overline{R} = 
    \frac{2\sigma}{\pi}\int_{\mathbb{C}}^{}\e^{k z-\overline{k}\overline{z}}\,\,\overline{q}(z)\,\upphi_{1}(z;k)\,\,L(\mathrm{d}z)\,.
\end{equation}
In \cite{KlSjSt23_02,KlSjSt23_01} asymptotic formulae for large $|k|$ 
(small $h:=|k|^{-1}$) were 
derived for the reflection coefficient for the characteristic 
function of a compact domain with smooth strictly convex boundary. 
In the case of the unit disk one gets 
\begin{equation}
	R \approx \frac{2}{\sqrt{\pi 
	|k|^{3}}}\left(\sin\Big(2|k|-\frac{\pi}{4}\Big)-\frac{5}{16|k|}\cos\Big(2|k|-\frac{\pi}{4}\Big)\right)\,,
	\label{Rasym}
\end{equation}
see \cite[Subsection 4.4, page 4935]{KlSjSt23_01}, and \cite[Subsection 6.3]{KlSjSt23_02} for more details.

\noindent
The reflection coefficient evolves in time by a trivial phase factor:
\begin{equation}\label{Rt}
R(k;t)=R(k;0)\mathrm{e}^{4\mathrm{i} t\mathrm{Re}(k^2)}\,.
\end{equation}
The inverse scattering transform for DS II, i.e. the reconstruction 
of $q(z,t)$ for given $R(k;t)$ is again given by 
(\ref{dbarphi}) and
(\ref{Phisasym}) after the exchange of $k$ and $z$ and $q$ and $R$ 
respectively. Since the system (\ref{dbarphi}) 
can be written with (\ref{Phisasym}) in the form
\begin{align}
	\upphi_{1}&=1-\frac{1}{2\pi} \int_{\mathbb{C}} 
	\frac{q}{\omega-z}\exp(\overline{k} \, \overline{\omega}-k \, \omega)\upphi_{2} \, {\ensuremath{L(\mathrm{d}\omega)}}
	\nonumber\\
	\upphi_{2}&=-\frac{1}{2\pi} \int_{\mathbb{C}} 
	\frac{\bar{q}}{\overline{\omega}-\bar{z}}\exp(k \, \omega-\overline{k} \, \overline{\omega})\upphi_{1} \, {\ensuremath{L(\mathrm{d}\omega)}}\,,
	\label{phi1inv}
\end{align}
it can be solved iteratively  as in \cite{KlSjSt23_02} with the 
lowest order terms being $\upphi_{1}^{0}=1$, $\upphi_{2}^{0}=0$. This 
leads to the  lowest order non-trivial contribution to $\upphi_{2}$
\begin{equation}
	\upphi_{2}^{1}=-\frac{1}{2\pi} \int_{\mathbb{C}} 
	\frac{\bar{q}}{\overline{\omega}-\bar{z}}\exp(k \, \omega-\overline{k} \, \overline{\omega}) \, {\ensuremath{L(\mathrm{d}\omega)}}\,.
	\label{phi2inv}
\end{equation}
For the example that $\bar{q}$ is given by the asymptotic expression 
(\ref{Rasym}) and (\ref{Rt}), we can write $\upphi_{2}^{1}= 
\bar{I}_{+}+\bar{I}_{-}$, where $I_{\pm}$ are in leading order integrals of the form
(\ref{int.1}) with 
$$	
a_{\pm}=\pm\frac{1}{\sqrt{\pi |\omega|^{3}}}\,,\qquad \phi_{\pm}=
	2\tau 
	(\omega^{2}+\overline{\omega}^{2})+\frac{k\omega-\bar{k}\overline{\omega}}{\mathrm{i}}\pm 2\sqrt{\omega\overline{\omega}}\,.
$$
Thus it will be important to obtain asymptotic expression for 
integrals of the form \eqref{int.1} for small $h$ for functions $a$ 
and $\phi$ including the above examples. 

\subsection{Main results of the paper}\label{sec1.2}

We assume in this paper that the phase $\phi$ in 
(\ref{int.1}) has a finite number $N$ of stationary points $\omega_{k}$ for 
$k=1,\ldots,N$ with non-degenerate Hessian. This allows in the 
following to concentrate on the case of a single such stationary 
point that we can assume without loss of generality at the origin. In 
this paper we concentrate on the case of a Hessian at this point with 
signature $(+,-)$.

\smallskip\par 
Let $V\subset \mathbb{R}^2_{\xi ,\eta }\simeq \mathbb{C}_\omega $, $\omega =\xi
+\mathrm{i}\eta $ be an open neighborhood of $0$, let
$\phi\in C^\infty(V;\mathbb{R})$ satisfy \eqref{1.le1}--\eqref{2,5.le1}, \i.e.
\begin{equation*}
 \phi(0)=0,\qquad \phi'(0)=0\,,
\end{equation*}
\begin{equation*}
 {\rm det}\big(\phi''(0)\big)\not=0\,,
\end{equation*}
\begin{equation*}
\phi'(\omega)\not=0,\hbox{ when } \,  \omega\in V\setminus\{0\}\,. 
\end{equation*}
We are interested in the asymptotic behaviour of the expression (\ref{3.le1})
\begin{equation*}\label{1.mr1}\begin{split}
    I(a\chi, \phi,\zeta;h)
&=\frac{1}{\pi}\iintJ \frac{1}{\zeta-(\xi +\mathrm{i}\eta )}\, a(\xi ,\eta
;h)\chi(\xi ,\eta )\, \e^{\frac{\mathrm{i}}{h}\phi(\xi ,\eta )}d\xi d\eta  
    \\
    &=-\frac{1}{\pi}\int \frac{1}{\omega-\zeta}\, a(\omega;h)\chi(\omega)\, \e^{\frac{\mathrm{i}}{h}\phi(\omega)}\, L(\mathrm{d}\omega)\, ,
\end{split}
\end{equation*}
where $a(\xi ,\eta ;h)=a(\omega;h)\in C^\infty(V\times ]0,h_0])$ is a classical symbol
of order $h^0$, i.e. such that we have (\ref{4.le1}):
\begin{equation*}
 a(\omega;h)\sim a_0(\omega)+h a_1(\omega)+\ldots\hbox{ in
 }C^\infty(V)\,,\quad  h\longrightarrow 0\,,
\end{equation*}
and $\chi(\omega)\in C_0^\infty(V)$ is a fixed cutoff function, equal
to 1 near $\omega=0$. Here we also assume that $\zeta \in
\mathrm{int}(\chi ^{-1}(1))$.

\smallskip
The asymptotic behaviour is generated in a small neighborhood of
$\{0,\zeta  \}$, where $0$ is the critical point of $\phi (\omega )$
and $\zeta $ the locus of the pole singularity of the integrand. When
$|\zeta |\ge 1/{\mathcal O}(1)$\footnote{Here we follow the convention that the expression “$\mathcal{O}(1)$” in a denominator denotes a
bounded positive quantity.}, we can decompose the integral with
a partition of unity and study the contribution from $\omega =0$ with
the stationary phase method and the one from $\omega =\zeta $,
observing that the corresponding contribution is of the form $c(\zeta
;h)\e^{\mathrm{i}\phi (\zeta )/h}$, where $c(\zeta ;h)\sim h(c_1(\omega
)+h c_2(\omega )+... )$ solves asymptotically
$$
\partial _{\overline{\omega  }}\Big(c\, \e^{\mathrm{i}\frac{\phi}{h}}\Big)=\big(a(\omega ;h)+{\mathcal
  O}(h^\infty )\big)\,\e^{\mathrm{i}\frac{\phi}{h}}\,,\quad \omega \in \mathrm{neigh}(\zeta;\mathbb{C})\,,
$$
by means of an iterated division procedure, starting with $c_1(\omega
)\mathrm{i}\partial_{\overline{\omega }}\phi (\omega )=a_0(\omega )$. 
Here $f(h)=\mathcal{O}(h^\infty)$ means that for all integer $N\in \mathbb{N}$ there exist $C_N>0$ such that $|f(h)|\leq C_N h^N$, $\forall h\in ]0,h_0]$.
This leads to (\ref{5.le2}):
\begin{equation*}
I(a\chi,\phi,\zeta;h)=c(\zeta;h)\, \e^{\frac{\mathrm{i}}{h}\phi(\zeta)}+
b(\zeta;h)+\mathcal{O}(h^\infty)\,,
\end{equation*}
where the symbols $b$, $c$ are order $h$ and satisfy (\ref{1.le2}), (\ref{3.le2}).

\medskip
In the following we study the case when $|\zeta |\ll 1$. Put
\begin{equation}\label{1.mr2}
\varepsilon :=\frac{|\zeta |}{\sqrt{h}}
\end{equation}
and assume first that $\varepsilon$ satisfies \eqref{5.le1}, i.e.
\begin{equation*}\label{2.mr2}
\varepsilon \gg 1\,.
\end{equation*}
Write $\zeta =\lambda \widetilde{\zeta }$, where $\lambda \asymp
|\zeta |$ is an extra parameter, and define $\widetilde{\omega }$ by
$\omega =\lambda \widetilde{\omega }$.
Write as in (\ref{2.le4}), 
\begin{equation*}
 \frac{1}{h}\phi(\omega)=\frac{\lambda^2}{h}\frac{1}{\lambda^2}\phi(\lambda\widetilde\omega)=\frac{1}{\widetilde h}\widetilde\phi(\widetilde\omega)\, ,
\end{equation*}
where as in (\ref{3.le4})
\begin{equation*}
\widetilde h=\frac{h}{\lambda^2}\asymp\frac{1}{\varepsilon ^2}\ll 1, \qquad 
\widetilde\phi(\widetilde\omega)=\frac{1}{\lambda^2}\phi(\lambda\widetilde \omega)\, .
\end{equation*}
Making the same changes of variables in (\ref{5.le2}), we will get
(\ref{6.le4}):
\begin{equation*}
 I(a\chi,\phi,\zeta;h)=\mathcal{O}\left(\widetilde h^\infty\right)+\lambda \left(\widetilde c(\widetilde\zeta;\widetilde h)\, \e^{\frac{\mathrm{i}}{\widetilde h}\tilde\phi(\widetilde\omega)}+\widetilde b(\widetilde\zeta;\widetilde h)+\mathcal{O}\big(\widetilde h^\infty\big)\right)\,,
\end{equation*}
where $\widetilde c(\widetilde\zeta;\widetilde h), \widetilde b(\widetilde\zeta;\widetilde h)$ are classical symbols of order $\widetilde h$, satisfying the analogues of \eqref{2,5.le2}, \eqref{9.le1}:
\begin{equation*}
\left \{
 \begin{array}{cc}
    \widetilde c(\widetilde\zeta;\widetilde h)\sim 
    \widetilde h\big(\widetilde c_1(\widetilde\zeta)+\widetilde h\widetilde c_2(\widetilde\zeta)+\ldots\,\big)\,,\\
     {}  {} \\
     \widetilde b(\widetilde\zeta;\widetilde h) \sim 
    \widetilde h\big(\widetilde b_1(\widetilde\zeta)+\widetilde h\widetilde b_2(\widetilde\zeta)+\ldots\,\big)\,.
\end{array}
\right.
\end{equation*}
On the formal asymptotic level we have (\ref{3.le5}):
\begin{equation*}
  \lambda\widetilde c(\widetilde\zeta;\widetilde h)=c(\zeta;h)\,,\qquad 
  \lambda\widetilde b(\widetilde\zeta;\widetilde h)=b(\zeta;h)\,.
\end{equation*}
See Section \ref{le}, in particular Proposition \ref{1le5}. 

\medskip
The main part of this work is about the case when
\begin{equation}\label{3.mr2}
\varepsilon \le h^{-\delta }
\end{equation}
for some fixed $\delta \in ]0,1/6[$. We then need to go out
to the complex domain and assume that $a(\xi ,\eta ;h)$ and $\phi (\xi
,\eta )$ extend to
holomorphic functions in a neighborhood of $0$ in ${\mathbb{C}}^2_{\xi
  ,\eta }$ that is independent of $h$. We assume that
$\mathrm{supp}\chi $ is contained in that neighborhood.

\medskip
We will use the method of polarization and start by
  identifying $\mathbb{R}^2_{\xi ,\eta }$ with the antidiagonal in ${\mathbb
    C}^2_{\omega ,\widetilde{\omega }}$ as in (\ref{2.sto3}),
  \begin{equation*}
    \mathrm{adiag\,}(\mathbb{C}_{\omega ,\tilde{\omega }}^2)=
    \{ (\omega ,\widetilde{\omega })\in \mathbb{C}^2_{\omega
      ,\tilde{\omega }};\,\,\, \widetilde{\omega }=\overline{\omega } \}\,.
  \end{equation*}
  Here we put as in (\ref{2,5.sto3})
  \begin{equation*}
  \omega =\xi +\mathrm{i}\eta\,, \hbox{ and also }
  \widetilde{\omega}=\xi -\mathrm{i}\eta\,. 
  \end{equation*} 
Since
$$
\mathrm{d}\xi \wedge \mathrm{d}\eta = \frac{\mathrm{d}\widetilde
  {\omega}\wedge \mathrm{d}\omega}{2\mathrm{i}}\,, 
$$
the integral in (\ref{1.sto3}) can be written as in (\ref{3.sto3})
\begin{equation*}
I(a\chi ,\phi ,\zeta ;h)=-\frac{1}{\pi}\iintJ_{\mathrm{adiag\,}({\mathbb{C}}^2)}\frac{1}{\omega-\zeta} \, a(\omega ,\widetilde{\omega };h)\chi (\omega
,\widetilde{\omega })\, \e^{\frac{\mathrm{i}}{h}\phi (\omega ,\widetilde{\omega })} \, \frac{\mathrm{d}\widetilde{\omega }\wedge \mathrm{d}\omega}{2\mathrm{i}}\,, 
\end{equation*}
abusing slightly the notation, letting $a$, $\chi$, $\phi $ denote
``the same'' functions on $\mathrm{adiag\,}(\mathbb{C}^2)$ using
(\ref{2,5.sto3}) and the inverse relations (\ref{4.sto3})
\begin{equation*}
\xi =\frac{\omega +\widetilde{\omega}}{2}\,,\quad \eta =\frac{\omega
  -\widetilde{\omega }}{2\mathrm{i}}\,.
\end{equation*}

\par
The idea is then to deform the contour
$\mathrm{adiag}({\mathbb{C}}^2_{\omega ,\tilde{\omega }})$ to a new contour in
${\mathbb{C}}^2_{\omega ,\tilde{\omega }}$ which is a Cartesian product
$\gamma _1\times \gamma _2$, where $\gamma _1\subset {\mathbb{C}}_\omega $ and
$\gamma _2\subset {\mathbb{C}}_{\tilde{\omega }}$ are curves. Stokes' formula
will then produce several terms to be analyzed, in particular the
integral along $\gamma _1\times \gamma _2$ which reduces to 1D
integrals (see (\ref{1.sto4})):
\begin{equation*}
\iintJ_{\gamma _1\times \gamma _2}(...)\mathrm{d}\widetilde{\omega} \wedge
\mathrm{d}\omega =\int_{\gamma _1}\left(\int_{\gamma _2}(...)
\mathrm{d}\widetilde{\omega } \right) \mathrm{d}\omega\,.
\end{equation*}

We define a smooth family of contours $\Gamma_\theta$, $\theta \in [0,1]$, of real dimension 2 in ${\mathbb{C}}^2_{\xi ,\eta }$ or equivalently in $\mathbb{C}^2_{\omega
  ,\tilde{\omega }}$. In the $\mathbb{C}^{2}_{\xi ,\eta }$ representation it is
given by
$$
\mathbb{R}^2\ni (t,s)\longmapsto (\xi (\theta ,t,s),\eta (\theta ,t,s))\in
\mathbb{C}^2\,,
$$
where as in (\ref{Sto.3})
\begin{equation*}
\begin{cases}
 \xi(\theta;t,s) = (1+\mathrm{i}\theta) t\\
 \eta(\theta;t,s) = (1-\mathrm{i}\theta)s\,,\\
\end{cases}
\quad 
(t,s)\in \mathbb{R}^2\,.
\end{equation*}
In the ${\mathbb{C}}_{\omega ,\widetilde{\omega }}$ representation it is
given as in (\ref{Sto.5}) by
\begin{equation*}
\mathbb{R}^2\ni(t,s)\,\longmapsto f(\theta;t,s)=\big(\omega(\theta;t,s)\,,\,\widetilde\omega(\theta;t,s)\big)\in\mathbb{C}^2_{\omega,\tilde\omega}\,,
\end{equation*}
where as in (\ref{Sto.4})
\begin{equation*}
\left\{
\begin{array}{rcl}
\omega(\theta;t,s)&:= &\xi(\theta;t,s)+\mathrm{i}\eta(\theta;t,s)\\
     &= &(t+\mathrm{i} s)+\mathrm{i}\theta(t-\mathrm{i} s)\\
\widetilde\omega(\theta;t,s)&:= &\xi(\theta;t,s)-\mathrm{i}\eta(\theta;t,s)\\
      &= &(t-\mathrm{i} s)+\mathrm{i}\theta(t+\mathrm{i} s)\\
\end{array}
\right.,\quad (t,s)\in \mathbb{R}^2\,.
\end{equation*}

\medskip\par 
In the $\mathbb{C}^2_{\xi ,\eta}$ representation, we have $\Gamma
_\theta =(1+\mathrm{i}\theta )\mathbb{R}\times (1-\mathrm{i}\theta )\mathbb{R}$, especially $\Gamma_0=\mathbb{R}^2$.
In the $\mathbb{C}^2_{\omega ,\tilde{\omega }}$ representation, we
get $\Gamma _0=\mathrm{adiag\,}(\mathbb{C}^2)$, and for $\theta =1$, we
have
$$
\omega (1,t,s)=(1+\mathrm{i})(t+s),\quad \widetilde{\omega }(1,t,s)=(1+\mathrm{i})(t-s)\,,
$$
so $\Gamma _1=(1+\mathrm{i})\mathbb{R}\times (1+\mathrm{i})\mathbb{R}$ is a
Cartesian product. 

\medskip\par 
In Proposition \ref{1sto7} and the subsequent identities we
establish a Stokes' formula for the transformation of (\ref{3.sto3})
and get as in (\ref{2.sto8})
\begin{align}
-\frac{1}{\pi}&\iintJ_{\mathrm{adiag\,}({\mathbb{C}}^2)} \frac{1}{\omega-\zeta} \, a(\omega ,\widetilde{\omega };h)\chi (\omega,\widetilde{\omega })\, \e^{\frac{\mathrm{i}}{h}\phi (\omega ,\widetilde{\omega })} 
\, \frac{\mathrm{d}\widetilde{\omega}\wedge \mathrm{d}\omega}{2\mathrm{i}}\nonumber\\ 
=&-\frac{1}{\pi}\iintJ_{\Gamma_{1-\delta}}\frac{1}{\omega-\zeta} \, 
a(\omega ,\widetilde{\omega };h)\chi(\omega,\widetilde{\omega}) \, \e^{\frac{\mathrm{i}}{h}\phi(\omega ,\widetilde{\omega })} \,
\frac{\mathrm{d}\widetilde{\omega }\wedge \mathrm{d}\omega}{2\mathrm{i}}\nonumber
\\
&-\int_0^{1-\delta} \frac{2\mathrm{i}\,\overline{\widetilde{\zeta}}}{1-\theta ^2} \, a(\zeta,\widetilde{\zeta };h)\chi (\zeta ,\widetilde{\zeta
})\, \e^{\frac{\mathrm{i}}{h}\phi (\zeta ,\widetilde{\zeta})} \, \mathrm{d}\theta \nonumber\\
&+\int_0^{1-\delta}\hskip-5pt\frac{1}{\pi}\hskip-5pt\iintJ_{\Gamma_\theta}
\frac{1}{\omega -\zeta } \, a(\omega ,\widetilde{\omega
};h)\, \e^{\frac{\mathrm{i}}{h}\phi (\omega ,\widetilde{\omega })}\widetilde\nu
(\theta;\omega,\widetilde\omega)\rfloor
\mathrm{d}_{\omega ,\tilde{\omega}}\big(\chi(\omega ,\widetilde{\omega}) \frac{\mathrm{d}\widetilde{\omega }\wedge
\mathrm{d}\omega}{2\mathrm{i}}\big)\mathrm{d}\theta \nonumber\\
&=:\mathrm{I}(\zeta ,1-\delta)+\mathrm{II}(\zeta ,1-\delta)+\mathrm{III}(\zeta,1-\delta)\,.\nonumber
\end{align}
Here $\widetilde{\nu }$ is a deformation vector field, as explained in
Proposition \ref{1sto7}. In the applications we will be able to pass to
the limit $\delta =0$, see (\ref{3.sto8}).

\par\medskip
We now assume that $\phi$ is holomorphic in
$\mathrm{neigh}((0,0);\mathbb{C}^2)$ (either in $\mathbb{C}^2_{\xi ,\eta
}$ or in $\mathbb{C}^2_{\omega ,\tilde{\omega }}$, depending on the
which coordinates we use). Assume (\ref{6.ft3})
\begin{equation*}
\phi (0,0)=0,\ \ \nabla \phi (0,0)=0\,.
\end{equation*}
Let $\phi_2(\xi ,\eta )$ be the 2nd order Taylor polynomial of $\phi
$ at $(0,0)$, then as in (\ref{1.ft4})
\begin{equation*}
\phi (\xi,\eta )=\phi_2(\xi ,\eta )+{\mathcal  O}(|(\xi ,\eta )|^3)\,.
\end{equation*}
We use the same notations for functions in the variables $\omega
,\widetilde{\omega }$,
\begin{equation*}
  \phi (\omega ,\widetilde{\omega } )=\phi_2(\omega ,\widetilde{\omega
  })+{\mathcal  O}(|(\omega ,\widetilde{\omega })|^3)\,.
\end{equation*}
Assume that $\phi_2$ is real-valued on $\mathbb{R}^2_{\xi ,\eta }$,
hence as in (\ref{2.ft4})
\begin{equation*}
\phi_2(\xi ,\eta )=\frac{\lambda }{2}\xi ^2-\frac{\mu }{2}\eta
^2+\rho\, \xi \eta ,\ \ \lambda ,\mu ,\rho \in \mathbb{R}\,.
\end{equation*}
Assume (\ref{3.ft4}):
\begin{equation*}
  \lambda,\, \mu >0. \end{equation*}

We assume (\ref{5.ft7}):
\begin{equation*}
a(\omega ,\widetilde{\omega })={\mathcal  O}(1)\hbox{ is holomorphic in }\mathrm{neigh}\big((0,0);\mathbb{C}^2\big)\,.
\end{equation*}
We allow $a$ to depend on $h$ and then assume that $a(\omega
,\widetilde{\omega };h)$ is a classical symbol of order
$h^0$, holomorphic in $(\omega ,\widetilde{\omega})$, in the sense that 
\begin{equation}\label{1.mr4}
a(\omega ,\widetilde{\omega }; h) \sim \sum_{j=0}^{+\infty} a_j(\omega ,\widetilde{\omega }) h^j \hbox{ in }  \, C^\infty\Big(\mathrm{neigh}\big((0,0);\mathbb{C}^2\big); \mathbb{C}\Big)\,,
\end{equation}
where $a, a_j$ are holomorphic in the same $\mathrm{neigh}\big((0,0);\mathbb{C}^2\big)$.
We will call such a symbol a  \textit{holomorphic classical symbol.}

From now on, we will work with the class of holomorphic classical symbols. Note that this framework lies between the frameworks of classical symbols and analytic symbols. This framework is sufficient in this paper, since our goal is to study the asymptotic behaviour of the integral \eqref{int.1}, using the complex contour deformation method, but it suffices to allow error terms of order $\mathcal{O}(h^\infty)$. Here we do not need a certain stability as in the framework of pseudo-differential operators with analytic symbols, for example for the composition or invertibility of these operators; in particular, there is no need to introduce a family of pseudo-norms. For more details on the framework of classical symbol spaces and  analytic symbol spaces, see \cite{DiSj99_01} and \cite{Sj82_01}, see also \cite{Bo76_01, BoKr67_01, Sj80_01, Sj80_02}.

\medskip\par
We assume that $\chi $ is supported in the neighborhood in
(\ref{5.ft7}) and that (\ref{6.ft7}) holds,
\begin{equation*}
\chi =1\hbox{ in }\mathrm{neigh}\big((0,0);\mathbb{C}^2_{\omega
  ,\tilde{\omega }}\big)\,.
\end{equation*}
Assume (\ref{3.mr2}):
$$
|\zeta |\le h^{\frac{1}{2}-\delta }\hbox{ for some fixed }\delta \in ]0,1/6[\,.
$$
Assume that
\begin{equation}\label{1.mr5}\begin{split}
&\zeta \ne 0\hbox{ and that }
\zeta \, \hbox{ is not in some fixed  conic}\\ &\hbox{neighborhood of }
\e^{\mathrm{i}\frac{\pi}{4}}\left( {\mathbb R}\setminus \{ 0 \}
\right)\hbox{ in }{\mathbb{C}}\setminus \{0 \}\,.
\end{split}
\end{equation}
Here we notice that if $\zeta \ne 0$ belongs to a small conic
neighborhood of $\e^{\mathrm{i}\frac{\pi}{4}}{\mathbb{R}}\setminus \{0 \}$ in ${\mathbb
  C}\setminus \{0 \}$, then we can make a small rotation: $\xi +\mathrm{i}\eta
=\e^{\mathrm{i}\theta }(\check{\xi }+i\check{\eta })$ (equivalently $\omega
=\e^{\mathrm{i}\theta }\check{\omega }$, $\widetilde{\omega }=\e^{-\mathrm{i}\theta
}\check{\widetilde{\omega }}$).
Replace $\omega $ with $\check{\omega }$ in (\ref{1.mr1}) and notice that
$$
\frac{1}{\omega -\zeta }=\e^{-\mathrm{i}\theta }\frac{1}{\check{\omega
  }-\check{\zeta  }}\,,\quad \check{\zeta }=\e^{-\mathrm{i}\theta }\zeta\,,
$$
so $\check{\zeta }$ avoids a small conic neighborhood of $\e^{\mathrm{i}\frac{\pi}
 {4}}{\mathbb{R}}$ in ${\mathbb{C}}\setminus \{0 \}$. The slightly rotated new
exponent $\phi (\e^{-\mathrm{i}\theta }\check{\omega })$ still satisfies
(\ref{2.ft4}).
\begin{prop}\label{1mr6}
Under the above assumptions the integrals $\mathrm{I}(\zeta ,1)$, $\mathrm{II}(\zeta
,1)$, $\mathrm{III}(\zeta ,1)$ are well defined and respectively equal to 
$\underset{\delta \to 0}{\lim}\, \mathrm{I}(\zeta ,1-\delta )$, $\underset{\delta \to 0}{\lim} \, \mathrm{II}(\zeta ,1-\delta )$, $\underset{\delta \to 0}{\lim} \, \mathrm{III}(\zeta ,1-\delta )$, so by (\ref{2.sto8}),
$$
I(a\chi ,\phi ,\zeta ;h)=\mathrm{I}(\zeta ,1)+\mathrm{II}(\zeta
,1)+\mathrm{III}(\zeta ,1)\,.
$$
\end{prop}

\par
To describe the asymptotic behavior of $I(a\chi ,\phi ,\zeta ;h)$ it suffices to do so for $\mathrm{I}(\zeta ,1)$, $\mathrm{II}(\zeta
,1)$, $\mathrm{III}(\zeta ,1)$.

\par\medskip\noindent
\textbf{\large 1.} {\bf The term $\mathrm{ I } (\zeta ,1)$} is studied in Section \ref{ft}. Since
$\Gamma_1$ is a Cartesian product, we can study the first integral
with respect to $\widetilde{\omega }$ and apply stationary phase. This
leads to (\ref{1.ft9})
\begin{equation*}
\mathrm{I}(\zeta ,1)=\frac{1}{2\pi \mathrm{i}}\int_{\e^{\mathrm{i}\frac{\pi}{4}}\mathbb{R}}
\frac{1}{\zeta -\omega }h^{\frac12}b(\omega ;h)\, \e^{\frac{\mathrm{i}}{h}\psi (\omega)} \, \mathrm{d}\omega\,,
\end{equation*}
where up to an exponentially small term $b$ is an analytic symbol
truncated to a neighborhood of $\omega =0$ and $\psi (\omega )$ is
holomorphic near $0$ with Taylor expansion
$$
\psi (\omega )=\psi _2(\omega )+{\mathcal O}(\omega ^3)\,,\quad \psi _2(\omega
)=f\, \frac{\omega^2}{2}\,,
$$
where $f$ is given by (\ref{4,5.ft10}):
\begin{equation*}
f=\frac{\lambda \mu +\rho ^2}{\lambda +\mu +2\mathrm{i}\rho }\,.
\end{equation*}

\noindent
In and around (\ref{1.ft12})--(\ref{5'.ft13}) we discuss the entire
special functions related to  Dawson's integral. In the result below, these
functions are used only in regions where they are not exponentially
large. Our main result for the term $\mathrm{I}(\zeta ,1)$ is given in
Proposition \ref{1ft15}:

\medskip\noindent
{\it Assume that $|\zeta |\le h^{\frac{1}{2}-\delta }$ and that
(\ref{1.mr5}) holds. Define $\widehat\zeta$ and $\zeta ^\dagger$ with $|\widehat\zeta|, \, |\zeta ^\dagger|\le
{\mathcal O}(h^{-\delta })$, by $\widehat{\zeta}:=h^{-1/2}\zeta =:(\mathrm{i}/f)^{1/2}\zeta
^\dagger$. Then there exists an asymptotic series
$c(\widehat{\omega };h)\sim c_0(\widehat{\omega
})+h c_1(\widehat{\omega })+...$, where $c_j(\widehat{\omega })$ is a
polynomial of degree $\le 3j$ (independent of $\zeta $) such that
(\ref{3.ft15}) holds,
\begin{equation*}
  \begin{split}
\mathrm{I}(\zeta ,1)  
&=h^{\frac12}c(f^{-1}D_{\widehat{\zeta
  }};h)\left(G_{\mathrm{i}\psi_2}(\widehat{\zeta}) \right)
+{\mathcal O}(h^\infty )\\
&=
  h^{1/2}c((\mathrm{i}/f)^{1/2}\partial _{\zeta^\dagger} ;h)\big(G^{l/r}_{-\frac{\omega^2}{2}}
(\zeta^\dagger )\big)+{\mathcal O}(h^\infty )\,.
  \end{split}
\end{equation*}
The functions $G_\rho =G_\rho ^{l/r}$ are defined by line integrals in
(\ref{8.ft10}), (\ref{2.ft12}) and we take $G_{..}=G_{..}^l$ when $\mathrm{Im}
(\e^{-\mathrm{i}\frac{\pi}{4}}\widehat{\zeta })>0$ and
$G_{..}=G_{..}^r$ when $\mathrm{Im}
(\e^{-\mathrm{i}\frac{\pi}{4}}\widehat{\zeta })<0$.}
\par\medskip


Assume that $|\zeta |\le h^{\frac{1}{2}-\delta }$. Then, by \eqref{9.ft7}, \eqref{9.ft5}, the leading term of the $\mathrm{I}(\zeta ,1)$ is given by:
\begin{equation}\label{leadI}
h^{\frac12}\,\frac{2\sqrt{2\pi} \, \e^{\mathrm{i}\frac{\pi}{4}}}{(\lambda+\mu+\mathrm{i}2\rho)^{\frac12}}\, a(0,0) \, G_{-\frac{\omega^2}{2}} \Big(h^{-\frac12}\big(\frac{\mathrm{i}}{f}\big)^{-\frac12}\zeta\Big) \,, 
\end{equation}
where 
\begin{equation*}
 G_{-\frac{\omega^2}{2}}(\zeta)=\begin{cases}
  G_{-\frac{\omega^2}{2}}^l(\zeta) \quad \mathrm{ when } \quad \mathrm{Im}
(\e^{-\mathrm{i}\frac{\pi}{4}}\zeta)>0\\
G_{-\frac{\omega^2}{2}}^r(\zeta) \quad \mathrm{ when } \quad \mathrm{Im}
(\e^{-\mathrm{i}\frac{\pi}{4}}{\zeta})<0
 \end{cases} \,.  
\end{equation*}

\par\bigskip\noindent
\textbf{\large 2.} {\bf We next discuss $\mathrm{II}(\zeta ,1)$.}
Proposition \ref{Inte.1} gives a representation of $\mathrm{II}(\zeta ,1)$:

\medskip\noindent
{\it Let $\zeta\in \mathbb{C}\setminus \e^{\mathrm{i}\frac{\pi}{4}}\mathbb{R}$. Then
$s_+=\mathrm{Re}(\e^{\mathrm{i}\frac{\pi}{4}} {\zeta}/{|\zeta|})\ne 0$, so $\widetilde{\zeta }(t)\to \infty $ when $t\to +\infty $
and $\mathrm{II}(\zeta ,1)$ is well defined and as in (\ref{II.22}) it
is equal to
$$
\frac{1}{\mathrm{i}}  \int_0^{+\infty}b(\zeta ,t)\e^{\frac{1}{h}\varphi (\zeta
    ,t)}\, {\ensuremath{\mathrm{d}t}}\,, 
$$
where
$$
  b(\zeta ,t)=a(\zeta,\widetilde
    \zeta(t);h)\,\chi(\zeta,\widetilde\zeta(t))\,
    \dfrac{\overline{\widetilde{\zeta}(t)}}{\sqrt{1+t^2}}\,,
$$
$$
\varphi (\zeta ,t)=\mathrm{i}\phi (\zeta ,\widetilde{\zeta }(\zeta ,t)),
$$
and
$\widetilde{\zeta}(t):=\widetilde{\zeta}(\zeta ,
t)$ are given by \eqref{II.20}, \eqref{II.20,1} and \eqref{II.20,2} respectively.
For $0\ne \zeta =|\zeta |\sigma $, $\sigma \in S^1$,
\begin{equation*}
\widetilde\zeta(t) 
= |\zeta|\left(2s_+\, \e^{\mathrm{i}\frac{\pi}{4}}t+\overline{\sigma }g_1(t) \right),
\end{equation*}
 where 
$\sigma=s_+\, \e^{-\mathrm{i}\frac{\pi}{4}}+s_-\, \e^{\mathrm{i}\frac{\pi}{4}}$ with 
$s_\pm=\mathrm{Re}\left(\e^{\pm \mathrm{i}\frac{\pi}{4}}\sigma\right),$ satisfying
 $s_+^2+s_-^2=1$,
and 
\begin{equation*}
 g_1(t):=\dfrac{1}{\sqrt{1+t^2}+t}\,,\ \ t\geq 0\,,
\end{equation*}
satisfies
\begin{equation*}
\partial_t^kg_1(t)={\mathcal  O}\left((1+t)^{-1-k}\right),
\end{equation*}
for every $k\in \mathbb{N}$.
}
We also have as in (\ref{5.nonq3})
\begin{equation*}
|\partial_t\varphi(\zeta ,t)|\asymp |\zeta |^2t,\hbox{ when }t\in
[T_0,+\infty [ \hbox{ and } (\zeta ,\widetilde{\zeta }(\zeta ,t))\in
\mathrm{supp}(\chi ),
\end{equation*}
when $T_0>0$ is large enough, allowing us to define as in (\ref{2.add2}),
\begin{equation*}
 c^{(N)}=\sum_{k=0}^{N}\Big(-\frac{h}{|\zeta|^2\partial_t\varphi_2}\partial_t\Big)^{k} \frac{h}{|\zeta|^2\partial_t\varphi_2} b =: \sum_{k=0}^{N} c_k\,. 
\end{equation*}
Here $\varphi_2(\zeta,t):=\mathrm{i} \phi_2(\zeta,\widetilde{\zeta}(\zeta,t))$, see \eqref{II.30}.

\medskip
Recalling that $\varepsilon \le h^{-\delta }$ we split this region into
two and start with the case of {\bf intermediate values of $\varepsilon $}
(cf.\ (\ref{II.39,6})), i.e.
\begin{equation*}
h^\delta \le \varepsilon \,    \le h^{-\delta }\,.
\end{equation*}
In this case we have Proposition \ref{1nonq3}:

\noindent
{\it 
  Assume that $\chi $ in (\ref{7.ft7})  has
  sufficiently small support.
For $T\ge T_0$, $T\gg \frac{\sqrt{h}}{|\zeta|}=\varepsilon^{-1}$, we have for every $b\in
|\zeta |\mathbb{S}^0([0,+\infty [)$ and $N\in \mathbb{N}$:
 \begin{equation*}
   \int_T^{+\infty} b(t) \e^{\frac{|\zeta|^2}{h}\varphi (\sigma ,t)}\,
   {\ensuremath{\mathrm{d}t}}=-\e^{\frac{|\zeta|^2}{h}\varphi (\sigma ,T)} \left(
     c^{(N)}(T)+\mathcal{O}(1)\left(\frac{h}{(T|\zeta|)^2}\right)^{N+1}\frac{h}{T|\zeta
       |}\right).
 \end{equation*}
 Here $c^{(N)}$ is defined as in (\ref{2.add2}), but with $\varphi_2$
 replaced by $\varphi$,  (\ref{4.add3}) remains valid. Here $\mathbb{S}^0([0,+\infty [)$ is a symbol space given by \eqref{4.add2}, see definition \ref{1add1}. 
}

\noindent
The integral $\int_0^T b(t) 
\e^{\frac{|\zeta|^2}{h}\varphi (\sigma ,t)}\, 
{\ensuremath{\mathrm{d}t}}$ can be computed with standard numerical 
techniques for the values of $\varepsilon$ and $T$ considered here.

We next turn to the {\bf the small $\varepsilon$ limit:}
$$
\varepsilon \le h^\delta\,.
$$
We can perform asymptotic expansions in $\varepsilon$  and get as in
(\ref{1.gt8}), (\ref{2.gt8}), with $M\geq 2L$,
\begin{multline*}
  \frac{1}{h^{\frac12}\varepsilon}\int_0^{+\infty} b(\zeta ,t;h)\, \e^{\frac{1}{h}\varphi (\zeta
    ,t)}\,{\ensuremath{\mathrm{d}t}}\\ ={\mathcal
    O}(\frac1\varepsilon )\left(h^{\frac{N}{2}-1}+\varepsilon ^{2L}\right)+
  \mathrm{an}(\varepsilon ) + \hbox{a linear
    combination of terms}\\
  \sum_{j=0}^{K} h^{\frac{j}{2}+\frac{\widetilde
      K}{2}-n}\varepsilon ^{2\nu}\left(\sum_{\mu =0}^{ j+\widetilde
      K} \varepsilon ^{j+\widetilde K-\mu} \frac{\widetilde
      d_\mu}{\varepsilon }+ \sum_{\mu=-M-1}^{ -1}
    \varepsilon ^{j+\widetilde K} \widetilde d_\mu
    \varepsilon ^{-\mu-1}\ln(\frac1\varepsilon )\right)
\end{multline*}
with 
\begin{equation*}
 K,\widetilde K, N, \nu\in \mathbb{N},\ \nu<L, \ \frac{\widetilde K}{2}-\left[\frac{\widetilde K}{3}\right]<\frac{N}{2}-1,\ n\leq \left[\frac{\widetilde K}{3}\right]\,.
\end{equation*}
Here $\widetilde d_\mu$ also depends on $j,\widetilde K,\nu$ and
$\sigma $.
We also have (\ref{3.gt8}),
\begin{align*}
  \frac{1}{h^{\frac12}\varepsilon}&\int_0^{+\infty}  b(\zeta ,t;h)\,\e^{\frac{1}{h}\varphi (\zeta
    ,t)}\,{\ensuremath{\mathrm{d}t}} ={\mathcal
    O}(\frac1\varepsilon )\left(h^{\frac{N}{2}-1}+\varepsilon ^{2L}\right)
  +\mathrm{an}(\varepsilon ) + \hbox{a linear
    } \nonumber\\
  &\hbox{ combination of terms } \varepsilon ^{2\nu} h^{\frac{\widetilde K}{2}-n} \sum_{j=0}^K
  h^{\frac{j}{2}}\left(p_j(\varepsilon )\frac{1}{\varepsilon }+\varepsilon ^{j+\widetilde
      K}q_j(\varepsilon )\ln(\frac1\varepsilon )\right),
\end{align*}
where $p_j$, $q_j$ are polynomials in $\varepsilon $
still with \eqref{2.gt8} valid. 
Here $\mathrm{an}(\varepsilon)$ is at most $\mathcal{O}(1)$ 
with respect to the two parameters $h$ and $\varepsilon$.

The leading term in $\mathrm{II}(\zeta ,1)$  is given by
(\ref{6.lt1}):
\begin{equation*}\begin{split}
\mathrm{II}(\zeta,1)&=\frac{1}{\mathrm{i}} \int_0^{+\infty}  b(\zeta ,t;h) \, \e^{\frac{1}{h}\varphi (\zeta ,t)} \, \mathrm{d}t\\
 &=\left[h^{\frac12}
 \frac{2\sqrt{2\pi }\, a(0,0) \e^{\mathrm{i}\frac{\pi}{4}}}{(\lambda +\mu +\mathrm{i} 2\rho)^{\frac12}}\times
 \begin{cases}
  \frac{1}{2} \, &\mathrm{if} \quad \mathrm{Im}\big(\e^{-\mathrm{i}\frac{\pi}{4}}\zeta\big)>0\\
-\frac{1}{2} \,  &\mathrm{if} \quad \mathrm{Im}\big(\e^{-\mathrm{i}\frac{\pi}{4}}\zeta\big)<0  
 \end{cases}\right]+R(h,\varepsilon)\,,
   \end{split}
\end{equation*}
where the remainder term $R(h,\varepsilon)$ is ${\mathcal O}(h^{\frac12}\varepsilon)+{\mathcal O}(h)+{\mathcal O}(\varepsilon^2)$.
\par\medskip\noindent
\textbf{\large 3.} {\bf As for $\mathrm{III}(\zeta ,1)$,} we show  in Section
\ref{tt} that
$$
\mathrm{III}(\zeta ,1)={\mathcal O}(h^\infty )\,.
$$

\par\medskip
The main result of this paper is the asymptotic behavior of the integral \eqref{int.1} (same as \eqref{3.le1}, 
\eqref{3.sto3}) in the regime $|\zeta|\in ]0,h^{\frac{1}{2}-\delta}[$ with $\delta>0$. Recalling \eqref{leadI}, \eqref{6.lt1} and \eqref{tt.24}, we obtain

\begin{theo}
 Let $V\subset{\mathbb{R}}^2\simeq{\mathbb{C}}$ be an open neighborhood of $0$,
  $\phi \in C^\infty (V;{\mathbb{R}})$. Assume (\ref{1.le1}),
  (\ref{2.le1}), (\ref{2,5.le1}). For $|\zeta|\in ]0,h^{\frac{1}{2}-\delta}[$ with $\delta\in ]0,\frac{1}{6}[$ satisfying \eqref{1.mr5}, define
  $I(a\chi ,\phi ,\zeta ;h)$ as in \eqref{3.le1} with $a$ is a holomorphic classical symbol, and $\chi \in C_0^\infty (V)$ equal
  to 1 near $0$. Then, the leading term $\big($modulo $\mathcal{O}(h^{\frac{1}{2}+\delta})\big)$ of $ I(a\chi, \phi,\zeta;h)$ is given by: 
\begin{equation*}
 h^{\frac12}\,\frac{2\sqrt{2\pi} \, \e^{\mathrm{i}\frac{\pi}{4}}a(0,0)}{(\lambda+\mu+\mathrm{i}2\rho)^{\frac12}}\,\, \times\begin{cases}
G_{-\frac{\omega^2}{2}}^l\Big(h^{-\frac12}\big(\frac{\mathrm{i}}{f}\big)^{-\frac12}\zeta\Big)+\frac{\mathrm{1}}{2}\quad \mathrm{ when } \quad \mathrm{Im}
(\e^{-\mathrm{i}\frac{\pi}{4}}\zeta)>0\\
{} \quad {} \, {}\\
G_{-\frac{\omega^2}{2}}^r\Big(h^{-\frac12}\big(\frac{\mathrm{i}}{f}\big)^{-\frac12}\zeta\Big)-\frac{\mathrm{1}}{2} \quad \mathrm{ when } \quad \mathrm{Im}
(\e^{-\mathrm{i}\frac{\pi}{4}}{\zeta})<0\,,
 \end{cases}   
\end{equation*}
where $f$ is given by \eqref{4,5.ft10}. 
\end{theo}

\par\smallskip
\begin{remark}
Note that if we further assume that $|\zeta |\ge \exp(-\exp (1/{\mathcal O}(h))$, then we prove that the three terms
$\mathrm{I}(\zeta,1)$, $\mathrm{II}(\zeta,1)$ and  $\mathrm{III}(\zeta,1)$ obtained via Stokes' formula are invariant under rotations modulo an exponentially small error. As a result, the dominant term of the integral \eqref{int.1} is then independent of the choice of parametrization of the deformation contours. 
\end{remark}

\noindent
\textsl{The paper is organized as follows}: In Section 2 we give a direct 
approach to deal with the case where the pole is far from the 
critical point of the phase. In 
Section 3 we propose the decomposition of the integral (\ref{int.1}) 
via Stokes' formula. The first term $\mathrm{I}(\zeta,1)$  of this decomposition is discussed 
in Section 4, the second in Section 5, and we prove that the third 
term is $\mathcal{O}(h^\infty),$ in Section 6. 
We add some concluding remarks in Section 7.
In Appendix A, we briefly review some elements of the theory of wavefront sets.
In Appendix B, we provide additional information on the special function that appears in the asymptotic expression for the first term $\mathrm{I}(\zeta,1)$ of the decomposition of the integral (\ref{int.1}) via Stokes' formula. Appendix C is devoted to the invariance of $\mathrm{I}(\zeta,1)$ under rotations.

\section{Direct approach when \texorpdfstring{$\frac{|\zeta|}{\sqrt{h}}\gg 1$}{}}\label{le}
\setcounter{equation}{0}

In this section we discuss the case that the 
singularity is well separated from the stationary 
point, $|\zeta|\gg \sqrt{h}$, where a standard steepest descent 
approach can be 
applied.

Let $V\subset \mathbb{R}^2\simeq \mathbb{C}$ be an open neighborhood of $0$, let $\phi\in C^\infty(V;\mathbb{R})$ satisfy
\begin{equation}\label{1.le1}
 \phi(0)=0,\qquad \phi'(0)=0\,,
\end{equation}
\begin{equation}\label{2.le1}
 {\rm det}\big(\phi''(0)\big)\not=0\,,
\end{equation}
\begin{equation}\label{2,5.le1}
\phi'(\omega)\not=0,\hbox{ when }  \omega\in V\setminus\{0\}\,. 
\end{equation}
We are interested in the asymptotic behaviour of 
\begin{equation}\label{3.le1}
 I(a\chi, \phi,\zeta;h)=-\frac{1}{\pi}\int_{\mathbb{C}} \frac{1}{\omega-\zeta}\, a(\omega;h)\chi(\omega)\, \e^{\frac{\mathrm{i}}{h}\phi(\omega)}\, L(\mathrm{d}\omega)\,,
\end{equation}
where $a(\omega;h)\in C^\infty(V\times ]0,h_0])$ is a classical symbol of order $h^0$, i.e. such that
\begin{equation}\label{4.le1}
 a(\omega;h)\sim a_0(\omega)+h a_1(\omega)+\ldots\hbox{ in
 }C^\infty(V),\quad h\to 0\,,
\end{equation}
and $\chi(\omega)\in C_0^\infty(V)$ is a fixed cutoff function, equal
to 1 near $\omega=0$.

\medskip\par 
The main contribution to the integral comes from a small neighborhood
of $\{0,\zeta  \}\in \mathbb{C}_\omega $, where $0$ is the critical
point of the exponent and $\zeta $ is the locus of the pole
singularity of the integrand.

\medskip\par
In this section we study the easy case when
\begin{equation}\label{5.le1}
 \varepsilon     :=\frac{|\zeta|}{\sqrt{h}} \gg 1\, .
\end{equation}
In this case it is possible to separate the contributions from $\omega
=0$ and from $\omega =\zeta $. We start with the even easier case when
\begin{equation}\label{6.le1}
 |\zeta|\geq \frac{1}{\mathcal{O}(1)}\,.
\end{equation}
Let $\chi_0\in C_0^\infty(V)$ have support in a small neighborhood of 
$\omega=0$ and be equal to 1 in another such neighborhood of 0. Then $\chi=\chi_0+\chi_1$, where $\chi_1\in C_0^\infty(V)$, $0\not\in \mathrm{supp}(\chi_1)$, $\zeta\not\in\mathrm{supp}(\chi_0)$. Write
\begin{equation}\label{7.le1}
 I(a\chi, \phi,\zeta;h)=I(a\chi_0, \phi,\zeta;h)+I(a\chi_1, \phi,\zeta;h)\, .
\end{equation}
The function
$\omega\longmapsto \frac{1}{\omega-\zeta}\, a(\omega;h)\chi_0(\omega)$
vanishes for $\omega$ close to $\zeta$ and is therefore a classical
symbol of order $h^0$. By stationary phase we get (cf. \cite[Chapter
5]{DiSj99_01})
\begin{equation}\label{8.le1}
 I(a\chi_0, \phi,\zeta;h)=b(\zeta;h)\,,
\end{equation}
where 
\begin{equation}\label{9.le1}
 b\sim h\big(b_1+h b_2+\ldots \big)
\end{equation}
is a classical symbol of order $h$ with
\begin{equation}\label{1.le2}
 b_1(\zeta)=\frac{2}{\zeta} \, a_0(0)  \, 
 \frac{\e^{\mathrm{i}\frac{\pi}{4}\mathrm{sgn}\phi''(0)}}{|\mathrm{det}\phi''(0)|^{\frac{1}{2}}}\,.
\end{equation}

\par 
Near $\mathrm{supp}(\chi_1)$, we have $\phi'(\omega)\not=0$ and it is
straight forward to construct a classical symbol $c(\omega;h)$ of
order $h$ such that
\begin{equation}\label{2.le2}
a(\omega;h)\chi_1(\omega)\e^{\frac{\mathrm{i}}{h}\phi(\omega)}=\partial_{\overline{\omega}}\Big(c(\omega;h) \e^{\frac{\mathrm{i}}{h}\phi(\omega)}\Big)+\mathcal{O}(h^\infty)\,,
\end{equation}
where ${\mathcal O}(h^\infty )$ denotes a remainder term which is ${\mathcal
  O}(h^N)$ uniformly in $\omega $, for every $N>0$.
Here 
\begin{equation}\label{2,5.le2}
c\sim h\big(c_1(\omega)+hc_2(\omega)+\ldots\big)
\end{equation}
with
\begin{equation}\label{3.le2}
 c_1(\omega)=\frac{a_0(\omega)\chi_1(\omega)}{\mathrm{i}\partial_{\overline{\omega}}\phi(\omega)}\, . 
\end{equation}
Moreover $\mathrm{supp}\big(c(\cdot;h)\big)\subset \mathrm{supp}(\chi)$.

\smallskip\par
Since formally $I(a\chi_1,\phi,\zeta;h)=\partial_{\overline{\zeta}}^{-1}\left(a\chi_1\, \e^{\frac{\mathrm{i}}{h}\phi}\right)$, we expect that
\begin{equation}\label{4.le2}
 I(a\chi_1,\phi,\zeta;h)=c(\zeta;h)\, \e^{\frac{\mathrm{i}}{h}\phi(\zeta)}+
 \mathcal{O}(h^\infty)\,.  
\end{equation}
One can prove this by substituting \eqref{2.le2} in the expression for $I(a\chi,\phi,\zeta;h)$ (see \eqref{3.le1}), integrate by parts:
\begin{equation*}
  I(a\chi,\phi,\zeta;h)=\mathcal{O}(h^\infty)+\frac{1}{\pi}\int_{\mathbb{C}} \partial_{\overline{\omega}}\Big(\frac{1}{\omega-\zeta}\Big) c(\omega;h) \, \e^{\frac{\mathrm{i}}{h}\phi(\omega)}\, L(\mathrm{d}\omega)
\end{equation*}
and use that
\begin{equation*}
  \partial_{\overline{\omega}}\left(\frac{1}{\omega-\zeta}\right)  =\pi \delta_\zeta(\omega)\, ,
\end{equation*}
where $\delta _\zeta (\omega )$ is the Dirac measure at $\zeta$. This gives \eqref{4.le2}.

\medskip\par
Combining \eqref{4.le2}, \eqref{8.le1} we get when 
$ |\zeta|\geq \frac{1}{\mathcal{O}(1)}$
\begin{equation}\label{5.le2}
I(a\chi,\phi,\zeta;h)=c(\zeta;h)\, \e^{\frac{\mathrm{i}}{h}\phi(\zeta)}+
b(\zeta;h)+\mathcal{O}(h^\infty)\, ,
\end{equation}
where the symbols $b,c$ are of order $h$ and satisfy \eqref{3.le2}, \eqref{1.le2}.

\medskip\par
We now allow $|\zeta|$ to be small and assume that
\begin{equation}\label{1.le3}
 |\zeta|\ll 1,\qquad \frac{|\zeta|}{\sqrt{h}}=:\varepsilon \gg 1\,.
\end{equation}

\par
We shall use a cutoff
\begin{equation}\label{2.le3}
 \widetilde\chi_0\left(\omega/|\zeta |\right)\,,
\end{equation}
where $\widetilde\chi_0\in C_0^\infty(\mathbb{R}^2)$ is equal to 1 on the disc
$D(0,2)$\footnote{For $r>0,$ $D(0,r)=\{x=(x_1,x_2)\in \mathbb{R}^2, \, x_1^2+x_2^2<r\}$.}, so that
$\widetilde\chi_0\left(\omega/|\zeta |\right)=1$ in the disc
$D(\zeta,2|\zeta |)$. Consider
\begin{equation}\label{3.le3}
J:=-\frac{1}{\pi}\int_{\mathbb{C}} \frac{\big(1-\widetilde\chi_0(\omega /|\zeta |)\big)}{\omega-\zeta}\,\,\, a(\omega;h)\chi(\omega)\,\e^{\frac{\mathrm{i}}{h}\phi(\omega)}\, L(\mathrm{d}\omega)  
\end{equation}
and make repeated integration by parts, using that
\begin{equation}\label{4.le3}
|\phi'(\omega)|\asymp |\omega| \,,\qquad 
\frac{h}{\mathrm{i}}\frac{1}{|\phi'(\omega)|^2}\phi'(\omega)\cdot\nabla_\omega\big(\e^{\frac{\mathrm{i}}{h}\phi(\omega)}\big)=\e^{\frac{\mathrm{i}}{h}\phi(\omega)}\,.
\end{equation}
We write $X \asymp Y$ for $X, Y \in \mathbb{R}$ if $X, Y$ have the same sign (or
vanish) and $X = \mathcal{O}(Y )$ and $Y = \mathcal{O}(X)$. Here we review $\omega\cong (\mathrm{Re}(\omega), \mathrm{Im}(\omega))$ as an element of $\mathbb{R}^2_{\mathrm{Re}(\omega), \mathrm{Im}(\omega)}.$

\noindent
The adjoint of
$\frac{1}{|\phi'(\omega)|^2}\phi'(\omega)\cdot\nabla_\omega$ is equal
to
$$
L:=\displaystyle -\sum_{j=1}^2 \partial_{\omega_j}\circ
\Big(\frac{\partial_{\omega_j}\phi (\omega)}{|\phi(\omega)|^2}\Big)\,, 
$$
where $w=\xi+\mathrm{i}\eta\cong(\xi,\eta)\in \mathbb{R}^2$.
We get 
$$
J=-\frac{1}{\pi}\int \e^{\frac{\mathrm{i}}{h}\phi(\omega)}
\Big(\frac{h}{\mathrm{i}}L\Big)^N\Big(\frac{\big(1-\widetilde\chi_0(\omega/|\zeta
  |)\big)}{\omega-\zeta}\, a(\omega;h)\chi(\omega)\Big)\, L(\mathrm{d}\omega)\,. 
$$
Here $|\zeta-\omega|\geq \frac{|\omega|}{2}$ on
$\mathrm{supp}\big(1-\widetilde\chi_0(\omega/|\zeta |)\big)$ and we see that
$$
J=\mathcal{O}(1)\int_{2|\zeta|\leq |\omega|\leq \frac{1}{\mathcal{O}(1)}} \,\, \frac{h^N}{|\omega|^{2N+1}}\, \, L(\mathrm{d}\omega)\,.
$$
Here $\displaystyle \frac{h}{|\omega|^2}\leq \mathcal{O}(1) \frac{h}{|\zeta|^2}=\mathcal{O}(1)\frac{1}{\varepsilon ^2}$ for $\omega$ in the domain of integration so 
\begin{equation}\label{5.le3}
J=\mathcal{O}(1)\varepsilon ^{-2N}\int_{2|\zeta |\leq |\omega|\leq \frac{1}{\mathcal{O}(1)}} \, \frac{1}{|\omega|}\, L(d\omega)=\mathcal{O}(1)\varepsilon ^{-2N}\,.  
\end{equation}
We conclude that for every $N\in \mathbb{N}$,
\begin{equation}\label{1.le4}
  I(a\chi,\phi,\zeta;h)=I(a\widetilde \chi_0(\cdot /|\zeta |),\phi,\zeta;h)+\mathcal{O}(\varepsilon ^{-2N})\,.
\end{equation}
In order to study
$(a\widetilde \chi_0(\cdot /|\zeta |),\phi,\zeta;h)$, we put
$\zeta=\lambda\widetilde\zeta$, where $\lambda \asymp |\zeta |$ is an
extra parameter. Make the change of variables
$\omega=\lambda \widetilde\omega$. Write
\begin{equation}\label{2.le4}
 \frac{1}{h}\phi(\omega)=\frac{\lambda^2}{h}\frac{1}{\lambda^2}\phi(\lambda\widetilde\omega)=\frac{1}{\widetilde h}\widetilde\phi(\widetilde\omega)\,,
\end{equation}
where
\begin{equation}\label{3.le4}
\widetilde h=\frac{h}{\lambda^2}\asymp\frac{1}{\varepsilon ^2}\ll 1, \quad 
\widetilde\phi(\widetilde\omega)=\frac{1}{\lambda^2}\phi(\lambda\widetilde \omega)\,,
\end{equation}
so the remainder in \eqref{1.le4} is $\mathcal{O}\big(\,{\widetilde {h}}^N\big)$.

\medskip\par
Clearly $\widetilde\phi(0)=0$, $\widetilde\phi'(0)=0$, $\widetilde\phi''(0)=\phi''(0)$ and we get
\begin{equation}\label{4.le4}
  I(a\widetilde \chi_0(\cdot /|\zeta |),\phi,\zeta;h)=-\frac{\lambda}{\pi}\int \frac{1}{\widetilde
    \omega-\widetilde \zeta}\, \widetilde a
  (\widetilde\omega)\widetilde\chi_0((\lambda /|\zeta |)\widetilde\omega)\, \e^{\frac{\mathrm{i}}{\tilde h}\tilde\phi(\widetilde\omega)} L(\mathrm{d}\widetilde\omega)\,,
\end{equation}
where
\begin{equation}\label{5.le4}
\widetilde a(\widetilde\omega;\widetilde h)=a(\omega;h)=a(\lambda\widetilde\omega;\lambda^2\widetilde h)\,.  
\end{equation}
Since $\lambda\leq 1$, $\widetilde h\geq h$, $\widetilde a$ is a holomorphic
classical symbol with $\widetilde h$ as the new semiclassical
parameter, uniformly with respect to $\lambda $. Moreover
$|\widetilde \zeta|\asymp 1$ so $\widetilde \zeta$ is
well separated from $\widetilde\omega=0$. We can therefore apply the asymptotic formula \eqref{5.le2} to
$ I(a\widetilde \chi_{0,\lambda},\phi,\zeta;h)$ in \eqref{4.le4} and
get for $ I(a\chi,\phi,\zeta;h)$ in \eqref{1.le4}:
\begin{equation}\label{6.le4}
 I(a\chi,\phi,\zeta;h)=\mathcal{O}\big(\widetilde h^\infty\big)+\lambda \Big(\widetilde c(\widetilde\zeta;\widetilde h)\, \e^{{\mathrm{i}}\tilde\phi(\widetilde\omega)/\widetilde h}+\widetilde b(\widetilde\zeta;\widetilde h)+\mathcal{O}\big(\widetilde h^\infty\big)\Big)\,,
\end{equation}
where $\widetilde c(\widetilde\zeta;\widetilde h), \widetilde b(\widetilde\zeta;\widetilde h)$ are classical symbols of order $\widetilde h$ satisfying the analogues of \eqref{2,5.le2}, \eqref{9.le1}:
\begin{equation}\label{7.le4}
\left \{
 \begin{array}{cc}
    \widetilde c(\widetilde\zeta;\widetilde h)\sim 
    \widetilde h\big(\widetilde c_1(\widetilde\zeta)+\widetilde h\widetilde c_2(\widetilde\zeta)+\ldots\,\big)\\
     {}  {} \\
     \widetilde b(\widetilde\zeta;\widetilde h) \sim 
    \widetilde h\big(\widetilde b_1(\widetilde\zeta)+\widetilde h\widetilde b_2(\widetilde\zeta)+\ldots\,\big)
\end{array}
\right.
\end{equation}
with
\begin{equation}\label{1.le5}
\widetilde b_1(\widetilde\zeta)=\frac{2}{\widetilde\zeta}\,\, \widetilde a_0(0)\,\,  
 \dfrac{\e^{\mathrm{i}\frac{\pi}{4}\mathrm{sgn}\phi''(0)}}{|\mathrm{det}\phi''(0)|^{\frac{1}{2}}} \,,
\end{equation}

\begin{equation}\label{2.le5}
 \widetilde c_1(\widetilde\zeta)=\dfrac{\widetilde a_0(\widetilde\zeta)}{\mathrm{i}\partial_{\overline{\tilde\zeta}}\tilde\phi(\tilde\zeta)}\, .  
\end{equation}

On the formal asymptotic level we have  
\begin{equation}\label{3.le5}
  \lambda\,\widetilde c(\widetilde\zeta;\widetilde h)=c(\zeta;h)\,,\qquad 
  \lambda\,\widetilde b(\widetilde\zeta;\widetilde h)=b(\zeta;h)\,.
\end{equation}
If $\lambda\geq \dfrac{1}{\mathcal{O}(1)} h^{\frac{1}{2}-\delta}$,
or equivalently 
\begin{equation}\label{4.le5}
\widetilde h\leq \mathcal{O}(1) h^{\frac{\delta}{2}}\,,
\end{equation}
for some $\delta >0$, then the remainders $\mathcal{O}\big((\widetilde h)^\infty\big)$ and $\lambda\mathcal{O}\big((\widetilde h)^\infty\big)$ are $\mathcal{O}(h^\infty)$.

\noindent
The following proposition summarizes the results of this section
\begin{prop}\label{1le5}
  Let $V\subset{\mathbb{R}}^2\simeq{\mathbb{C}}$ be an open neighborhood of $0$,
  $\phi \in C^\infty (V;{\mathbb{R}})$. Assume (\ref{1.le1}),
  (\ref{2.le1}), (\ref{2,5.le1}). For $\zeta ={\mathcal O}(1)$, define
  $I(a\chi ,\phi ,\zeta ;h)$ as in (\ref{3.le1}) with the classical
  symbol $a$ as in (\ref{4.le1}) and $\chi \in C_0^\infty (V)$ equal
  to 1 near $0$.

  In any open set where $|\zeta |\asymp 1$, we have
  (\ref{5.le2}), where $c$, $b$ are classical symbols of order $h$,
  satisfying (\ref{9.le1}), (\ref{1.le2}), (\ref{2,5.le2}), 
  (\ref{3.le2}).
  
  Assume now (\ref{5.le1}), let the parameter $\lambda $ be
  $\asymp |\zeta |$ and put
  $\widetilde{h}=h/\lambda ^2\asymp\varepsilon ^{-2}$, so that
  $\widetilde{h}\ll 1$. Define the rescaled variable
  $\widetilde{\omega }$ by $\omega =\lambda \widetilde{\omega }$ and
  the rescaled functions $\widetilde{\phi }(\widetilde{\omega })$,
  $\widetilde{a}(\widetilde{\omega };\widetilde{h})$ by (\ref{3.le4}),
  (\ref{5.le4}). Put $\zeta =\lambda \widetilde{\zeta }$. Then we have (\ref{6.le4}), where
  $\widetilde c(\widetilde\zeta;\widetilde h), \widetilde
  b(\widetilde\zeta;\widetilde h)$ are classical symbols of order
  $\widetilde h$, satisfying (\ref{7.le4}), (\ref{1.le5}),
  (\ref{2.le5}).  If
  $ \widetilde h\leq \mathcal{O}(1) h^{\frac{\delta}{2}}$, for some
  $\delta >0$, then the remainders
  $\mathcal{O}\big(\widetilde h^\infty\big)$ and
  $\lambda\mathcal{O}\big(\widetilde h^\infty\big)$ are
  $\mathcal{O}(h^\infty)$.  (\ref{6.le4}) coincides with (\ref{5.le2})
  in the sense of formal asymptotic expansions and is therefore
  independent of the choice of $\lambda \asymp |\zeta |$.
\end{prop}

\section{Decomposition by Stokes' formula}\label{Sto}
\setcounter{equation}{0}

In this section, we apply Stokes' formula to the 
integral on steepest descent contours in $\mathbb{C}^{2}$. This leads 
to a decomposition of the studied integral into three integrals.

We recall a contour deformation version of Stokes' formula (cf.
\cite[lemme 12.2]{Sj82_01}). Let $0<k<n$ be integers. Let $f:\, I\times
W\ni (t,x)\mapsto f(t,x)\in \mathbb{R}^n$ be a smooth map, where
$I\Subset \mathbb{R}$, $W\Subset \mathbb{R}^k$ are open and connected. We
introduce the integration contour
\begin{equation}\label{1.sto1}
W\ni x\mapsto f(t,x)=:f_t(x)\in \mathbb{R}^n. 
\end{equation}
When $f_t$ is injective with injective differential we can view
$\Gamma_t:=f_t(W)$ as a submanifold of $\mathbb{R}^n$, at least locally.

\par Let $\Omega $ be a smooth differential $k$-form on some open
subset of $\mathbb{R}^n$ such that
\begin{equation}\label{2.sto1}
\pi _x(f^{-1}(\mathrm{supp\,}\Omega ))\Subset W\,,
\end{equation}
where $\pi _x:I\times W\to W$ is the natural projection. Let
\begin{equation}\label{3.sto1}
  \nu _t=\partial _tf(t,x)\in T_{f_t(x)}\mathbb{R}^n\,.
  \end{equation}
  \begin{lemma}\label{1sto1}
    $\displaystyle \int_{\Gamma_t}\hskip-3pt\Omega$ is a smooth function of $t$ and we have
    \begin{equation}\label{4.sto1}
\partial _t\int_{\Gamma _t}\hskip-3pt\Omega =\int_{\Gamma _t}\hskip-3pt\nu _t\rfloor \mathrm{d}\Omega\,.
    \end{equation}
Here $\nu_t\rfloor \mathrm{d}\Omega$ denotes the contraction of $\mathrm{d}\Omega$ with $\nu_t$.   
  \end{lemma}
  
  \noindent
 \textbf{Proof.} 
With the differential form $\widetilde{\Omega }=f^*\Omega $ on
$I\times W$ and
$$
\widetilde{\Gamma }_t=\{ (t,x);\, x\in W \}
$$
(\ref{4.sto1}) amounts to
\begin{equation}\label{5.sto1}
\partial _t\int_{\widetilde{\Gamma } _t}\widetilde{\Omega}
=\int_{\widetilde{\Gamma} _t}\frac{\partial }{\partial t}\rfloor
\mathrm{d}\widetilde{\Omega } \,,
\end{equation}
which follows from a simple calculation.\hfill$\square$

 \begin{remark}\label{2sto1}
Assume that $f$ is injective with injective differential, so that
$f$ is a diffeomorphism (at least locally): $I\times W\to \Gamma $,
where $\Gamma $ is a smooth submanifold of dimension $k+1$. Then there is a
smooth vector field $\nu $ on $\Gamma $, such that $\nu (x)=\nu _t(x)$
when $x\in \Gamma _t$ and $\Gamma _t$ are $k$-dimensional smooth
submanifolds of $\Gamma $. Clearly,
\begin{equation}\label{1.sto2}
\Gamma_s=\big(\exp (s-t)\nu\big) \big(\Gamma_t\big)\,,\quad s,t\in I\,.
\end{equation}
It follows that
\begin{equation}\label{2.sto2}
\int_{\Gamma_s}\Omega =\int_{\Gamma_t}(\exp (s-t)\nu )^*\Omega\,,
\end{equation}
and taking the derivative with respect to $s$ at $s=t$,
\begin{equation}\label{3.sto2}
\partial_t\int_{\Gamma _t}\Omega =\int_{\Gamma_t}{\mathcal  L}_\nu (\Omega )\,,
\end{equation}
where ${\mathcal  L}_\nu $ denotes Lie derivative with respect to $\nu
$. Using Cartan's identity, we recover (\ref{4.sto1}):
\begin{equation}\label{4.sto2}
\int_{\Gamma_t}{\mathcal  L}_\nu (\Omega )=\int_{\Gamma_t}\nu \rfloor
\mathrm{d}\Omega +\int_{\Gamma_t}\mathrm{d}(\nu \rfloor \Omega )\,.
\end{equation}
Here ${{\nu \rfloor \Omega }_\vert}_{\Gamma_t}$ has compact support
so the integral of its differential, i.e. the last term in
(\ref{4.sto2}), vanishes.
  \end{remark}

  \begin{remark}\label{1sto2}
We make the same geometric assumptions as in Remark \ref{2sto1} and
relax the smoothness assumption on $\Omega $: Assume now that $\Omega
$ is a $k$-form with distribution coefficients satisfying
(\ref{2.sto1}) and
\begin{equation}\label{5.sto2}
\mathrm{WF}(\Omega )\cap N^*(\Gamma _t)=\emptyset ,\ t\in I\,,
\end{equation}
where $N^*(\Gamma _t)$ denotes the conormal bundle of $\Gamma
_t\subset \mathbb{R}^n$. Then ${{\Omega }_\vert}_{\Gamma _t}$ is a well
defined $k$ form on $\Gamma _t$ with distribution coefficients and
$$I\ni t\mapsto {{\Omega }_\vert}_{\Gamma _t}$$ is smooth with values in
the distributional $k$-forms in the sense that $      I \ni t\mapsto
\int_{\Gamma_t}\phi \Omega $ is smooth for all $\phi \in C^\infty
(\mathbb{R}^n)$ and (\ref{4.sto1}) holds. See Appendix \ref{wf} for more
details and a review of wavefront sets.
\end{remark}

\par Let $\phi (\xi ,\eta )$ be analytic near $(0,0)\in \mathbb{R}^2$
  and let $a(\xi ,\eta ;h)$ be a holomorphic classical symbol defined
  near the same point. Let $\chi \in C_0^\infty (\mathbb{R}^2)$ be equal
  to 1 near $(0,0)$ and with support contained in the domains of
  definition of $\phi$ and $a$. We will continue the study of the integral 
  \begin{equation}\label{1.sto3}
I(a\chi ,\phi ,\zeta ;h)=-\frac{1}{\pi}\iintJ_{\mathbb{R}^2}\frac{1}{
(\xi +\mathrm{i}\eta )-\zeta}\,a(\xi ,\eta ;h)\chi (\xi ,\eta )\,\e^{\frac{\mathrm{i}}{h}\phi (\xi
,\eta)}\,\mathrm{d}\xi \mathrm{d}\eta 
  \end{equation}
  for $\zeta \in \mathbb{C}$.
  For this we will use the method of polarization and start by
  identifying $\mathbb{R}^2_{\xi ,\eta }$ with the antidiagonal
  \begin{equation}\label{2.sto3}
    \mathrm{adiag\,}(\mathbb{C}_{\omega ,\tilde{\omega }}^2)=
    \{ (\omega ,\widetilde{\omega })\in \mathbb{C}^2_{\omega
      ,\tilde{\omega }};\,\,\, \widetilde{\omega }=\overline{\omega } \}\,.
  \end{equation}
  Here we put
  \begin{equation}\label{2,5.sto3}
  \omega =\xi +\mathrm{i}\eta , \hbox{ and also }
  \widetilde{\omega}=\xi -\mathrm{i}\eta\,. 
  \end{equation} 
Since
$$
\mathrm{d}\xi \wedge \mathrm{d}\eta = \frac{\mathrm{d}\widetilde
  {\omega}\wedge \mathrm{d}\omega}{2\mathrm{i}}, 
$$
the integral in (\ref{1.sto3}) can be written
\begin{equation}\label{3.sto3}
I(a\chi ,\phi ,\zeta ;h)=-\frac{1}{\pi}\iintJ_{\mathrm{adiag\,}({\mathbb{C}}^2)}\frac{ \e^{\frac{\mathrm{i}}{h}\phi (\omega ,\widetilde{\omega })}}{\omega-\zeta} \, a(\omega ,\widetilde{\omega };h)\chi (\omega
,\widetilde{\omega }) \, \frac{\mathrm{d}\widetilde{\omega }\wedge \mathrm{d}\omega}{2\mathrm{i}}\,, 
\end{equation}
abusing slightly the notation, letting $a$, $\chi$, $\phi $ denote
``the same'' functions on $\mathrm{adiag\,}(\mathbb{C}^2)$ using
(\ref{2,5.sto3}) and the inverse relations
\begin{equation}\label{4.sto3}
\xi =\frac{\omega +\widetilde{\omega}}{2},\quad \eta =\frac{\omega
  -\widetilde{\omega }}{2\mathrm{i}}\,.
\end{equation}

\par Later we will assume that $\phi (\xi ,\eta )$ is real valued on the
real domain with $(0,0)$ as a non-degenerate critical point and such
that this is the only critical point of $\phi $ on the real domain.
One would like to apply the method of stationary phase, but the
singularity at $\omega =\zeta $ poses an obstacle for that, especially
when $|\zeta |\le {\mathcal  O}(\sqrt{h})$. (When $|\zeta |\gg \sqrt{h}$,
we separate the contributions from
the critical point at $\omega =0$ and from the singularity at $\omega
=\zeta $. See Section \ref{le}.)

\par
The idea is then to deform the contour $\mathrm{adiag}({\mathbb{C}}^2_{\omega ,\tilde{\omega }})$ to a new contour in ${\mathbb{C}}^2_{\omega ,\tilde{\omega }}$ which is a Cartesian product
$\gamma _1\times \gamma _2$, where $\gamma _1\subset {\mathbb{C}}_\omega $
and $\gamma _2\subset {\mathbb{C}}_{\tilde{\omega }}$ are
curves. Stokes' formula will then produce several terms to be
analyzed, in particular the integral along $\gamma _1\times \gamma _2$ which
reduces to 1D integrals:
\begin{equation}\label{1.sto4}
\iintJ_{\gamma _1\times \gamma _2}(...)\mathrm{d}\widetilde{\omega} \wedge
\mathrm{d}\omega =\int_{\gamma _1}\left(\int_{\gamma _2}(...)
\mathrm{d}\widetilde{\omega } \right) \mathrm{d}\omega\,.
\end{equation}
  
\par
We will apply Lemma \ref{1sto1} with 
\begin{equation}\label{Sto.2}
 \Omega=-\frac{1}{\pi}\frac{1}{\omega-\zeta} \,\,a(\omega,\widetilde \omega;h)\,\chi(\omega,\widetilde\omega)\, \e^{\frac{\mathrm{i}}{h}\phi(\omega,\widetilde \omega)}\, \frac{\mathrm{d}\widetilde\omega\wedge \mathrm{d}\omega}{2\mathrm{i}}\,,
\end{equation}
where $\chi\in C^\infty_0\big(\mathbb{C}^2;[0,1]\big)$ is a
smooth cutoff function equal to $1$ near (0,0) with support inside the
domains where $a$, $\phi $ are holomorphic.

\par
We define a smooth family of contours $\Gamma_\theta$, $\theta \in [0,1]$, of real dimension 2 in ${\mathbb{C}}^2_{\xi ,\eta }$ or equivalently in $\mathbb{C}^2_{\omega
  ,\tilde{\omega }}$. In the $\mathbb{C}^{2}_{\xi ,\eta }$ representation it is
given by
$$
\mathbb{R}^2\ni (t,s)\mapsto (\xi (\theta ,t,s),\eta (\theta ,t,s))\in
\mathbb{C}^2,$$
where
\begin{equation}\label{Sto.3}
\begin{cases}
 \xi(\theta;t,s) = t+\mathrm{i}\theta t=(1+\mathrm{i}\theta) t\\
 \eta(\theta;t,s) = s-\mathrm{i}\theta s=(1-\mathrm{i}\theta)s\\
\end{cases},\quad (t,s)\in \mathbb{R}^2\,.
\end{equation}
In the ${\mathbb{C}}_{\omega ,\widetilde{\omega }}$ representation it is
given by
\begin{equation}\label{Sto.5}
\mathbb{R}^2\ni(t,s)\,\mapsto f(\theta;t,s)=\big(\omega(\theta;t,s)\,,\,\widetilde\omega(\theta;t,s)\big)\in\mathbb{C}^2_{\omega,\tilde\omega},
\end{equation}
where
\begin{equation}\label{Sto.4}
\begin{cases}
\omega(\theta;t,s):&\hskip-9pt= \xi(\theta;t,s)+\mathrm{i}\eta(\theta;t,s)\\
     & \hskip-9pt= (1+\mathrm{i}\theta) t+\mathrm{i}(1-\mathrm{i}\theta)s\\
     & \hskip-9pt= (t+\mathrm{i} s)+\mathrm{i}\theta(t-\mathrm{i} s)\\
\widetilde\omega(\theta;t,s):&\hskip-9pt= \xi(\theta;t,s)-\mathrm{i}\eta(\theta;t,s)\\
      &\hskip-9pt=(1+\mathrm{i}\theta) t-\mathrm{i}(1-\mathrm{i}\theta)s\\
      & \hskip-9pt= (t-\mathrm{i} s)+\mathrm{i}\theta(t+\mathrm{i} s)\\
\end{cases}\,,\quad (t,s)\in \mathbb{R}^2,
\end{equation}
hence
\begin{equation}\label{Sto.6}
f(\theta;t,s)=\big((1+\mathrm{i}\theta) t+\mathrm{i}(1-\mathrm{i}\theta)s, (1+\mathrm{i}\theta) t-\mathrm{i}(1-\mathrm{i}\theta)s\big)\in\mathbb{C}^2_{\omega,\tilde\omega}\,.
\end{equation}
In the $\mathbb{C}^2_{\xi ,\eta}$ representation, we have $\Gamma
_\theta =(1+\mathrm{i}\theta )\mathbb{R}\times (1-\mathrm{i}\theta )\mathbb{R}$, especially $\Gamma
_0=\mathbb{R}^2$.
In the $\mathbb{C}^2_{\omega ,\tilde{\omega }}$ representation, we
get $\Gamma _0=\mathrm{adiag\,}(\mathbb{C}^2)$, and for $\theta =1 $, we
have
$$
\begin{cases}
\omega (1,t,s)&=(1+\mathrm{i})(t+s)=\sqrt{2}\e^{\mathrm{i}\frac{\pi}{4}}(t+s)\\
\widetilde{\omega }(1,t,s)&=(1+\mathrm{i})(t-s)=\sqrt{2}\e^{\mathrm{i}\frac{\pi}{4}}(t-s)  
\end{cases}\,,
$$

\par\medskip\noindent
so $\Gamma _1=(1+\mathrm{i})\mathbb{R}\times (1+\mathrm{i})\mathbb{R}=\e^{\mathrm{i}\frac{\pi}{4}}\mathbb{R}\times \e^{\mathrm{i}\frac{\pi}{4}}\mathbb{R}$ is a
Cartesian product. 

\par\medskip
Identifying $\mathbb{C}^2_{\omega ,\tilde{\omega}}$ with the
underlying real space $\mathbb{R}^4_{\mathrm{Re}\, \omega ,\mathrm{Im}\, \omega ,
\mathrm{Re}\, \tilde{\omega },\mathrm{Im}\, \tilde{\omega }}$, we see that
$\mathrm{WF\,}(\Omega )$ is contained in the conormal bundle of the
set $\{(\omega ,\widetilde{\omega });\,\,\, \omega =\zeta\}$ viewed as
a subset of $\mathbb{R}^4$. We can therefore define the restriction
${{\Omega}_\vert}_{\Gamma_\theta}$ for a given fixed $\theta \in
[0,1]$ if the conormal bundle of $\Gamma_\theta$ has zero
intersection with the conormal bundle of $\{(\omega ,\widetilde{\omega
});\,\,\, \omega =\zeta\}$. Since our manifolds are affine spaces of the
same dimension 2, this is equivalent to the property that they intersect at a single point or still that
$$
\Gamma _\theta \ni (\omega ,\widetilde{\omega } )\longmapsto \omega \in {\mathbb 
  C} \hbox{ is injective}\,.
$$
In view of (\ref{Sto.4}) this is equivalent to the implication
$$
\begin{cases}
t+\mathrm{i} s+\mathrm{i}\theta (t-\mathrm{i} s)=0\\
  t,s\in \mathbb{R}
\end{cases}
\Longrightarrow \,\,\,\, t=s=0\,,
$$
which, by separating real and imaginary parts, amounts to
$$
\det \begin{pmatrix}1 &\theta \\ \theta & 1\end{pmatrix}\ne 0\,,
$$
i.e.\  $\theta ^2\ne 1$. In conclusion,\begin{equation}\label{1.sto5}
{{\Omega }_\vert}_{\Gamma
  _\theta }\hbox{ is well defined when }0\le \theta <1.
\end{equation} 

\par
The deformation vector field in (\ref{3.sto1}),
$\nu(\theta,\cdot):=\partial_\theta f\big(\theta,
\omega(\theta;\cdot),\widetilde \omega(\theta,\cdot)\big)$, is given
along the contour $\Gamma_\theta$ by
\begin{equation}\label{2.sto5}\begin{split}
  \nu (\theta )
  &
  =\big(\partial_{\theta}\omega\big)\cdot \partial_{\omega}
  +\big(\partial_{\theta}\widetilde\omega\big)\cdot \partial_{\widetilde\omega}
  +\big(\partial_{\theta}\overline{\omega}\big)\cdot \partial_{\overline\omega}
  +\big(\partial_{\theta}\overline{\widetilde\omega}\big)\cdot
  \partial_{\overline{\tilde\omega}}\\ &=:
  \nu _\omega (\theta )\cdot \partial _\omega +
  \nu _{\tilde{\omega }} (\theta )\cdot \partial _{\tilde{\omega }}+
\nu _{\overline{\omega } } (\theta )\cdot \partial _{\overline{\omega } }+
  \nu _{\overline{\tilde{\omega }}} (\theta )\cdot \partial
  _{\overline{\tilde{\omega }}} \,,
\end{split}
\end{equation}
and more explicitly by
\begin{equation}\label{Sto.7}
\nu(\theta,\cdot) =\mathrm{i}(t-\mathrm{i} s)\cdot \partial_{\omega}+\mathrm{i}(t+\mathrm{i} s)\cdot \partial_{\tilde\omega}
-\mathrm{i}(t+\mathrm{i} s)\cdot \partial_{\overline\omega}-\mathrm{i}(t-\mathrm{i} s)\cdot
\partial_{\overline{\tilde\omega}} \,.
\end{equation}

\par
Now, fix $\zeta\in \mathbb{C}$ and apply Lemma \ref{1sto1} to obtain
for $0\le\theta <1$,
\begin{equation}\label{Sto.8}
\frac{\mathrm{d}I}{\mathrm{d}\theta}
 =-\frac{1}{\pi}\iintJ_{\Gamma_\theta}\nu(\theta ;\omega,\widetilde\omega)\mathord{\rfloor}  \mathrm{d}_{\omega,\tilde\omega}\left(\frac{1}{\omega-\zeta} \, a(\omega,\widetilde \omega;h)\,\chi(\omega,\widetilde\omega)\, \e^{\frac{\mathrm{i}}{h}\phi(\omega,\widetilde \omega)}\, \frac{\mathrm{d}\widetilde\omega\wedge \mathrm{d}\omega}{2\mathrm{i}}\right),
\end{equation}
where $I(\theta):=I(\Gamma_{\theta};\zeta)$ is given by
\begin{equation}\label{Sto.9}
\hskip-7ptI(\theta):=I(\Gamma_{\theta};\zeta)=-\frac{1}{\pi}
\iintJ_{\Gamma_{\theta}}
\frac{1}{\omega-\zeta} \, a(\omega,\widetilde
w;h)\,\chi(\omega,\widetilde\omega)\,
\e^{\frac{\mathrm{i}}{h}\phi(\omega,\widetilde \omega)
}\, \frac{\mathrm{d}\widetilde\omega\wedge \mathrm{d}\omega}{2\mathrm{i}}\,.
\end{equation}
Further, in the sense of differential forms with coefficients that
are distributions on $\mathbb{C}^2_{\omega ,\tilde{\omega }}$ with
singular support contained in $\{(\omega ,\widetilde{\omega });\,\,\,
\omega =\zeta  \}$,

\begin{equation}\label{Sto.11}\begin{split}
\mathrm{d}_{\omega,\tilde\omega}&  \left(\frac{1}
{\omega-\zeta}\, a(\omega,\widetilde \omega;h)\,\chi(\omega,\widetilde\omega)\, \e^{\frac{\mathrm{i}}{h}\phi(\omega,\widetilde \omega)} \frac{\mathrm{d}\widetilde \omega\wedge \mathrm{d}\omega}{2\mathrm{i}}\right)\\
  =&\frac{\partial}{\partial\overline{\omega}}\left(
  \frac{1}{\omega-\zeta} \, a(\omega,\widetilde \omega;h)\,\chi(\omega,\widetilde\omega)\, \e^{\frac{\mathrm{i}}{h}\phi(\omega,\widetilde \omega) }\right) \frac{\mathrm{d}\overline{\omega}\wedge \mathrm{d}\widetilde\omega\wedge \mathrm{d} \omega}{2\mathrm{i}}\\
     &
     +\frac{\partial}{\partial\overline{\widetilde\omega}}\left(
     \frac{1}{\omega-\zeta}\,
       a(\omega,\widetilde \omega;h)\,\chi(\omega,\widetilde\omega)\,
       \e^{\frac{\mathrm{i}}{h}\phi(\omega,\widetilde \omega)
       }\right)\,\frac{\mathrm{d}\overline{\widetilde\omega}\wedge \mathrm{d}\widetilde
     \omega\wedge \mathrm{d}\omega}{2\mathrm{i}}\,,\\
     =&\pi \delta_{\zeta}(\omega) \, a(\omega,\widetilde \omega;h)\,\chi(\omega,\widetilde\omega)\, \e^{\frac{\mathrm{i}}{h}\phi(\omega,\widetilde \omega) }\, \frac{\mathrm{d}\overline{\omega}\wedge \mathrm{d}\widetilde\omega\wedge \mathrm{d}\omega}{2\mathrm{i}}\\
    &+\frac{1}{\omega-\zeta} \, a(\omega,\widetilde \omega;h)\,\frac{\partial\chi}{\partial\overline{\omega}}(\omega,\widetilde\omega)\, \e^{\frac{\mathrm{i}}{h}\phi(\omega,\widetilde \omega) }\, \frac{\mathrm{d}\overline{\omega}\wedge \mathrm{d}\widetilde\omega\wedge \mathrm{d}\omega}{2\mathrm{i}}\\
    &+ \frac{1}{\omega-\zeta} \, a(\omega,\widetilde
    w;h)\,\frac{\partial \chi}{\partial
      \overline{\widetilde\omega}}(\omega,\widetilde\omega)\,
    \e^{\frac{\mathrm{i}}{h}\phi(\omega,\widetilde \omega)
    }\, \frac{\mathrm{d}\overline{\widetilde\omega}\wedge \mathrm{d}\widetilde\omega\wedge \mathrm{d}\omega}{2\mathrm{i}}\,.
\end{split}
\end{equation}
Here we use that $a$ and $\phi$ are holomorphic.

\par
Equality (\ref{Sto.11}) can be written more briefly as
\begin{align}\label{Sto.12}
 \mathrm{d}_{\omega,\tilde\omega} & \Big(\frac{1}{\omega-\zeta}\, a(\omega,\widetilde \omega;h)\,\chi(\omega,\widetilde\omega)\, \e^{\frac{\mathrm{i}}{h}\phi(\omega,\widetilde \omega)} \frac{\mathrm{d}\widetilde\omega\wedge \mathrm{d}\omega}{2\mathrm{i}}\Big)\nonumber\\
 =&\pi \delta_{\zeta}(\omega)\, a(\omega,\widetilde \omega;h)\,\chi(\omega,\widetilde\omega)\, \e^{\frac{\mathrm{i}}{h}\phi(\omega,\widetilde \omega) }\, \frac{\mathrm{d}\overline{\omega}\wedge \mathrm{d}\widetilde\omega\wedge \mathrm{d}\omega}{2\mathrm{i}}\\
&+\frac{1}{\omega-\zeta} \, a(\omega,\widetilde \omega;h)\, \e^{\frac{\mathrm{i}}{h}\phi(\omega,\widetilde \omega) }\,\, \mathrm{d}_{\omega,\widetilde\omega}\Big( \chi(\omega,\widetilde\omega)\, \frac{\mathrm{d}\widetilde\omega\wedge \mathrm{d}\omega}{2\mathrm{i}}\Big)
\,.\nonumber
\end{align}

\par
Pointwise in $\mathbb{C}^2_{\omega ,\tilde{\omega }}$ we have that
$$
\partial _\omega ,\, \partial _{\overline{\omega }},\, \partial
_{\tilde{\omega }},\, \partial _{\overline{\tilde{\omega
    }}}\quad \hbox{ and }\quad
\mathrm{d}\omega ,\, \mathrm{d}{\overline{\omega }},\, \mathrm{d}{\widetilde{\omega }},\, \mathrm{d}{\overline{\widetilde{\omega
    }}}
$$
are dual bases in
$$
\mathbb{C}\otimes T\mathbb{C}^2_{\omega ,\tilde{\omega }}\hbox{ and }
\mathbb{C}\otimes T^*\mathbb{C}^2_{\omega ,\tilde{\omega }}\hbox{ respectively.}
$$
This gives
\begin{equation}\label{Sto.13}
  \begin{split}
    & {\partial_ \omega}\mathord{\rfloor} \mathrm{d}\omega=1,\quad {\partial
      _ \omega}\mathord{\rfloor} \mathrm{d}\overline{\omega}=0, \quad
    \partial _\omega \mathord{\rfloor} \mathrm{d}\widetilde{\omega }=0, \quad
    \partial _\omega \mathord{\rfloor} \mathrm{d}\overline{\widetilde{\omega
      }}=0\,,
    \\
    & \partial_{ \overline{\omega}}\mathord{\rfloor} \mathrm{d}\omega=0,
    \quad {\partial _{\overline{\omega}}}\mathord{\rfloor}
    \mathrm{d}\overline{\omega}=1, \quad \partial _{\overline{\omega}}
    \mathord{\rfloor} \mathrm{d}\widetilde{\omega }=0, \quad \partial
    _{\overline{\omega}} \mathord{\rfloor}
    \mathrm{d}\overline{\widetilde{\omega }}=0\,,
    \\
    & \partial _{\tilde{\omega}}\mathord{\rfloor}
    \mathrm{d}\omega=0,\quad
    \partial_{\tilde{\omega}}
    \mathord{\rfloor}
    \mathrm{d}\overline{{\omega}}=0,
    \quad
     {\partial _{\tilde\omega}}\mathord{\rfloor}
   \mathrm{d}\widetilde\omega=1,\quad {\partial
      _{\tilde\omega}}\mathord{\rfloor}
    \mathrm{d}\overline{\widetilde\omega}=0\,,
    \quad\\
    & {\partial _{\overline{\tilde\omega}}}\mathord{\rfloor}
    \mathrm{d}\omega=0, \quad {\partial_{
        \overline{\tilde{\omega}}}}\mathord{\rfloor}
    \mathrm{d}\overline{\omega}=0,\quad
    {\partial _{\overline{\tilde\omega}}}\mathord{\rfloor}
    \mathrm{d}\widetilde\omega=0, \quad {\partial_{
        \overline{\tilde\omega}}}\mathord{\rfloor}
    \mathrm{d}\overline{\widetilde\omega}=1\,.
\end{split}
\end{equation}
Then by \eqref{Sto.12}, \eqref{Sto.7} and \eqref{Sto.13}, we obtain
\begin{align}\label{Sto.14}
 \nu(\theta;\omega,\widetilde\omega)\mathord{\rfloor} \mathrm{d}_{\omega,\widetilde\omega}& \left(\frac{1}{\omega-\zeta} \, a(\omega,\widetilde \omega;h)\,\chi(\omega,\widetilde\omega)\, \e^{\frac{\mathrm{i}}{h}\phi(\omega,\widetilde \omega) }\,\, \frac{\mathrm{d}\widetilde\omega\wedge \mathrm{d}\omega}{2\mathrm{i}}\right)\\
  =& \mathrm{i}(t-\mathrm{i} s)\pi\delta_{\zeta}(\omega)a(\omega,\widetilde \omega;h)\,\chi(\omega,\widetilde\omega)\, \e^{\frac{\mathrm{i}}{h}\phi(\omega,\widetilde \omega) }\,\,\,\frac{\mathrm{d}\overline{\omega}\wedge \mathrm{d}\widetilde \omega}{2\mathrm{i}}\nonumber\\
 & -\mathrm{i}(t+\mathrm{i} s)\pi\delta_{\zeta}(\omega) a(\omega,\widetilde \omega;h)\,\chi(\omega,\widetilde\omega)\, \e^{\frac{\mathrm{i}}{h}\phi(\omega,\widetilde \omega) }\,\,\,\frac{\mathrm{d}\widetilde\omega\wedge \mathrm{d}\omega}{2\mathrm{i}}\nonumber\\
 &-\mathrm{i}(t+\mathrm{i} s) \pi\delta_{\zeta}(\omega) a(\omega,\widetilde \omega;h)\,\chi(\omega,\widetilde\omega)\, \e^{\frac{\mathrm{i}}{h}\phi(\omega,\widetilde \omega) }\,\,\,\frac{\mathrm{d}\overline{\omega}\wedge \mathrm{d}\omega}{2\mathrm{i}}\nonumber\\
 &+0\nonumber\\
 &+\frac{1}{\omega-\zeta} \, a(\omega,\widetilde \omega;h)\, \e^{\frac{\mathrm{i}}{h}\phi(\omega,\widetilde \omega) } 
 \nu(\theta;\omega,\widetilde\omega)\mathord{\rfloor}
 \mathrm{d}_{\omega,\tilde\omega}\Big(\chi(\omega,\widetilde\omega)\frac{\mathrm{d}\widetilde\omega\wedge \mathrm{d}\omega}{2\mathrm{i}}\Big)
\,.\nonumber
\end{align}
Here $(t, s)\in\mathbb{R}^2$ is given in terms of $(\omega,\widetilde \omega)$ as follows. Recall \eqref{Sto.4} and take complex conjugates:
\begin{equation}\label{Sto.15}
 \begin{cases}
  \omega=(t+\mathrm{i} s)+\mathrm{i}\theta (t-\mathrm{i} s) \\
  \overline{\omega}=-\mathrm{i}\theta(t+\mathrm{i} s)+(t-\mathrm{i} s)
 \end{cases}\,.
\end{equation}
Then
\begin{equation}\label{Sto.16}
 \begin{cases}
  t+\mathrm{i} s=\dfrac{1}{1-\theta^2}\big(\omega-\mathrm{i}\theta\overline{\omega}\big)\\
  t-\mathrm{i} s=\dfrac{1}{1-\theta^2}\big(\mathrm{i}\theta\omega+\overline{\omega}\big)
 \end{cases}\,.
\end{equation}
Combining \eqref{Sto.4}, \eqref{Sto.16} gives
\begin{equation}\label{Sto.17}
\begin{cases}
 \widetilde\omega=\dfrac{1}{1-\theta^2}\big(2\mathrm{i}\theta\omega+(1+\theta^2)\overline{\omega}\big)\\
 \overline{\widetilde\omega}=\dfrac{1}{1-\theta^2}\big((1+\theta^2){\omega}-2\mathrm{i}\theta\overline{\omega}\big)
\end{cases}\,.
\end{equation}
Equivalently,
\begin{equation}\label{Sto.17.5}
\begin{cases}
 \overline{\omega} =\dfrac{-2\mathrm{i}\theta}{1+\theta^2}\omega+\dfrac{1-\theta^2}{1+\theta^2}\widetilde\omega\\
 {}\\
 \overline{\widetilde\omega}=\dfrac{1-\theta^2}{1+\theta^2}{\omega}-\dfrac{2\mathrm{i}\theta}{1+\theta^2}\widetilde{\omega}
\end{cases}\,.
\end{equation}

\par
By (\ref{Sto.16}), the expression \eqref{Sto.14} becomes 
\begin{equation}\label{Sto.18}
  \begin{split}
 \nu(\theta;\omega,&\widetilde\omega)\mathord{\rfloor}  \mathrm{d}_{\omega,\tilde\omega} \left(\frac{1}{\omega-\zeta} \, a(\omega,\widetilde \omega;h)\,\chi(\omega,\widetilde\omega)\, \e^{\frac{\mathrm{i}}{h}\phi(\omega,\widetilde \omega) }\,\, \frac{\mathrm{d}\widetilde\omega\wedge\mathrm{d}\omega}{2\mathrm{i}}\right)\\
 =&\dfrac{\mathrm{i}}{1-\theta^2}\big(\mathrm{i}\theta\omega+\overline{\omega}\big)\pi\delta_{\zeta}(\omega) \, a(\omega,\widetilde
 w;h)\,\chi(\omega,\widetilde\omega)\,
 \e^{\frac{\mathrm{i}}{h}\phi(\omega,\widetilde \omega)
 }\,\,\,\frac{\mathrm{d}\overline{\omega}\wedge \mathrm{d}\widetilde \omega}{2\mathrm{i}}
 \\
 & -\dfrac{\mathrm{i}}{1-\theta^2}\big(\omega-\mathrm{i}\theta\overline{\omega}\big)\pi\delta_{\zeta}(\omega) a(\omega,\widetilde \omega;h)\,\chi(\omega,\widetilde\omega)\, \e^{\frac{\mathrm{i}}{h}\phi(\omega,\widetilde \omega) }\,\,\,\frac{\mathrm{d}\widetilde\omega\wedge \mathrm{d}\omega}{2\mathrm{i}}\\
 & -\dfrac{\mathrm{i}}{1-\theta^2}\big(\omega-\mathrm{i}\theta\overline{\omega}\big) \pi\delta_{\zeta}(\omega) a(\omega,\widetilde \omega;h)\,\chi(\omega,\widetilde\omega)\, \e^{\frac{\mathrm{i}}{h}\phi(\omega,\widetilde \omega) }\,\,\,\frac{\mathrm{d}\overline{\omega}\wedge \mathrm{d}\omega}{2\mathrm{i}}\\
 &+\frac{1}{\omega-\zeta} \, a(\omega,\widetilde \omega;h)\, \e^{\frac{\mathrm{i}}{h}\phi(\omega,\widetilde \omega) } 
 \nu(\theta;\omega,\widetilde\omega)\mathord{\rfloor}
 \mathrm{d}_{\omega,\tilde\omega}\Big(\chi(\omega,\widetilde\omega)\frac{\mathrm{d}\widetilde\omega\wedge \mathrm{d}\omega}{2\mathrm{i}}\Big)\,,
 \end{split}
\end{equation}
or equivalently
\begin{align}\label{Sto.19}
 \hskip-6pt\nu(\theta;\omega,\widetilde\omega)\mathord{\rfloor}  \mathrm{d}_{\omega,\tilde\omega} \left(\frac{1}{\omega-\zeta} \, a(\omega,\widetilde \omega;h)\,\chi(\omega,\widetilde\omega)\, \e^{\frac{\mathrm{i}}{h}\phi(\omega,\widetilde \omega) }\,\, \frac{\mathrm{d}\widetilde\omega\wedge \mathrm{d}\omega}{2\mathrm{i}}\right)
 =S_{l}+S_{r}
\end{align}
with
\begin{align}\label{Sto.20}
  S_{l} :=\dfrac{\mathrm{i}}{1-\theta^2}&\pi\delta_{\zeta}(\omega)a(\omega,\widetilde \omega;h)\,\chi(\omega,\widetilde\omega)\, \e^{\frac{\mathrm{i}}{h}\phi(\omega,\widetilde \omega) }\times\\
 &\left[\big(\mathrm{i}\theta\omega+\overline{\omega}\big)\,\,\,\frac{\mathrm{d}\overline{\omega}\wedge \mathrm{d}\widetilde \omega}{2\mathrm{i}}
 -\big(\omega-\mathrm{i}\theta\overline{\omega}\big)\,\,\,\Big(\frac{\mathrm{d}\widetilde\omega\wedge \mathrm{d}\omega}{2\mathrm{i}}
 +\frac{\mathrm{d}\overline{\omega}\wedge \mathrm{d}\omega}{2\mathrm{i}}\Big)\right]\nonumber
\end{align}
and
\begin{equation}\label{Sto.21}
 S_{r}:=\frac{1}{\omega-\zeta} \, a(\omega,\widetilde \omega;h)\, \e^{\frac{\mathrm{i}}{h}\phi(\omega,\widetilde \omega) } 
 \nu(\theta;\omega,\widetilde\omega)\mathord{\rfloor}
\mathrm{d}_{\omega,\tilde\omega}\Big(\chi(\omega,\widetilde\omega)\,\frac{\mathrm{d}\widetilde\omega\wedge \mathrm{d}\omega}{2\mathrm{i}}\Big)\,.
\end{equation}

\par 
In order to study $\displaystyle\iintJ_{\Gamma_\theta}\hskip-4pt S_\ell$ we shall express the
restrictions to $\Gamma_\theta$ of $\mathrm{d}\overline{\omega } \wedge
\mathrm{d}\widetilde{\omega }$, $\mathrm{d}\widetilde{\omega }\wedge \mathrm{d}\omega $, $\mathrm{d}\overline{\omega }\wedge \mathrm{d}\omega  $ as
multiples of $L(\mathrm{d}\omega )$.
Identities \eqref{Sto.17} and \eqref{Sto.4} give
\begin{align}\label{Sto.22}
& L(\mathrm{d}\widetilde\omega):=\frac{\mathrm{d}\overline{\widetilde\omega}\wedge \mathrm{d}\widetilde\omega}{2\mathrm{i}}=-L(\mathrm{d}\omega )=-(1-\theta^2)\,\mathrm{d}t\wedge \mathrm{d}s\nonumber\\
&  \frac{\mathrm{d}\widetilde\omega\wedge \mathrm{d}{\omega}}{2\mathrm{i}}=\frac{1+\theta^2}{1-\theta^2}L(\mathrm{d}\omega)=(1+\theta^2)\,\mathrm{d}t\wedge \mathrm{d}s\\
&\frac{\mathrm{d}\overline{\omega}\wedge \mathrm{d}\widetilde \omega}{2\mathrm{i}}=\frac{2\mathrm{i}\theta}{1-\theta^2}L(\mathrm{d}\omega)=2\mathrm{i}\theta \,\mathrm{d}t\wedge \mathrm{d}s\,.\nonumber
\end{align}
Then \eqref{Sto.20} becomes
\begin{equation}\label{Sto.23}
 \hskip-9pt S_{l}=\dfrac{-2\mathrm{i}\left[(1+\theta^2)\omega-2\mathrm{i}\theta\overline{\omega}\right]}{(1-\theta^2)^2}\Big(\pi\delta_{\zeta}(\omega)a(\omega,\widetilde
 \omega ;h)\,\chi(\omega,\widetilde\omega)\, \e^{\frac{\mathrm{i}}{h}\phi(\omega,\widetilde \omega) }\Big)
 \,\, L(\mathrm{d}\omega )\,.
\end{equation}
Recalling \eqref{Sto.17}, we get
\begin{equation}\label{Sto.24}
 S_{l}=\dfrac{-2\mathrm{i}\,\overline{\widetilde{\omega}}}{1-\theta^2}\Big(\pi\delta_{\zeta}(\omega) \, a(\omega,\widetilde
 \omega ;h)\,\chi(\omega,\widetilde\omega)\, \e^{\frac{\mathrm{i}}{h}\phi(\omega,\widetilde \omega) }\Big)
 \,\, L(\mathrm{d}\omega )\,.
\end{equation}
This gives
\begin{equation}\label{Sto.25}
 -\frac{1}{\pi} \iintJ_{\Gamma _\theta}S_l=
\dfrac{2\mathrm{i}\,\overline{\widetilde{\zeta}}}{1-\theta^2}
\, a(\zeta,\widetilde \zeta;h)\,\chi(\zeta,\widetilde\zeta)\, \e^{\frac{\mathrm{i}}{h}\phi(\zeta,\widetilde \zeta)}\,,
\end{equation}
where $\widetilde\zeta$ is determined by the condition that $(\zeta,\widetilde\zeta)\in \Gamma_{\theta}$, i.e. by \eqref{Sto.17} with $\omega,\widetilde\omega$ replaced by $\zeta,\widetilde\zeta$:
\begin{equation}\label{Sto.26}
 \widetilde\zeta=\dfrac{1}{1-\theta^2}\,\,\Big(2\mathrm{i}\theta\zeta+(1+\theta^2)\overline{\zeta}\Big)\,,
\end{equation}
i.e.
\begin{align}\label{Sto.25.5}
 -\frac{1}{\pi} \iintJ_{\Gamma_\theta}\hskip-1pt S_l =
 &\frac{2\mathrm{i}\Big((1+\theta^2){\zeta}-2\mathrm{i}\theta \overline{\zeta}\Big)}{(1-\theta^2)^2}\,\,\,\,
a\Big(\zeta,\frac{1}{1-\theta^2}\,\big(2\mathrm{i}\theta\zeta+(1+\theta^2)\overline{\zeta}\big);h\Big)\,\times\\
& \chi\Big(\zeta,\frac{1}{1-\theta^2}\,\big(2\mathrm{i}\theta\zeta+(1+\theta^2)\overline{\zeta}\big)\Big)\, \e^{\frac{\mathrm{i}}{h}\phi\big(\zeta,\frac{1}{1-\theta^2}\,\,(2\mathrm{i}\theta\zeta+(1+\theta^2)\overline{\zeta}\big)}\,.\nonumber
\end{align}

\par
We next show that the restriction of $\nu\mathord{\rfloor} \mathrm{d}\Big(\chi(\omega,\widetilde\omega)\mathrm{d}\widetilde\omega\wedge \mathrm{d}\omega\Big)$ to $\Gamma_\theta$
for $0\leq \theta\leq 1$, is independent of the choice of smooth deformation field $\nu$ for the family $\Gamma_{\theta}$: Indeed let $\widetilde\nu$ be a second deformation field for the same family. Then $\mu:=\widetilde\nu-\nu$ is a real vector field which is tangent to $\Gamma_\theta$ for any fixed $\theta\in[0,1]$ and it suffices to show that
\begin{equation}\label{1.sto6}
\left(\mu \mathord{\rfloor}\mathrm{d}\Big(\chi(\omega,\widetilde\omega)\mathrm{d}\widetilde\omega\wedge \mathrm{d}\omega\Big)\right)_{\big|\displaystyle\Gamma_{\theta}}=0\,.
\end{equation}

\noindent
\textbf{Proof} of (\ref{1.sto6}). Since $\chi(\omega,\widetilde\omega)\mathrm{d}\widetilde\omega\wedge \mathrm{d}\omega_{\big|\displaystyle\Gamma_{\theta}}$ is a 2-form and hence of maximal degree, we have 
$\mathrm{d}_{\omega,\tilde\omega}\Big( \chi(\omega,\widetilde\omega)\mathrm{d}\widetilde\omega\wedge \mathrm{d}\omega\Big)_{\big|\displaystyle\Gamma_{\theta}}=0$ and since $\mu:=\nu-\widetilde\nu$ is tangential to $\Gamma_\theta$ we have 
\begin{equation}\label{Sto.28}
  \begin{split}
 \left(\mu \mathord{\rfloor}\mathrm{d}\Big(\chi(\omega,\widetilde\omega)\mathrm{d}\widetilde\omega\wedge \mathrm{d}\omega\Big)\right)_{\big|\displaystyle\Gamma_{\theta}} &=\mu \mathord{\rfloor}\left(\mathrm{d}\big(\chi(\omega,\widetilde\omega)\mathrm{d}\widetilde\omega\wedge \mathrm{d}\omega\big)\right)_{\big|\displaystyle\Gamma_{\theta}}\\
 &=\mu \mathord{\rfloor}\mathrm{d}\Big(\chi(\omega,\widetilde\omega)\mathrm{d}\widetilde\omega\wedge \mathrm{d}\omega_{\big|\displaystyle\Gamma_{\theta}}\Big)\\
 &=0\,.
 \end{split}
\end{equation}
Here we used that if $\mu $ is a vector field tangent to $\Gamma _\theta
$ and $\Omega $ is a $k$-form, then
$$
\mu \rfloor \left({{\Omega}_\vert}_{\Gamma _\theta }  \right)=
{{(\mu \rfloor \Omega )}_\vert}_{\Gamma _\theta }.
$$
Indeed, this is clear when $\Omega $ is a 1-form and in the general
case we get it when $\Omega$ is an exterior product of 1-forms
and
more generally when $\Omega $ is a linear combination of such exterior
products.\hfill{$\square$}

\par\medskip
This gives
\begin{lemma}\label{Err}
  Let $0\le \theta <1$,  let $\nu (\theta )$ be as in (\ref{Sto.7})
  and $\widetilde{\nu }(\theta )=\nu (\theta )+\mu (\theta )$ where
  $\mu (\theta )$ is a smooth vector field tangent to $\Gamma_\theta$.
Then on $\Gamma_\theta$, 
 \begin{equation}\label{Sto.27}
 S_{r}= \frac{1}{\omega-\zeta} \, a(\omega,\widetilde \omega;h)\, \e^{\frac{\mathrm{i}}{h}\phi(\omega,\widetilde \omega) } 
 \widetilde\nu(\theta;\omega,\widetilde\omega)\mathord{\rfloor}
 \mathrm{d}_{\omega,\tilde\omega}\Big( \chi(\omega,\widetilde\omega) \,\frac{\mathrm{d}\widetilde\omega\wedge \mathrm{d}\omega}{2\mathrm{i}}\Big)\,.
\end{equation}
\end{lemma}

\par We review some results of this section in the following
proposition.
\begin{prop}\label{1sto7}
Define the contours $\Gamma_\theta \subset \mathbb{C}^2_{\omega
  ,\tilde{\omega }}$ by (\ref{Sto.5}), (\ref{Sto.4}). ${{\Omega
  }_\vert}_{\Gamma _\theta }$ is well defined for $0\le \theta <1$
(see (\ref{1.sto5})). Let $\nu $ be the deformation vector field on
$\Gamma _\theta $ in (\ref{2.sto5}), (\ref{Sto.7}). With $I=I(\Gamma
_\theta ,\zeta )$ defined as in (\ref{Sto.9}), $\frac{\mathrm{d}I}{\mathrm{d}\theta}$ is given by (\ref{Sto.8}) for $0\le\theta <1$. The
integrand in (\ref{Sto.8}) is of the form $S_l+S_r$, (see \eqref{Sto.19}--\eqref{Sto.21}): For
$0\le \theta <1$, we have
  \begin{align}\label{1.sto8}
  \frac{\mathrm{d}I}{\mathrm{d}\theta}
 =&-\frac{1}{\pi}\iintJ_{\Gamma_\theta}\nu(\theta ;\omega,\widetilde\omega)\mathord{\rfloor}  \mathrm{d}_{\omega,\tilde\omega}\left(\frac{1}{\omega-\zeta} \, a(\omega,\widetilde \omega;h)\,\chi(\omega,\widetilde\omega)\, \e^{\frac{\mathrm{i}}{h}\phi(\omega,\widetilde \omega)}\, \frac{\mathrm{d}\widetilde\omega\wedge \mathrm{d}\omega}{2\mathrm{i}}\right),\nonumber\\
=&-\frac{1}{\pi}\iintJ_{\Gamma_\theta} S_l-\frac{1}{\pi}\iintJ_{\Gamma_\theta}S_r\nonumber\\
=& \frac{2\mathrm{i} \, \overline{\widetilde{\zeta }}}{1-\theta ^2} \, a(\zeta
,\widetilde{\zeta };h)\chi (\zeta ,\widetilde{\zeta
})\, \e^{\frac{\mathrm{i}}{h}\phi (\zeta ,\widetilde{\zeta})}\\
&-\frac{1}{\pi}\iintJ_{\Gamma_\theta}
\frac{1}{\omega -\zeta} \, a(\omega ,\widetilde{\omega
};h)\, \e^{\frac{\mathrm{i}}{h}\phi(\omega ,\widetilde{\omega })} \nu (\theta;\omega,\widetilde\omega)\rfloor
\mathrm{d}_{\omega ,\tilde{\omega }}\Big(\chi(\omega ,\widetilde{\omega})\, \frac{\mathrm{d}\widetilde{\omega }\wedge \mathrm{d}\omega}{2\mathrm{i}}\Big)\,,\nonumber
\end{align}
where $\widetilde{\zeta }$ is determined by the condition that $(\zeta
,\widetilde{\zeta })\in \Gamma _\theta $ (see \eqref{Sto.26}). Here
we also used \eqref{Sto.25} and \eqref{Sto.27}. $\nu $ can be
replaced by any other smooth vector field of the form $\widetilde{\nu
}=\nu +\mu $ with $\mu $ tangent to $\Gamma _\theta $.
\end{prop}

\par 
We can integrate (\ref{1.sto8}) from $\theta =0$ to $\theta
=1-\delta $, $0<\delta <1$ and get

\begin{align}\label{2.sto8}
-\frac{1}{\pi}&\iintJ_{\mathrm{adiag\,}({\mathbb{C}}^2)} \frac{1}{\omega-\zeta} \, a(\omega ,\widetilde{\omega };h)\chi (\omega,\widetilde{\omega })\, \e^{\frac{\mathrm{i}}{h}\phi (\omega ,\widetilde{\omega })} 
\, \frac{\mathrm{d}\widetilde{\omega}\wedge \mathrm{d}\omega}{2\mathrm{i}}\nonumber\\ 
=&-\frac{1}{\pi}\iintJ_{\Gamma_{1-\delta}}\frac{1}{\omega-\zeta} \, 
a(\omega ,\widetilde{\omega };h)\chi(\omega,\widetilde{\omega}) \, \e^{\frac{\mathrm{i}}{h}\phi(\omega ,\widetilde{\omega })} \,
\frac{\mathrm{d}\widetilde{\omega }\wedge \mathrm{d}\omega}{2\mathrm{i}}
\\
&-\int_0^{1-\delta} \frac{2\mathrm{i}\,\overline{\widetilde{\zeta}}}{1-\theta ^2} \, a(\zeta,\widetilde{\zeta };h)\chi (\zeta ,\widetilde{\zeta
})\, \e^{\frac{\mathrm{i}}{h}\phi (\zeta ,\widetilde{\zeta})} \, \mathrm{d}\theta \nonumber\\
&+\int_0^{1-\delta}\frac{1}{\pi}\iintJ_{\Gamma_\theta}
\frac{1}{\omega -\zeta } \, a(\omega ,\widetilde{\omega
};h)\, \e^{\frac{\mathrm{i}}{h}\phi (\omega ,\widetilde{\omega })}\widetilde\nu
(\theta;\omega,\widetilde\omega)\rfloor
\mathrm{d}_{\omega ,\tilde{\omega}}\left(\chi(\omega ,\widetilde{\omega}) \frac{\mathrm{d}\widetilde{\omega }\wedge
\mathrm{d}\omega}{2\mathrm{i}}\right)\mathrm{d}\theta \nonumber\\
&=:\mathrm{I}(\zeta ,1-\delta)+\mathrm{II}(\zeta ,1-\delta)+\mathrm{III}(\zeta,1-\delta)\,.\nonumber
\end{align}

\par
If we could pass to the limit $\delta =0$, we would get
\begin{align}\label{3.sto8}
-\frac{1}{\pi}&\iintJ_{\mathrm{adiag\,}({\mathbb{C}}^2)} \frac{1}{\omega-\zeta} \, a(\omega ,\widetilde{\omega };h)\chi (\omega,\widetilde{\omega })\, \e^{\frac{\mathrm{i}}{h}\phi (\omega ,\widetilde{\omega })} 
\, \frac{\mathrm{d}\widetilde{\omega}\wedge \mathrm{d}\omega}{2\mathrm{i}}\nonumber\\ 
=&-\frac{1}{\pi}\iintJ_{\Gamma_{1}}\frac{1}{\omega-\zeta} \, 
a(\omega ,\widetilde{\omega };h)\chi(\omega,\widetilde{\omega}) \, \e^{\frac{\mathrm{i}}{h}\phi(\omega ,\widetilde{\omega })} \,
\frac{\mathrm{d}\widetilde{\omega }\wedge \mathrm{d}\omega}{2\mathrm{i}}
\\
&-\int_0^{1} \frac{2\mathrm{i}\,\overline{\widetilde{\zeta}}}{1-\theta ^2} \, a(\zeta,\widetilde{\zeta };h)\chi (\zeta ,\widetilde{\zeta
})\, \e^{\frac{\mathrm{i}}{h}\phi (\zeta ,\widetilde{\zeta})} \, \mathrm{d}\theta \nonumber\\
&+\int_0^{1}\frac{1}{\pi}\iintJ_{\Gamma_\theta}
\frac{1}{\omega -\zeta } \, a(\omega ,\widetilde{\omega
};h)\, \e^{\frac{\mathrm{i}}{h}\phi (\omega ,\widetilde{\omega })}\widetilde\nu
(\theta;\omega,\widetilde\omega)\rfloor
\mathrm{d}_{\omega ,\tilde{\omega}}\left(\chi(\omega ,\widetilde{\omega}) \frac{\mathrm{d}\widetilde{\omega }\wedge
\mathrm{d}\omega}{2\mathrm{i}}\right)\mathrm{d}\theta \nonumber
\end{align}
and for the first term in the right hand side, we
could use that $\Gamma_1$ is a Cartesian product and apply Fubini's
theorem.

\section{The term \texorpdfstring{$\mathrm{I}(\zeta ,1-\delta)$}{} in (\ref{2.sto8})}\label{ft}
\setcounter{equation}{0}

In this section the leading order contributions to 
the first term in the decomposition of Section \ref{Sto} are 
identified.

We shall work in coordinates that are well adapted to $\Gamma
_\theta$.
\begin{lemma}\label{1ft1}
Let $\theta \in [0,1]$. Then
\begin{equation}\label{1.ft1}
  \begin{split}
&e_+:=\frac{1}{\sqrt{2}}\left(\frac{1-\theta }{\sqrt{1+\theta
      ^2}}\, \e^{-\mathrm{i} \frac{\pi}{4}} \, , \, \frac{1+\theta }{\sqrt{1+\theta^2}}\, \e^{\mathrm{i} \frac{\pi}{4}}\right)\in \Gamma _\theta \, ,\\
&e_-:=\frac{1}{\sqrt{2}}\left(\frac{1+\theta }{\sqrt{1+\theta
      ^2}} \, \e^{\mathrm{i}\frac{\pi}{4}} \, , \, \frac{1-\theta }{\sqrt{1+\theta^2}}\, \e^{-\mathrm{i}\frac{\pi}{4}}\right)\in \Gamma_\theta
  \end{split}
\end{equation}
form an orthonormal basis in $\mathbb{C}^2$,
hence also an orthonormal basis in $\Gamma _\theta $ as a Euclidean
space.
\end{lemma}

\noindent
\textbf{Proof.} 
  By (\ref{Sto.4}), $\Gamma _\theta $ is given by
  \[\begin{cases}
\omega =\tau +\mathrm{i}\theta \overline{\tau }\\ \widetilde{\omega
}=\overline{\tau }+\mathrm{i}\theta \tau \end{cases}\,, \, \tau =t+\mathrm{i} s
\,\,\,\, \mathrm{with}\,\, \ t,s\in \mathbb{R}.
\]
Removing a common positive prefactor, put
\[\begin{split}
\widetilde{e}_+&=((1-\theta ) \, \e^{-\mathrm{i}\frac{\pi}{4}},(1+\theta ) \, \e^{\mathrm{i}\frac{\pi}{4}})=:(\omega _+,\widetilde{\omega }_+)\, ,\\
\widetilde{e}_-&=((1+\theta) \, \e^{\mathrm{i}\frac{\pi}{4}},(1-\theta) \, 
\e^{-\mathrm{i}\frac{\pi}{4}})=:(\omega _-,\widetilde{\omega }_-)\,.
\end{split}
\]
Then
$$
\begin{cases}
\omega_+=\e^{-\mathrm{i}\frac{\pi}{4}}-\theta \e^{-\mathrm{i}\frac{\pi}{4}}=\e^{-\mathrm{i}\frac{\pi}{4}}+\mathrm{i}\theta
\e^{\mathrm{i}\frac{\pi}{4}}=\tau +\mathrm{i}\theta \overline{\tau }\\
\widetilde{\omega}_+=\e^{\mathrm{i}\frac{\pi}{4}}+\theta \e^{\mathrm{i}\frac{\pi}{4}}=\e^{\mathrm{i}\frac{\pi}{4}}+\mathrm{i}\theta \e^{-\mathrm{i}\frac{\pi}{4}}=\overline{\tau }+\mathrm{i}\theta \tau
\end{cases}\,,
\qquad \tau =\e^{-\mathrm{i}\frac{\pi}{4}}\, ,
$$
$$
\begin{cases}
\omega _-=\e^{\mathrm{i}\frac{\pi}{4}}+\theta \e^{\mathrm{i}\frac{\pi}{4}}=\e^{\mathrm{i}\frac{\pi}{4}}+\mathrm{i}\theta
\e^{-\mathrm{i}\frac{\pi}{4}}=\tau +\mathrm{i}\theta \overline{\tau}\\
\widetilde{\omega}_-=\e^{-\mathrm{i}\frac{\pi}{4}}-\theta \e^{-\mathrm{i}\frac{\pi}{4}}=\e^{-\mathrm{i}\frac{\pi}{4}}+\mathrm{i}\theta \e^{\mathrm{i}\frac{\pi}{4}}=\overline{\tau}+\mathrm{i}\theta \tau
\end{cases}\,,
\quad \tau =\e^{\mathrm{i}\frac{\pi}{4}}\,,
$$
so $\widetilde{e}_\pm$ and hence also $e_\pm$ belong to $\Gamma
_\theta $. By direct calculation, we see that $e_+$, $e_-$ form an
orthonormal basis in $\mathbb{C}^2$ and the lemma follows.
\hfill$\square$

\par\medskip
Every $(\omega ,\widetilde{\omega} )\in \Gamma _\theta $
has a unique decomposition
\begin{equation}\label{2.ft1}
(\omega ,\widetilde{\omega } )=w_+e_++w_-e_-\,,
\end{equation}
where $w_\pm =((\omega ,\widetilde{\omega })|e_\pm)_{\mathbb{C}^2}$ are
real (and this still holds for $(\omega ,\widetilde{\omega })\in {\mathbb 
  C}^2$ now with $w_\pm\in \mathbb{C}$).
We have more direct formulae for $w_\pm$ by looking at the
$\omega $-component of  (\ref{1.ft1}), (\ref{2.ft1}),
\begin{equation}\label{3.ft1}
  \omega = w_+\frac{1-\theta}{\sqrt{2(1+\theta ^2)}} \e^{-\mathrm{i}\frac{\pi}{4}}+w_-\frac{1+\theta }{\sqrt{2(1+\theta ^2)}} \e^{\mathrm{i}\frac{\pi}{4}}
\end{equation}
and using that $(\e^{-\mathrm{i}\frac{\pi}{4}}, \e^{\mathrm{i}\frac{\pi}{4}})$ is an orthonormal basis in $\mathbb{C}\simeq \mathbb{R}^2$, leading to
\begin{equation}\label{4.ft1}
w_+=\frac{\sqrt{2(1+\theta ^2)}}{1-\theta }\mathrm{Re} (\omega \e^{\mathrm{i}\frac{\pi}{4}}),\quad
w_-=\frac{\sqrt{2(1+\theta ^2)}}{1+\theta }\mathrm{Re} (\omega \e^{-\mathrm{i}\frac{\pi}{4}})\,.
\end{equation}

\par 
Similarly, by taking the $\widetilde{\omega }$ components in
(\ref{1.ft1}), (\ref{2.ft1}),
\begin{equation}\label{4,5.ft1}
\widetilde{\omega}=w_+\frac{1+\theta }{\sqrt{2(1+\theta ^2)}}\e^{\mathrm{i}\frac{\pi}{4}}
+w_-\frac{1-\theta}{\sqrt{2(1+\theta ^2)}}\e^{-\mathrm{i}\frac{\pi}{4}}\,,
\end{equation}
\begin{equation}\label{5.ft1}
w_+=\frac{\sqrt{2(1+\theta^2)}}{1+\theta}\mathrm{Re} (\widetilde{\omega}
\e^{-\mathrm{i}\frac{\pi}{4}}),\quad
w_-=\frac{\sqrt{2(1+\theta^2)}}{1-\theta}\mathrm{Re} (\widetilde{\omega}
\e^{\mathrm{i}\frac{\pi}{4}})\,.
\end{equation}

\par 
Taking the differentials of $\omega $,
  $\widetilde{\omega }$ with respect to $w_+$, $w_-$ in
(\ref{3.ft1}), (\ref{4,5.ft1}), for any
fixed $\theta $, we get 
\begin{equation}\label{6.ft1}
\mathrm{d}\omega \wedge \mathrm{d}\widetilde{\omega }=\frac{1}{\mathrm{i}}\mathrm{d}w_+\wedge \mathrm{d}w_-\, ,
\end{equation}
and we can view this as a relation between differential forms on
$\Gamma _\theta $.

\par 
As for $\zeta \in \mathbb{C}$ in (\ref{2.sto8}) we have no
$\widetilde{\zeta }$ and we choose the simpler decomposition
\begin{equation}\label{7.ft1}
\zeta =z_+\e^{-\mathrm{i}\frac{\pi}{4}}+z_-\e^{\mathrm{i}\frac{\pi}{4}},
\end{equation}
\begin{equation}\label{8.ft1}
  z_+=\mathrm{Re} (\zeta \e^{\mathrm{i}\frac{\pi}{4}}),\quad z_-=\mathrm{Re} (\zeta \e^{-\mathrm{i}\frac{\pi}{4}})\,.
\end{equation}

\par The first term $\mathrm{I}(\zeta ,1-\delta)$ in (\ref{2.sto8}) becomes
\begin{equation}\label{9.ft1}
  \mathrm{I}(\zeta ,1-\delta)=\frac{1}{2\pi }
  \iintJ \frac{a\chi \e^{\frac{\mathrm{i}}{h}\phi} \mathrm{d}w_+\wedge \mathrm{d}w_-}{
    \e^{-\mathrm{i}\frac{\pi}{4}}(z_+-\frac{1-\theta}{\sqrt{2(1+\theta ^2)}}w_+)+\e^{\mathrm{i}\frac{\pi}{4}}(z_--\frac{1+\theta }{\sqrt{2(1+\theta ^2)}}w_-)
  }
\end{equation}
with $\theta =1-\delta $.

For $\zeta $ in a fixed bounded set in $\mathbb{C}$ or equivalently
for $(z_+,z_-)$ in a fixed bounded set in $\mathbb{R}^2$, we have
\begin{equation}\label{1.ft2}\begin{split}
\mathrm{I}(\zeta ,1-\delta)=&{\mathcal  O}(1) \e^{\frac{F}{h}}\hskip-3pt
\iintJ \frac{|\chi (\omega ,\widetilde{\omega })|\mathrm{d}w_+\mathrm{d}w_-}
{|z_+-\frac{1-\theta }{\sqrt{2(1+\theta ^2)}}w_+|
+|z_--\frac{1+\theta }{\sqrt{2(1+\theta ^2)}}w_-|}\,,\\
&\hbox{ where }F:=\sup_{\Gamma _\theta }(-\mathrm{Im}\, \phi)\,.
\end{split}
\end{equation}
Here $w_\pm ={\mathcal  O}(1)$ on $\mathrm{supp\,}\chi $ by (\ref{4.ft1}), (\ref{5.ft1}) so the last integral is
\begin{equation}\label{2.ft2}
{\mathcal  O}(1)\int_J
\Big(1+\big|\ln |z_+-\frac{\delta}{\sqrt{2(1+\theta ^2)}}w_+|
\big|\Big)\, \mathrm{d}w_+\,,
\end{equation}
where $J$ is a fixed bounded interval, independent of 
$z_\pm$ and $\theta =1-\delta$. This can be written
\begin{equation}\label{3.ft2}
{\mathcal  O}\left(\frac{1}{\delta} \right)\int_{\widetilde{J}_\delta
}\Big(1+\big|\ln |z_+-\widetilde{w}_+|\big| \Big) \, \mathrm{d}\widetilde{w}_+\,,
\end{equation}
where
\begin{equation}\label{4.ft2}
\widetilde{J}_\delta =\frac{\delta}{\sqrt{2(1+\theta ^2)}}\,J\,.
\end{equation}
Here $\widetilde{J}_\delta$ is an interval of length ${\mathcal 
  O}(\delta)$. If $\mathrm{dist\,}(z_+,\widetilde{J}_\delta)\ge
\delta$, we can estimate the expression (\ref{3.ft2}) by
$$
\le {\mathcal  O}(1)\left( 1+|\ln
  (\mathrm{dist\,}(z_+,\widetilde{J}_\delta))| \right).
$$

If $\mathrm{dist\,}(z_+,\widetilde{J}_\delta)<\delta$, we
compare with the ``worst'' case when $z_+$ is the midpoint of
$\widetilde{J}_\delta$ and estimate (\ref{3.ft2}) by
$$
\le \frac{{\mathcal  O}(1)}{\delta}\int_0^{{\mathcal  O}(\delta)} \big(1+|\ln
t|\big)\, \mathrm{d}t=\frac{{\mathcal  O}(1)}{\delta}\Big[2t-t\ln t\Big]_0^{{\mathcal O}(\delta)}
\le {\mathcal  O}(1)\ln \Big(\frac{1}{\delta}\Big)\,.
$$

In conclusion,
\begin{equation}\label{5.ft2}
  |\mathrm{I}(\zeta ,1-\delta)|\le {\mathcal  O}(1)\e^{\frac{F}{h}}\times \begin{cases}
\left( 1+|\ln(\mathrm{dist\,}(z_+,\widetilde{J}_\delta))| \right),
\hbox{ when } \, \, \mathrm{dist\,}(z_+,\widetilde{J}_\delta)\ge \delta\,,\\
\ln (\frac{1}{\delta}),\hbox{ when } \, \, \mathrm{dist\,}(z_+,\widetilde{J}_\delta)\le \delta\,.
  \end{cases}
\end{equation}
This gives\begin{equation}\label{6.ft2}
|\mathrm{I}(\zeta ,1-\delta)|\le {\mathcal  O}(1)\e^{\frac{F}{h}}\big(1+|\ln(|z_+|+\delta)|\big)\,.
\end{equation}

We next consider $\underset{\delta \to 0}{\lim}\mathrm{I}(\zeta,1-\delta)$. Notice that when $\delta=0$, (\ref{3.ft1}), (\ref{4,5.ft1}) reduce to
\begin{equation}\label{1.ft3}
\omega =w_-\e^{\mathrm{i}\frac{\pi}{4}},\quad \widetilde{\omega}=w_+\e^{\mathrm{i}\frac{\pi}{4}}\,.
\end{equation}
When $z_+\ne 0$, we can let $\theta \to 1$ ($\delta \to 0$) in
(\ref{9.ft1}) and get
\begin{equation}\label{2.ft3}
\lim_{\delta \to 0}\mathrm{I}(\zeta ,1-\delta)=\mathrm{I}(\zeta,1)
=\frac{1}{2\pi}\iintJ \frac{a\chi \e^{\mathrm{i}\frac{\phi}{h}}}
{\e^{-\mathrm{i}\frac{\pi}{4}}z_+ +\e^{\mathrm{i}\frac{\pi}{4}}(z_--w_-)}\, \mathrm{d}w_+\mathrm{d}w_-\,.
\end{equation}

We sum up the results so far without using the coordinates
$w_\pm$, $z_\pm$. From (\ref{7.ft1}) we get
\begin{equation}\label{3.ft3}
|z_+|=\mathrm{dist\,}(\zeta,\e^{\mathrm{i}\frac{\pi}{4}}\mathbb{R})\,.
\end{equation}
\begin{prop}\label{1ft3}
Define $\mathrm{I}(\zeta ,1-\delta)$ as in Proposition \ref{1sto7}
and let $\zeta $ vary in a compact set $K\subset \mathbb{C}$. With $F$ as
in (\ref{1.ft2}) we have
\begin{equation}\label{4.ft3}
  |\mathrm{I}(\zeta ,1-\delta)|\le {\mathcal O}(1)\e^{\frac{F}{h}}\left(
1+|\ln (\mathrm{dist\,}(\zeta,\e^{\mathrm{i}\frac{\pi}{4}}\mathbb{R})+\delta)|
  \right).
\end{equation}
When $\zeta \in K\setminus (\e^{\mathrm{i}\frac{\pi}{4}}\mathbb{R})$,
\begin{equation}\label{4,5.ft3}
\mathrm{I}(\zeta,1)=-\frac{1}{\pi}\iintJ_{\Gamma_{1}}\frac{1}{\omega -\zeta}\, (a\chi)(\omega ,\widetilde{\omega}) \e^{\frac{\mathrm{i}}{h}\phi (\omega,\widetilde{\omega})}\,\frac{\mathrm{d}\widetilde{\omega}\wedge \mathrm{d}\omega}{2\mathrm{i}} 
\end{equation}
is finite and
\begin{equation}\label{5.ft3}
\mathrm{I}(\zeta ,1-\delta)\longrightarrow \mathrm{I}(\zeta,1),\quad \delta \to 0\,.
\end{equation}
\end{prop}

\par
We now assume that $\phi$ is holomorphic in
$\mathrm{neigh}((0,0);\mathbb{C}^2)$ (either in $\mathbb{C}^2_{\xi ,\eta
}$ or in $\mathbb{C}^2_{\omega ,\tilde{\omega }}$, depending on the
which coordinates we use). Assume that
\begin{equation}\label{6.ft3}
\phi (0,0)=0,\ \ \nabla \phi (0,0)=0.
\end{equation}
Let $\phi_2(\xi ,\eta )$ be the 2nd order Taylor polynomial of $\phi
$ at $(0,0)$,
\begin{equation}\label{1.ft4}
\phi (\xi,\eta )=\phi_2(\xi ,\eta )+{\mathcal  O}(|(\xi ,\eta )|^3)\,.
\end{equation}
We use the same notations for functions in the variables $\omega
,\widetilde{\omega }$,
\begin{equation}\label{1,5.ft4}
  \phi (\omega ,\widetilde{\omega } )=\phi_2(\omega ,\widetilde{\omega
  })+{\mathcal  O}(|(\omega ,\widetilde{\omega })|^3)\,.
\end{equation}
Assume that $\phi_2$ is real-valued on $\mathbb{R}^2_{\xi ,\eta }$,
hence
\begin{equation}\label{2.ft4}
\phi_2(\xi ,\eta )=\frac{\lambda }{2}\xi ^2-\frac{\mu }{2}\eta
^2+\rho\, \xi \eta ,\ \ \lambda ,\mu ,\rho \in \mathbb{R}\,.
\end{equation}
Assume that
\begin{equation}\label{3.ft4}
  \lambda,\, \mu >0. \end{equation}
Then
\begin{equation}\label{4.ft4}
\det \phi_2'' =\det \begin{pmatrix}\lambda &\rho \\ \rho
  &-\mu \end{pmatrix}
=
-(\lambda \mu +\rho ^2)<0\,,
\end{equation}
so the quadratic form $\phi_2$ is non-degenerate of signature $(+,-)$.
Conversely, if $\phi_2$ is a real quadratic form which is
non-degenerate of signature $(+,-)$, then we have (\ref{2.ft4}) and
$\lambda \mu +\rho ^2>0$. Of course, in this case we can perform a rotation
in $\mathbb{R}^2$ so that we get (\ref{2.ft4}) in the new coordinates
with $\rho =0$, $\lambda,\, \mu  >0$.

\par 
Recall that $\Gamma_\theta $ is given by (\ref{Sto.3}),
\begin{equation}\label{4,5.ft4}
\xi =(1+\mathrm{i}\theta )t,\quad \eta =(1-\mathrm{i}\theta )s,\quad t,s\in \mathbb{R}\,,
\end{equation}
so that
$$
|\xi |^2=(1+\theta ^2)t^2,\quad |\eta |^2=(1+\theta ^2)s^2\hbox{ on
}\Gamma_\theta\,.
$$
We have
\begin{equation}\label{5.ft4}
\phi_2(\xi ,\eta )=\frac{\lambda }{2}(1+\mathrm{i}\theta)^2t^2-\frac{\mu
}{2}(1-\mathrm{i}\theta)^2s^2+\rho (1+\theta ^2)ts\hbox{ on }\Gamma_\theta\,.
\end{equation}
In particular,
\begin{equation}\label{6.ft4}
\mathrm{Im} \, \phi_2(\xi ,\eta )=\theta (\lambda t^2+\mu s^2),\ \ (\xi ,\eta
)\in \Gamma_\theta 
\end{equation}
and from (\ref{3.ft4}) we conclude that
\begin{equation}\label{7.ft4}
\mathrm{Im} \, \phi_2(\xi ,\eta )=\frac{\theta }{1+\theta ^2}( \lambda |\xi
  |^2+\mu |\eta |^2)\asymp \theta (|\xi |^2+|\eta |^2)\hbox{ on
}\Gamma_\theta\,.
\end{equation}

\par 
We now switch to the variables $\omega ,\widetilde{\omega
}$. Since
\begin{equation}\label{8.ft4}
\xi =\frac{1}{2}(\omega +\widetilde{\omega })\,,\quad \eta
=\frac{1}{2\mathrm{i}}(\omega -\widetilde{\omega })\,,
\end{equation}
we get from (\ref{2.ft4}),
\begin{equation}\label{1.ft5}
\phi_2(\omega ,\widetilde{\omega })=\frac{\lambda }{8}(\omega
+\widetilde{\omega })^2+\frac{\mu }{8}(\widetilde{\omega }-\omega
)^2-\mathrm{i}\frac{\rho}{4}(\omega +\widetilde{\omega })(\omega-\widetilde{\omega
} )\,.
\end{equation}
It follows that
\begin{equation}\label{2.ft5}
\partial_{\tilde{\omega }}\phi_2=\left(\frac{\mu +\lambda
  }{4}+\mathrm{i}\frac{\rho}{2} \right)\widetilde{\omega }-\frac{\mu -\lambda
}{4}\omega\,,
\end{equation}
so $\mathbb{C}\ni \widetilde{\omega }\longmapsto \phi_2(\omega
,\widetilde{\omega })$ has the unique and non-degenerate critical point
\begin{equation}\label{3.ft5}
\widetilde{\omega }_c^0(\omega )=\frac{\mu -\lambda }{\mu +\lambda
  +\mathrm{i} 2\rho} \, \omega \,.
\end{equation}
We have
\begin{equation}\label{4.ft5}
\omega +\widetilde{\omega }_c^0(\omega )=\frac{2(\mu +\mathrm{i}\rho)}{\mu
  +\lambda +\mathrm{i} 2\rho}\omega\,,
\end{equation}
\begin{equation}\label{4,5.ft5}
\omega -\widetilde{\omega }_c^0(\omega )=\frac{2(\lambda  +\mathrm{i}\rho )}{\mu+\lambda +\mathrm{i} 2\rho}\omega \,.
\end{equation}
Using this in (\ref{1.ft5}), we get
$$
  \phi_2(\omega ,\widetilde{\omega }_c^0(\omega ))=
  \frac{\frac{\lambda }{2}(\mu +\mathrm{i} \rho)^2+\frac{\mu }{2}(\lambda +\mathrm{i} \rho)^2-\mathrm{i}\rho (\mu +\mathrm{i}\rho)(\lambda +\mathrm{i}\rho)}
{\Big((\mu +\mathrm{i}\rho)+(\lambda +\mathrm{i}\rho)\Big)^2}\,\omega^2\,,
$$
which simplifies to 
\begin{equation}\label{6.ft5}
\phi_2(\omega ,\widetilde{\omega }_c^0(\omega ))=\frac{\lambda \mu
  +\rho ^2}{2(\lambda +\mu +\mathrm{i} 2\rho)}\, \omega ^2\,.
\end{equation}
When $\rho =0$ this gives
\begin{equation}\label{7.ft5}
\phi_2(\omega ,\widetilde{\omega }_c^0(\omega ))=
\frac{\lambda \mu }{2(\mu +\lambda )}\,\omega ^2\,,
\end{equation}
and when $\lambda =\mu $,
\begin{equation}\label{8.ft5}
\phi_2(\omega ,\widetilde{\omega}_c^0(\omega))=\frac{\lambda +\mu
  -\mathrm{i} 2\rho}{8}\,\omega^2=\frac{\lambda
  -\mathrm{i} \rho}{4}\,\omega^2\,.
\end{equation}

\par By (\ref{2.ft5}),
\begin{equation}\label{8,5.ft5}
\partial_{\tilde{\omega }}^2\phi_2=\frac{\mu +\lambda
}{4}+\mathrm{i}\frac{\rho}{2}\,,
\end{equation}
hence,
\begin{equation}\label{9.ft5}
  \phi_2(\omega ,\widetilde{\omega })=\frac{1}{2}\left( \frac{\mu +\lambda
    }{4}+\mathrm{i}\frac{\rho}{2} \right) (\widetilde{\omega
  }-\widetilde{\omega }_c^0(\omega ))^2+\phi_2(\omega
  ,\widetilde{\omega}_c^0(\omega ))\,.
\end{equation}

\par
We now restrict the attention to $\Gamma_1$. By (\ref{5.ft4}),
we have
\begin{equation*}\begin{split}
&\phi_2(\xi ,\eta )=\mathrm{i}(\lambda t^2+\mu s^2)+2\rho \, ts\,,  \\& 
\hbox{ for}  \,\, (\xi ,\eta
)=((1+\mathrm{i})t,(1-\mathrm{i})s)=(\sqrt{2} \e^{\mathrm{i}\frac{\pi}{4}}t,\sqrt{2}
\e^{-\mathrm{i}\frac{\pi}{4}}s) \in \Gamma_1\,.
\end{split}
\end{equation*}
By (\ref{7.ft4}) we have
\begin{equation}\label{1.ft6}
\mathrm{Im} \,\phi_2(\omega ,\widetilde{\omega })\asymp |\omega
|^2+|\widetilde{\omega }|^2\hbox{ on }\Gamma_1\,.
\end{equation}
The parametrization (\ref{Sto.4}) becomes
\begin{equation}\label{2.ft6}
  \Gamma_1:\ \begin{cases}
    \omega =(1+\mathrm{i})t+(1+\mathrm{i})s=\sqrt{2} \e^{\mathrm{i}\frac{\pi}{4}}(t+s)\\
    \widetilde{\omega }=(1+\mathrm{i})t-(1+\mathrm{i})s=\sqrt{2} \e^{\mathrm{i}\frac{\pi}{4}}(t-s)\,,
  \end{cases}\quad t,s\in \mathbb{R}\,,
\end{equation}
so $\Gamma_1=\e^{\mathrm{i}\frac{\pi}{4}}\mathbb{R}_\omega \times \e^{\mathrm{i}\frac{\pi}{4}}{\mathbb{R}}_{\widetilde{\omega }}$ is a Cartesian product. From the fact
that $\phi_2(\omega ,\widetilde{\omega })$ is a quadratic polynomial in $\widetilde{\omega }$ for
every fixed $\omega \in \mathbb{C}$, with $\mathrm{Im}\, \phi_2(\omega
,\widetilde{\omega })\ge |\widetilde{\omega }|^2/{\mathcal  O}(1)$, when
$\widetilde{\omega}\in \e^{\mathrm{i}\frac{\pi}{4}}\mathbb{R}$,
$|\widetilde{\omega }|\gg |\omega |$ and ``the fundamental lemma''
(see \cite[Lemma 2.1, p. 145]{MeSj75_01}), we deduce that
\begin{equation}\label{3.ft6}
\mathrm{Im}\, \phi_2(\omega ,\widetilde{\omega }_c^0(\omega ))\ge
\inf_{\widetilde{\omega }\in \e^{\mathrm{i}\frac{\pi}{4}}\mathbb{R}}
\mathrm{Im}\, \phi_2(\omega,\widetilde{\omega })\,,
\end{equation}
for every $\omega \in \mathbb{C}$. In particular, by (\ref{1.ft6}),
\begin{equation}\label{4.ft6}\mathrm{Im}\, \phi_2(\omega ,\widetilde{\omega }_c^0(\omega ))\asymp |\omega 
|^2,\quad \omega \in \e^{\mathrm{i}\frac{\pi}{4}}\mathbb{R}\,.
\end{equation}

\par
We now turn to the holomorphic function $\phi$, defined in
$\mathrm{neigh}\big((0,0);\mathbb{C}^2_{\omega ,\tilde{\omega }}\big)$ (or
in $\mathrm{neigh}\big((0,0);\mathbb{C}^2_{\xi ,\eta }\big)$, depending on the choice of coordinates). The function $\phi (\omega ,\cdot )$ has a  non-degenerate critical point $\widetilde{\omega }_c(\omega )$ for $\omega \in \mathrm{neigh}(0;\mathbb{C})$ such that
\begin{equation}\label{5.ft6}
\widetilde{\omega }_c(\omega )=\widetilde{\omega }_c^0(\omega )+{\mathcal O}(\omega ^2)\,,
\end{equation}
so $\widetilde{\omega }_c(\cdot )$ is a holomorphic function with
differential $\widetilde{\omega }^0_c$ at $\omega =0$. As in the
quadratic case, we have
\begin{equation}\label{6.ft6}
\mathrm{Im}\, \phi (\omega ,\widetilde{\omega }_c(\omega ))\ge
\inf_{\widetilde{\omega }\in \e^{\mathrm{i}\frac{\pi}{4}}\mathbb{R}\cap \mathrm{neigh}(0;\mathbb{C})}\mathrm{Im}\, \phi (\omega
,\widetilde{\omega })\,,
\end{equation}
for $\omega \in \mathrm{neigh}(0;\mathbb{C})$,
\begin{equation}\label{1.ft7}
\mathrm{Im}\, \phi (\omega ,\widetilde{\omega }_c(\omega ))\asymp  
|\omega |^2,\ \ \omega \in \mathrm{neigh}(0;\e^{\mathrm{i}\frac{\pi}{4}}\mathbb{R})\,.
\end{equation}

\par Put
\begin{equation}\label{1,5.ft7}
  \psi (\omega )=\phi (\omega ,\widetilde{\omega }_c(\omega )),\ \
  \psi_2 (\omega )=\phi_2 (\omega ,\widetilde{\omega }^0_c(\omega ))
\end{equation}
and notice that
\begin{equation}\label{1,7.ft7}
\psi (\omega )=\psi _2(\omega )+{\mathcal O}(\omega ^3)\,.
\end{equation}

\par 
We now study $\mathrm{I}(\zeta ,1)$ in (\ref{4,5.ft3}). To start with, we
assume only that $a$ is smooth and ${\mathcal  O}(1)$ in
$\mathrm{neigh}((0,0);\Gamma_1)$ and that $\chi \in C_0^\infty $
has its support in a small neighborhood of $(0,0)$. (Similarly we could here weaken the assumption on $\phi $.) The infimum in (\ref{6.ft6}) is attained at a point $\widetilde{\omega}_{\mathrm{inf}}(\omega)\in
\e^{\mathrm{i}\frac{\pi}{4}}\mathbb{R}\cap \mathrm{neigh}(0;\mathbb{C})$ and we have
\begin{equation}\label{2.ft7}
\mathrm{Im}\, \phi (\omega ,\widetilde{\omega })-\mathrm{Im}\, \phi (\omega
,\widetilde{\omega }_{\mathrm{inf}}(\omega ))\asymp |\widetilde{\omega
}-\widetilde{\omega }_{\mathrm{inf}}|^2,\ \ \widetilde{\omega }\in
\e^{\mathrm{i}\frac{\pi}{4}}\mathbb{R}\cap \mathrm{neigh}(0;\mathbb{C})\,.
\end{equation}
It follows that
\begin{equation}\label{3.ft7}
\int_{\e^{\mathrm{i}\frac{\pi}{4}}\mathbb{R}}(a\chi)(\omega ,\widetilde{\omega })\,\e^{\frac{\mathrm{i}}{h}\phi(\omega ,\widetilde{\omega })}\mathrm{d}\widetilde{\omega}={\mathcal O}(h^{\frac{1}{2}})\, \e^{-\frac{1}{h}\mathrm{Im}\, \phi (\omega ,\widetilde{\omega}_{\mathrm{inf}}(\omega ))}\,,
\end{equation}
for $\omega \in \mathrm{neigh}(0;\mathbb{C})$. Observe here that
\begin{equation}\label{4.ft7}
\mathrm{Im}\, \phi (\omega ,\widetilde{\omega }_{\mathrm{inf}(\omega )})\asymp
|\omega |^2,\quad \omega \in \e^{\mathrm{i}\frac{\pi}{4}}\mathbb{R}\cap \mathrm{neigh}(0;\mathbb{C})\,.
\end{equation}

\par We now strengthen the assumption on $a$ and assume that
\begin{equation}\label{5.ft7}
a(\omega ,\widetilde{\omega })={\mathcal  O}(1)\hbox{ is holomorphic in }\mathrm{neigh}\big((0,0);\mathbb{C}^2\big)\,.
\end{equation}
We allow $a$ to depend on $h$ and then assume that $a(\omega
,\widetilde{\omega };h)$ is a holomorphic classical symbol of order
$h^0$.

\par
We assume that $\chi $ is supported in the neighborhood in
(\ref{5.ft7}) and that
\begin{equation}\label{6.ft7}
\chi =1\hbox{ in }\mathrm{neigh}\big((0,0);\mathbb{C}^2_{\omega
  ,\tilde{\omega }}\big)\,.
\end{equation}
The method of steepest descent (and analytic stationary phase, see \cite[Théorème 2.8, p. 16]{Sj82_01}) gives
\begin{equation}\label{7.ft7}
\int_{\e^{\mathrm{i}\frac{\pi}{4}}\mathbb{R}}(a\chi)(\omega ,\widetilde{\omega })\, \e^{\frac{\mathrm{i}}{h}\phi(\omega ,\widetilde{\omega })}\mathrm{d}\widetilde{\omega}=h^{\frac12}\,b(\omega
;h)\, \e^{\frac{\mathrm{i}}{h}\phi(\omega ,\widetilde{\omega}_c(\omega))}\,,
\end{equation}
for $\omega \in \mathrm{neigh}(0;\mathbb{C})$ where $b(\omega ;h)$ is a holomorphic
 classical symbol of order $0$\,,
\begin{equation}\label{8.ft7}
b(\omega ;h)\sim b_0(\omega )+h b_1(\omega )+... \,,
\end{equation}
where
\begin{equation}\label{9.ft7}
  b_0(\omega )=\frac{\sqrt{2\pi}\, \e^{\mathrm{i}\frac{\pi}{4}}}{\sqrt
  {\partial_{\tilde{\omega}}^2\phi (\omega ,\widetilde{\omega
    }_c(\omega ))}}\, a_0(\omega ,\widetilde{\omega }_c(\omega))\,.
\end{equation}
(By (\ref{9.ft5}) we know that $\mathrm{Re} (\partial_{\tilde{\omega
  }}^2\phi)  >0$ and we choose the branch of the square root with the
same property.)
Here $a_0$ is equal to $a$ when $a$ is independent of $h$ and more
generally equal to the leading term in the asymptotic expansion
of $a$ as a holomorphic classical symbol. Further,
\begin{equation}\label{1.ft8}
\partial_{\overline{\omega } }b={\mathcal  O}(\e^{-\frac{1}{Ch}}),\hbox{ where }C>0.
\end{equation}

\par
Since $\Gamma_1=\e^{\mathrm{i}\frac{\pi}{4}}\mathbb{R}\times \e^{\mathrm{i}\frac{\pi}{4}}\mathbb{R}$, we
can apply Fubini's theorem to the integral (\ref{4,5.ft3}) to get
\begin{equation}\label{2.ft8}\begin{split}
  \mathrm{I}(\zeta ,1)&=-\frac{1}{2\pi \mathrm{i}}\int_{\e^{\mathrm{i}\frac{\pi}{4}}{\mathbb{R}}}\frac{1}{\omega-\zeta}\,J(\omega) \, \mathrm{d}\omega\, ,\hbox{ where }\\
    J(\omega )&=\int_{\e^{\mathrm{i}\frac{\pi}{4}}\mathbb{R}}(a\chi)(\omega,\widetilde{\omega})\, \e^{\frac{\mathrm{i}}{h}\phi(\omega
      ,\widetilde{\omega })} \, \mathrm{d}\widetilde{\omega}
    \end{split}
\end{equation}
and then apply (\ref{3.ft7}), (\ref{4.ft7}), (\ref{7.ft7}) to the
inner integral.

We sum up the discussion about $\mathrm{I}(\zeta ,1)$ in
(\ref{4,5.ft3}), (\ref{5.ft3}).
\begin{prop}\label{1ft8}
Assume that $\phi $ is holomorphic in $\mathrm{neigh}((0,0);{\mathbb{C}}^2)$ with a critical point at $(0,0)$ and Taylor expansion as in (\ref{6.ft3}), (\ref{1.ft4}), (\ref{1,5.ft4}), (\ref{2.ft4}) and (\ref{3.ft4}). Then
\begin{equation}\label{3.ft8}
\mathrm{Im}\, \phi (\omega ,\widetilde{\omega })\asymp |\omega
|^2+|\widetilde{\omega }|^2\hbox{ on }\mathrm{neigh}\big((0,0);\Gamma_1\big) 
\end{equation}
and $\phi(\omega,\cdot)$ has a non-degenerate critical point
$\widetilde{\omega }_c(\omega )$ near $0$ for $\omega \in
\mathrm{neigh}(0;\mathbb{C})$ which satisfies (\ref{1.ft7}).

Let $a(\omega ,\widetilde{\omega })$ be holomorphic in
$\mathrm{neigh}((0,0);\mathbb{C}^2)$, either independent of $h$ or a
holomorphic classical symbol of order $h^0$,
$$a(\omega ,\widetilde{\omega })\sim a_0(\omega ,\widetilde{\omega
})+ha_1(\omega ,\widetilde{\omega })+...\,\,.
$$
In view of (\ref{1.ft7}) $\e^{\frac{\mathrm{i}}{h}\psi (\omega )}$ behaves
like a Gaussian along $\e^{\mathrm{i}\frac{\pi}{4}}{\mathbb R}$ with its
peak at $\omega=0$. Recall that $|\zeta |$ is small.

Let $\chi \in C_0^\infty (\mathbb{C}^2)$ be supported in a sufficiently
small neighborhood of $(0,0)$, where the assumptions on $a$ and $\phi
$ are valid. Assume that $\chi =1$ in another small neighborhood of
$(0,0)$. Then $\mathrm{I}(\zeta ,1)$ is given by an iterated integral
as in (\ref{2.ft8}), where the inner integral $J$ satisfies
\begin{equation}\label{4.ft8}
\partial_{\overline{\omega }}J={\mathcal O}(\e^{-\frac{1}{{\mathcal  O}(h)}}),\quad
\omega \in \mathrm{neigh}(\e^{\mathrm{i}\frac{\pi}{4}}\mathbb{R};\mathbb{C})\,.
\end{equation}
Moreover there is a neighborhood $\Omega $ of $0$ in ${\mathbb{C}}$, where we
have (\ref{7.ft7}), (\ref{8.ft7}), (\ref{9.ft7}), (\ref{1.ft8}). If
$V\Subset \Omega $ is another neighborhood of $0$, then
$$
J={\mathcal  O}(\e^{-\frac{1}{{\mathcal  O}(h)}})\hbox{ in }\,\,\mathrm{neigh}\big((\e^{\mathrm{i}\frac{\pi}{4}}\mathbb{R})\setminus V;\mathbb{C}\big)\,.
$$
\end{prop}

\par 
We next look at the asymptotics of $\mathrm{I}(\zeta ,1)$ in
(\ref{2.ft8}). Inserting (\ref{7.ft7}) there gives
\begin{equation}\label{1.ft9}
\mathrm{I}(\zeta ,1)=-\frac{1}{2\pi \mathrm{i}}\int_{\e^{\mathrm{i}\frac{\pi}{4}}\mathbb{R}}
\frac{1}{\omega -\zeta }h^{\frac12}b(\omega ;h)\, \e^{\frac{\mathrm{i}}{h}\psi (\omega)} \, \mathrm{d}\omega\,,
\end{equation}
where $b$ satisfies (\ref{8.ft7})--(\ref{1.ft8}) and we have defined $\psi
$ in (\ref{1,5.ft7}).

\par For $\zeta \in {\mathbb{C}}$ small and outside a conic neighborhood of
$\e^{\mathrm{i}\frac{\pi}{4}}{\mathbb{R}}$, we have
$$
|\mathrm{I}(\zeta ,1)|\le {\mathcal O}(1)\int_0^{+\infty
}\frac{h^{\frac{1}{2}}}{|\zeta |+h^{1/2}+t} \, \e^{-\frac{t^2}{h}}\, \mathrm{d}t\,, 
$$
where we also used that we can deform $\e^{\mathrm{i}\frac{\pi}{4}}{\mathbb{R}}$ in a
$h^{1/2}$ neighborhood of $0$ in a such a way that $|\zeta -\omega
|\ge \frac{h^{1/2}}{{\mathcal O}(1)}$ along the deformed contour. Put
$t=\sqrt{h}\, s$:
\[
|\mathrm{I}(\zeta ,1)|={\mathcal
  O}(h^{\frac{1}{2}})\int_0^{+\infty }\frac{1}{\lambda
  +s} \, \e^{-s^2}\, \mathrm{d}s, \, \hbox{ where } \, \lambda :=\frac{|\zeta |}{\sqrt{h}}+1\ge 1\,.
\]
Here the integral is
$$
\le \int_0^{+\infty }\frac{1}{\lambda } \, \e^{-s^2}\, \mathrm{d}s={\mathcal
  O}(1)/\lambda 
$$ 
and also
$$
\ge \int_0^1 \frac{1}{\lambda +s} \, \e^{-s^2}\, \mathrm{d}s \ge \frac{1}{{\mathcal
    O}(\lambda )} \,. 
    $$
We retain that
\begin{equation}\label{2.ft9}
|\mathrm{I}(\zeta ,1) |=\frac{h^{1/2}{\mathcal O}(1)}{\lambda }=\frac{{\mathcal
    O}(h)}{|\zeta |+h^{1/2}} \, .
\end{equation}
 
\par 
Let $\chi \in C_0^\infty (D(0, 1))$ be $=1$ in $D(0,1/2)$. Let
$0<\delta\ll 1
 $ be a constant. We shall study
\begin{equation}\label{3.ft9}
\mathrm{I}_{\mathrm{new}}(\zeta )=-\frac{1}{2\pi
  \mathrm{i}}\int_{\e^{\mathrm{i}\frac{\pi}{4}}{\mathbb
    R}}\frac{h^{\frac12}}{\omega -\zeta} \, b(\omega ;h)\chi
(h^{\delta -1/2}\omega) \, \e^{\frac{\mathrm{i}}{h}\psi (\omega )} \, \mathrm{d}\omega\,,
\end{equation}
where we restrict $\zeta $ to $D(0,h^{-\delta +1/2}/4)$ from now on. 
Observe that
$$
\mathrm{I}_{\mathrm{new}}(\zeta )-\mathrm{I}(\zeta ,1)={\mathcal O}(\e^{-h^{-2\delta
  }/{\mathcal O}(1)}),\hbox{ for }\zeta \in D\Big(0,\frac{h^{\frac{1}{2}-\delta
  }}{4}\Big)\,.
$$

\par 
We make the dilation $\omega =h^{\frac12}\widehat{\omega}$ and
define $\widehat{\zeta }$ by $\zeta =h^{\frac12}\widehat{\zeta }$. Put
\begin{equation}\label{4.ft9}
\widehat{\psi }(\widehat{\omega })=\frac{1}{h}\psi
(h^{\frac12}\widehat{\omega })\,.
\end{equation}
We have the Taylor series expansion, with $\psi_2=\varphi_2$,
\begin{equation}\label{5.ft9}
\psi (\omega )=\psi _2(\omega )+r(\omega )\sim \psi _2(\omega )+\psi
_3(\omega )+\psi _4(\omega )+\ldots\,,
\end{equation}
where $\psi _k(\omega )=\dfrac{\psi ^{(k)}(0)}{k!}\omega^k$ is homogeneous of
degree $k$ and $r(\omega )={\mathcal O}(\omega ^3)$.
It follows that
\begin{equation}\label{6.ft9}
\widehat{\psi }(\widehat{\omega };h)=\psi _2(\widehat{\omega
})+h^{\frac12}\widehat{r}(\widehat{\omega };h),\quad
h^{\frac{1}{2}}\widehat{r}(\widehat{\omega
};h)=\frac{1}{h}r(h^{\frac12}\widehat{\omega })
\end{equation}
for $\widehat{\omega }={\mathcal O}(h^{-\frac{1}{2}})$, where
\begin{equation}\label{6,5.ft9}
h^{\frac{1}{2}}\widehat{r}(\widehat{\omega };h)= \sum_{j=3}^{+\infty} h^{-1+\frac{j}{2}}\psi_j(\widehat{\omega }) \,.
\end{equation}
Since $\psi_j$ is homogeneous of degree $j$, this can be viewed as an asymptotic formula in $D(0,h^{-\delta})$. More precisely, put
$$
h^{\frac{1}{2}}\widehat{r}^{(N)}(\widehat{\omega
})=\sum_{j=3}^{N-1}h^{-1+\frac{j}{2}}\psi _j(\widehat{\omega })\,. 
$$
Then
\begin{equation}\label{7.ft9}
  h^{\frac{1}{2}}\widehat{r}(\widehat{\omega })-
h^{\frac{1}{2}}\widehat{r}^{(N)}(\widehat{\omega })
  ={\mathcal O}\big(h^{-1+N(\frac{1}{2}-\delta
  )}\big)\hbox{ in
} D(0,h^{-\delta} )\,.
\end{equation}
In particular,
\begin{equation}\label{1.nyft10}
  h^{\frac{1}{2}}\widehat{r}(\widehat{\omega }),\
  h^{\frac{1}{2}}\widehat{r}^{N}(\widehat{\omega })
  ={\mathcal O}(1)h^{3(\frac{1}{2}-\delta )-1}={\mathcal
    O}(1)h^{\frac{1}{2}-3\delta }=o(1)\,,
\end{equation}
where from now on we assume that $\delta <1/6$.

\par 
Put $\exp^{(N)}(t)=\overset{N-1}{\underset{k=0}{\sum}} \frac{t^k}{k!}$. When $t={\mathcal O}(1)$ 
we have 
$$
\exp^{(N)}(t)-\exp t={\mathcal O}(1)t^N\,.
$$
We estimate
\begin{align*}
\exp^{(N)}\big(h^{\frac{1}{2}}\widehat{r}^{(N)}(\widehat{\omega })\big)-\exp
\big(h^{\frac{1}{2}}\widehat{r}(\widehat{\omega })\big)
=&\left(\exp^{(N)}\big(h^{\frac{1}{2}}\widehat{r}^{(N)}(\widehat{\omega })\big)-\exp
\big(h^{\frac{1}{2}}\widehat{r}^{(N)}(\widehat{\omega })\big)\right)\\
&+\left(\exp \big(h^{\frac{1}{2}}\widehat{r}^{(N)}(\widehat{\omega })\big)-\exp
\big(h^{\frac{1}{2}}\widehat{r}(\widehat{\omega })\big)  \right)\\
=&:A+B\,.
\end{align*}
Here for $\widehat{\omega }\in D(0,h^{-\delta })$,
$$
A={\mathcal
  O}(1)\left(h^{\frac{1}{2}}\widehat{r}^{(N)}(\widehat{\omega })
\right)^N={\mathcal O}(1)\left( h^{(\frac{1}{2}-3\delta )} \right)^N
={\mathcal O}(1) h^{(\frac{1}{2}-3\delta )N}\, ,
$$
$$
B=
{\mathcal O}(1)\left(h^{\frac{1}{2}}\widehat{r}^{(N)}(\widehat{\omega })-
h^{\frac{1}{2}}\widehat{r}(\widehat{\omega })
\right)={\mathcal O}(1)h^{N(\frac{1}{2}-\delta )-1}\,,
$$
hence
\begin{equation}\label{2.nyft10}
\exp^{(N)}(h^{\frac{1}{2}}\widehat{r}^{(N)}(\widehat{\omega }))-\exp
(h^{\frac{1}{2}}\widehat{r}(\widehat{\omega }))={\mathcal
  O}(1)h^{(\frac{1}{2}-3\delta )N-1}\hbox{ in }D(0,h^{-\delta })\,.
\end{equation}

\par 
We next expand
\begin{align*}
\exp^{(N)}\left(h^{\frac{1}{2}}\widehat{r}^{(N)}(\widehat{\omega })
\right)& =\sum_0^{N-1}\frac{1}{k!}\left(h^{\frac{1}{2}}\widehat{r}^{(N)}(\widehat{\omega
  }) \right)^k
=\sum_0^{N-1}\frac{1}{k!}\left( \sum_{j=3}^{N-1}\frac{1}{h}\psi
  _j(h^{\frac{1}{2}}\widehat{\omega }) \right)^k\\
&=\sum_{k=0}^{N-1}\sum_{\hskip10pt  3\le j_1,...,j_k\le
  N-1}\widetilde{c}_{k,j_1,..,j_k}\left(\frac{1}{h}\psi _{j_1}(\widehat{\omega })
\right)...
\left(\frac{1}{h}\psi _{j_k}(\widehat{\omega }) \right)\\
&=\sum_{k=0}^{N-1}\sum_{\hskip10pt 3\le j_1,...,j_k\le
  N-1}c_{k,j_1,..,j_k}h^{-k+\frac{j_1}{2}+...+\frac{j_k}{2}}\widehat{\omega }^{j_1+...+j_k}\,.
\end{align*}
With $l =j_1+...+j_k\ge 3k$, we get with new coefficients
$c_{k,l}$:
\begin{equation}\label{3.nyft10}
  \exp^{(N)}\left(h^{1/2}\widehat{r}^{(N)}(\widehat{\omega })) \right)
  =
  \sum_{k\ge 0}\sum_{l=3k}^{3(N-1)}c_{k,l}h^{-k+l/2}\widehat{\omega }^l\,.
\end{equation}
With $-k+l/2=:j/2$, $j\in {\mathbb Z}$ we get
$$
\frac{l }{2}=\frac{j}{2}+k\le \frac{j}{2}+\frac{l}{3} \Longleftrightarrow 
\frac{l}{6}\le \frac{j}{2}\,,
$$
then
\begin{equation}\label{4.nyft10}
l \le 3j\,.
\end{equation}
Hence,
\begin{equation}\label{1.nyft11}
\exp^{(N)}\left(h^{\frac{1}{2}}\widehat{r}^{(N)}(\widehat{\omega })
\right) =\sum_{j\ge 0}h^{\frac{j}{2}}p_j(\widehat{\omega })\, ,
\end{equation}
where $p_j$ is a polynomial of degree $\le 3j$ and the sum is
finite. The $p_j$ also depend on $N$; $p_j=p_j^{(N)}$, but we can
apply (\ref{2.nyft10}) with $\delta =0$, to see that
\begin{equation}\label{2.nyft11}
  p_j^{(N)}(\widehat{\omega })=p_j^{(M)}(\widehat{\omega }),\hbox{
    when }M\ge N
\end{equation}
and $h^{j/2}>h^{-1+N/2},$ i.e. when $j<N-2$. In other words,
$p_j^{N}=p_j$ is independent of $N$ when $N>j+2$.

It follows that
\begin{equation}\label{2.ft10}
\e^{\mathrm{i} h^{\frac12} \widehat r(\widehat{\omega };h)}\sim
1+h^{\frac12}p_1(\widehat{\omega})+hp
_2(\widehat{\omega})+\ldots \hbox{ in }D(0,h^{-\delta }),
\end{equation}
where $p _j$ is a polynomial of degree $\le 3j$.

\par Similarly, define
\begin{equation}\label{10.ft9}
\widehat{b}(\widehat{\omega };h)=b(h^{\frac12}\widehat{\omega };h) 
\end{equation}

\par\noindent 
then $\mathrm{I}_{\mathrm{new}}(\zeta ;h)$ becomes
\begin{equation}\label{1.ft10}
\widehat{\mathrm{I}}(\widehat{\zeta };h)=-\frac{h^{\frac12}}{2\pi
  \mathrm{i}}\int_{\e^{\mathrm{i}\frac{\pi}{4}}{\mathbb R}}\, \frac{\chi
(h^\delta\widehat{\omega })}{\widehat{\omega}-\widehat{\zeta}} \, \widehat{b}(\widehat{\omega
};h) \e^{\mathrm{i}\widehat{\psi}(\widehat{\omega};h)} \mathrm{d}\widehat{\omega }\,.
\end{equation}

From the Taylor expansion of $b$ in (\ref{10.ft9}), we get
$$
\widehat{b}(\widehat{\omega };h)\sim
\widehat{b}_0(\widehat{\omega
};h)+h^{\frac{1}{2}}\widehat{b}_1(\widehat{\omega };h)+...
\hbox{ in }D(0,h^{-\delta })\,,
$$
where $b_j$ is a polynomial of degree $\le j$,
hence using also (\ref{2.ft10})
\begin{equation}\label{3.ft10}
\widehat{b} (\widehat{\omega};h) \, \e^{\mathrm{i}\widehat{\psi}(\widehat{\omega};h)}=d(\widehat{\omega};h) \, \e^{\mathrm{i}\psi _2(\widehat{\omega})} \,,
\end{equation}
where
$$
d(\widehat{\omega };h)\sim d_0(\widehat{\omega
})+h^{\frac12}d_1(\widehat{\omega })+\ldots \, \hbox{ in } \, D(0,h^{-\delta})\,,
$$
and $d_j(\widehat{\omega })$ is a computable polynomial of degree
$\le 3j$.

\par 
Recall that $\psi _2(\omega )$ is given by (\ref{1,5.ft7}), (\ref{6.ft5}), hence
$\psi_2(\widehat{\omega })=f\dfrac{\widehat{\omega }^2}{2}$,
\begin{equation}\label{4,5.ft10}
f=\frac{\lambda \mu +\rho ^2}{\lambda +\mu +2i\rho}\,.
\end{equation}
 We have
\begin{equation}\label{5.ft10}
  \widehat{\omega }=\frac{1}{f}\partial _{\widehat{\omega }}\psi
  _2(\widehat{\omega })\,,
\end{equation}
and asymptotically,
\begin{equation}\label{6.ft10}
  d(\widehat{\omega };h) \e^{\mathrm{i}\psi_2(\widehat{\omega })}
  =
  c(f^{-1}D_{\hat{\omega }};h)\left( \e^{\mathrm{i}\psi
      _2(\widehat{\omega })} \right),
\end{equation}
where $c$ has the same properties and the same leading term as $d$,
\begin{equation}\label{4.ft10}
c(\widehat{\omega };h)\sim c_0(\widehat{\omega
})+h^{\frac12}c_1(\widehat{\omega })+\ldots \, \hbox{ in } \, D(0,h^{-\delta})\,, 
\end{equation}
and $c_j(\widehat{\omega })$ is a computable polynomial of degree
$\le 3j$.

\par In fact,
$$
\widehat{\omega }^k \e^{\mathrm{i}\psi _2(\widehat{\omega })}=
(f^{-1}D_{\widehat{\omega }})^k \e^{\mathrm{i}\psi _2(\widehat{\omega
  })}+g_k(\widehat{\omega }) \, \e^{\mathrm{i}\psi _2(\widehat{\omega })}\,,
$$
where $g_k$ is a polynomial of degree $\le k-1$ and by recursion, we
find a polynomial $c_k$ of the same degree as $d_k$ such that
$$
c_k(f^{-1}D_{\widehat{\omega }})(\e^{\mathrm{i}\psi _2(\widehat{\omega
  })})=d_k(\widehat{\omega }) \, \e^{\mathrm{i}\psi _2(\widehat{\omega })}\,.
$$

\par Now use (\ref{1.ft10}), (\ref{3.ft10}), (\ref{6.ft10}), to get
\begin{equation}\label{7.ft10}
  \begin{split}
\widehat{\mathrm{I}}(\widehat{\zeta };h)&=-\frac{h^{\frac12}}{2\pi
  \mathrm{i}}\int_{\e^{\mathrm{i}\frac{\pi}{4}}{\mathbb R}} \frac{\chi (h^\delta
  \widehat{\omega })}{\widehat{\omega} - \widehat{\zeta
  }} \, d(\widehat{\omega};h) \, \e^{\mathrm{i}\psi
  _2(\widehat{\omega })} \, \mathrm{d}\widehat{\omega}\\
&=-\frac{h^{\frac12}}{2\pi
  \mathrm{i}}\int_{\e^{\mathrm{i}\frac{\pi}{4}}{\mathbb R}} \frac{\chi (h^\delta
  \widehat{\omega })}{\widehat{\omega } - \widehat{\zeta
  }} \, c(f^{-1}D_{\widehat{\omega }};h)\left( \e^{\mathrm{i}\psi_2(\widehat{\omega })}\right) \, \mathrm{d}\widehat{\omega }\\
&=h^{\frac12}c(f^{-1}D_{\widehat{\zeta
  }};h)\left(G_{\mathrm{i}\psi_2}(\widehat{\zeta}) \right)
+{\mathcal O}(h^\infty ), \, \hbox{ in }D(0,h^{-\delta }/4)\,,
  \end{split}
\end{equation}
where
\begin{equation}\label{8.ft10}
  G_{\mathrm{i}\psi_2}(\widehat{\zeta})=
  \frac{1}{2\pi \mathrm{i}}\int_{\e^{\mathrm{i}\frac{\pi}{4}}{\mathbb R}} 
  \frac{1}{\widehat{\zeta
  }-\widehat{\omega}}\left( \e^{\mathrm{i}\psi
    _2(\widehat{\omega })}\right) \, \mathrm{d}\widehat{\omega}, \, \hbox{ in }D(0,h^{-\delta }/4)\,.
\end{equation}
The ${\mathcal O}(h^\infty )$ remainder comes from the cutoff and from the asymptotic nature of $c$. We also used (when necessary) that the
contour $\e^{\mathrm{i}\frac{\pi}{4}}{\mathbb{R}}$ can be deformed in a fixed neighborhood
of $0$.
Here we recall that
\begin{equation}\label{9.ft10}
\widehat{\zeta }\ne 0\hbox{ and }
\widehat{\zeta } \, \hbox{ is not in a conic neighborhood of }
\e^{\mathrm{i}\frac{\pi}{4}}\Big( {\mathbb R}\setminus \{ 0 \} \Big)\,.
\end{equation}
\par 
When $\mathrm{i}\psi_2$ is replaced by $-\frac{1}{2}\widehat{\omega }^2$ and
$\e^{\mathrm{i}\frac{\pi}{4}}{\mathbb R}$ by ${\mathbb R}$, we have, modulo the multiplicative factor $2\pi \mathrm{i}$, the special functions $G_{l/r}$ discussed in \cite{KlSjSt23_02} 
and as there we are in the presence of two special functions 
$G_{\mathrm{i}\psi_2}^{l/r}$ where $G_{\mathrm{i}\psi_2}^l$ is
defined as in (\ref{8.ft10}) after possibly deforming  the contour
$\e^{\mathrm{i}\frac{\pi}{4}}{\mathbb R}$ slightly so that the deformed contour avoids $\widehat{\zeta }$ and passes with this point to the left when oriented in the direction of increasing $\mathrm{Re} \,\widehat{\omega}$. Likewise  
$G_{\mathrm{i}\psi_2}^r$ is defined, with the deformed contour
passing with $\widehat{\zeta }$ to the right. See Appendix \ref{spe}.

\par\smallskip 
When the contour is equal to $\e^{\mathrm{i}\frac{\pi}{4}}{\mathbb R}$ and
$\widehat{\zeta}$ satisfies (\ref{9.ft10}), we get
$G_{\mathrm{i}\psi_2}^l(\widehat\zeta)$ when $\mathrm{Im} (\e^{-\mathrm{i}\frac{\pi}{4}}\widehat{\zeta })>0$ and
$G_{\mathrm{i}\psi_2}^r(\widehat\zeta)$ when $\mathrm{Im} (\e^{-\mathrm{i}\frac{\pi}{4}}\widehat{\zeta })<0$.

\par\smallskip 
By the residue theorem we have for $\widehat{\zeta }\in D(0,h^{-\delta
}/4)$,
\begin{equation}\label{1.ft11}
G^r_{\mathrm{i}\psi_2}(\widehat{\zeta})-G^l_{\mathrm{i}\psi_2}(\widehat{\zeta })= \,\e^{\mathrm{i}\psi_2(\widehat{\zeta} )}\,.
\end{equation}

\noindent
Similarly we can define $\widehat{I}^{l,r}(\widehat{\zeta};h)$
directly from (\ref{1.ft10}) and then
\begin{equation}\label{2.ft11}
  \widehat{\mathrm{I}}^r(\widehat{\zeta})-\widehat{\mathrm{I}}^l(\widehat{\zeta})
=h^{\frac12} \, \widehat{b}(\widehat{\zeta};h) \, 
\e^{\mathrm{i}\widehat{\psi}(\widehat{\zeta};h)} \hbox{ in }D(0,h^{-\delta }/4)\,.
\end{equation}

\par 
From (\ref{7.ft10}), (\ref{8.ft10}), we get
under the assumption (\ref{9.ft10})
\begin{equation}\label{3.ft11}
\widehat{\mathrm{I}}(\widehat{\zeta };h)=h^{\frac12} \,
c(f^{-1}D_{\widehat{\zeta}};h)\Big(G_{\mathrm{i}\psi_2}(\widehat{\zeta})\Big)+\mathcal{O}(h^\infty) \hbox{
in }D(0,h^{-\delta }/4)\,,
\end{equation}
choosing the branch ``$l$'' for $\widehat{\mathrm{I}}$ and $G$ when
$\mathrm{Im} \big(\e^{-\mathrm{i}\frac{\pi}{4}}\widehat{\zeta}\big)>0$ and the branch ``$r$'' when
$\mathrm{Im} \big(\e^{-\mathrm{i}\frac{\pi}{4}}\widehat{\zeta}\big)<0$.

\par We rewrite this expansion in terms of the Dawson-type functions
$G^l$, $G^r$ (more precisely these are the Hilbert transforms of the 
Dawson integral, see e.g., \cite{AbSt84_01}) and start with a general remark (dropping the hats): Let
\begin{equation}\label{1.ft12}
\rho (\omega )=r\frac{\omega ^2}{2}\,,\quad  0\ne r\in {\mathbb{C}}\,.
\end{equation}
Then $\{\omega \in {\mathbb{C}\setminus\{0\}};\, \mathrm{Re} \rho (\omega )<0 \}$ is the
disjoint union of two sectors $R_-$ and $R_+$. Put
\begin{equation}\label{2.ft12}
G^{l / r}_\rho (\zeta )=\frac{1}{2\pi \mathrm{i}}\int_\gamma \frac{1}{\zeta
  -\omega } \e^{\rho (\omega )}\, \mathrm{d}\omega\,,
\end{equation}
where $\gamma =a+b{\mathbb{R}}$, $b\ne 0$ is a straight line, oriented by
the natural orientation of ${\mathbb{R}}$, such that for $\gamma (t)=a+bt$,
\begin{equation}\label{3.ft12}
\gamma (t)\in R_\pm \hbox{ and }\mathrm{dist}(\gamma (t),{\mathbb
  C}\setminus R_\pm)\asymp |t|,\hbox{ when }\pm t\gg 1\,,
\end{equation}
\begin{equation}\label{4.ft12}\begin{split}
  &\gamma (t) \hbox{ avoids }\zeta \hbox{ and passes with that point to
  the}\\ &\hbox{ left/right, when }t\hbox{ increases from }-\infty \hbox{ to
}+\infty .
\end{split}
\end{equation}
The value of $G_\rho ^{l/r}(\zeta )$ only depends on the homotopy
class of $\gamma $ with these properties, so $G_{\rho }^{l/r}$ are
entire functions related via the the Cauchy formula,
\begin{equation}\label{5.ft13}
G_\rho ^l(\zeta )-G_\rho ^r(\zeta )=-\e^{\rho (\zeta )}\,.
\end{equation}
These definitions extend those of $G^{l/r}=G^{l/r}_{-\zeta
  ^2/2}$, $G^{l/r}_{\mathrm{i} \psi _2}$ above, if we adopt the natural
orientations of ${\mathbb{R}}$ and $\e^{\mathrm{i}\frac{\pi}{4}}{\mathbb{R}} $ respectively.

\par Let
\begin{equation}\label{1.ft13}
\widetilde{\rho }(\widetilde{\omega })=\widetilde{r} \, \frac{\widetilde{\omega
}^2}{2}\,,\quad 0\ne \widetilde{r}\in {\mathbb{C}}\,.
\end{equation}
Then $\{\widetilde{\omega }\in \dot{{\mathbb{C}}};\, \mathrm{Re} \widetilde{\rho
}(\widetilde{\omega })<0 \}$ is the disjoint union of two sectors
$\widetilde{R}_-$ and $\widetilde{R}_+$. Define $G^{l
  /r}_{\widetilde{\rho }}$ as in (\ref{2.ft12}). We have
\begin{equation}\label{2.ft13}
\widetilde{\rho }(\widetilde{\omega })=\rho (\omega )\hbox{ when
}\widetilde{\omega }=(r/\widetilde{r })^{1/2}\omega ,
\end{equation}
and we choose the branch of the square root for which
$\widetilde{R}_\pm=(r/\widetilde{r})^{1/2}R_\pm$. If $\widetilde{\zeta
}=(r/\widetilde{r})^{1/2}\zeta $, we get
\begin{equation}\label{3.ft13}
  G^{l/r}_{\widetilde{\rho }}(\widetilde{\zeta })=
\frac{1}{2\pi
  \mathrm{i}}\int _{\widetilde{\gamma }} \frac{1}{\widetilde{\zeta}
  -\widetilde{\omega} } \,
\e^{\widetilde{\rho} (\widetilde{\omega } )}\, \mathrm{d} \widetilde{\omega }
=
  \frac{1}{2\pi
    \mathrm{i}}\int _\gamma \frac{1}{\zeta -\omega } \, \e^{\rho (\omega )}\, \mathrm{d} \omega
  =G_\rho (\zeta )\,.
\end{equation}

\par\noindent
We apply this to the case $\widehat{\rho }(\widehat{\omega
})=\mathrm{i}\psi _2(\widehat{\omega })$ (replacing tildes by hats) $\rho
(\omega )=~-\frac{\omega ^2}{2}$. We have $\mathrm{i}\psi _2(\widehat{\omega
})=\mathrm{i} f\frac{\widehat{\omega} ^2}{2}=:\widehat{r} \, \frac{\widehat{\omega }^2}{2}$ 
with $f$ as in (\ref{4,5.ft10}), $\rho (\omega )=r \, \frac{\omega ^2}{2}$, $r=-1$. The
change of variables in (\ref{2.ft13}) is
\begin{equation}\label{4.ft13}
\widehat{\omega }=(r/\widehat{r})^{1/2}\omega =(\mathrm{i} /f)^{1/2}\omega\,.
\end{equation}
We have $\mathrm{Im} (\mathrm{i}/f)>0$ and choose the branch of $(\mathrm{i}/f)^{1/2}$ in
the first quadrant, $0<\arg \big((\mathrm{i}/f)^{1/2}\big) < \pi /2$. Then
(\ref{3.ft13}) gives
\begin{equation}\label{5'.ft13}
G_{\mathrm{i}\psi _2}^{l / r}(\widehat{\zeta })=G_{-\frac{\omega^2}{2}}^{l /
  r}(\zeta ), \hbox{ when } \, \widehat{\zeta }=(\mathrm{i} /f)^{1/2}\zeta \,. 
\end{equation}

\par
We now make the same change of variable in (\ref{6.ft10}). From (\ref{4.ft13}) we get
\begin{equation}\label{1.ft14}
D_{\widehat{\omega }}=(f/\mathrm{i} )^{1/2}D_\omega =\mathrm{i} (f/\mathrm{i}
)^{1/2}\partial _\omega = (\mathrm{i} f)^{1/2}\partial _\omega \,,
\end{equation}
and (\ref{6.ft10}) gives
\begin{equation}\label{2.ft14}
  d(\widehat{\omega };h)\e^{\mathrm{i}\psi _2(\widehat{\omega})}=
  c\big((\mathrm{i}/f)^{1/2}\partial _\omega ;h\big)(\e^{-\frac{\omega^2}{2}})\,.
\end{equation}
Substitution in (\ref{7.ft10}) gives
\begin{equation}\label{3.ft14}
\hskip-15pt  \widehat{\mathrm{I}}^{l/r}(\widehat{\zeta };h)=
  h^{1/2}c\big((\frac{\mathrm{i}}{f})^{1/2}\partial _{\zeta^\dagger };h\big)\Big(G_{-\frac{(\zeta
    ^\dagger)^2}{2}}^{l/r}(\zeta^\dagger )\Big)+{\mathcal O}(h^\infty ),\hbox{ in
  }D_{\widehat{\zeta }}(0,\frac{h^{-\delta }}{4})\,,
\end{equation}
where
\begin{equation}\label{4.ft14}
\widehat{\zeta }=(\mathrm{i}/f)^{1/2}\zeta^\dagger \, .
\end{equation}

\begin{prop}\label{1ft15} Recall the expression (\ref{1.ft9}) for
  $\mathrm{I}(\zeta ,1)$,
where $b$, $\psi $ fulfill (\ref{1,5.ft7}), (\ref{8.ft7})--(\ref{1.ft8}). 
Let $\psi _2(\omega )=f\frac{\omega ^2}{2}$ be the leading
term in the Taylor expansion (\ref{5.ft9}) of $\psi $. Here $f$ is
given in (\ref{4,5.ft10}).

\par
We make the dilation $\omega =h^{1/2}\widehat{\omega }$ and
define $\widehat{\zeta }$ by $\zeta =h^{1/2}\widehat{\zeta }$. Define
$\widehat{\psi }$ as in (\ref{4.ft9}). Let $0\le \delta <1/6 $. Then
$$
\widehat{\psi }(\widehat{\omega })=\frac{1}{h}\psi
(h^{\frac{1}{2}}\widehat{\omega })=\psi _2(\widehat{\omega })+{\mathcal
  O}(h^{\frac{1}{2}-3\delta })\hbox{ in }D(0,h^{-\delta })\,.
$$
Put $\widehat{b}(\widehat{\omega };h)=b(h^{1/2}\widehat{\omega };h)$
and define $\widehat{\mathrm{I}}(\widehat{\zeta };h)$ as in (\ref{1.ft10}). Then
\begin{equation}\label{1.ft15}
\widehat{\mathrm{I}}(\widehat{\zeta };h)-\mathrm{I}(\zeta ;h)={\mathcal
  O}\big(\e^{-h^{-2\delta }/{\mathcal O}(1)}\big), \quad \zeta \in
D(0,\frac{1}{4}h^{\frac{1}{2}-\delta })\,.
\end{equation}
We have asymptotically in $D(0,h^{-\delta })$,
\begin{equation}\label{2.ft15}
\widehat{b}(\widehat{\omega };h) \e^{\mathrm{i}\widehat{\psi }(\widehat{\omega
  };h)}=c(f^{-1}D_{\widehat{\omega }};h)\left( \e^{\mathrm{i}\psi
    _2(\widehat{\omega })} \right)\,,
\end{equation}
where $\widehat{c}(\widehat{\omega };h)\sim c_0(\widehat{\omega
})+h^{1/2}c_1(\widehat{\omega })+... \, $ in $D(0,h^{-\delta })$ and $c_j(\widehat{\omega })$ is a polynomial of degree
$\le 3j$ (as in
(\ref{4.ft10})).

\par 
Uniformly for $\widehat{\zeta }\in D(0,h^{-\delta }/4)$
satisfying (\ref{9.ft10}),
\begin{equation}\label{3.ft15}\begin{split}
\mathrm{I}(\zeta,1)  &=h^{\frac12}c(f^{-1}D_{\widehat{\zeta}};h)
\left(G_{\mathrm{i}\psi_2}(\widehat{\zeta}) \right)
+{\mathcal O}(h^\infty)\\
&=h^{1/2}c((\mathrm{i}/f)^{1/2}\partial _{\zeta^\dagger} ;h)
\Big(G^{l/r}_{-\frac{\omega^2}{2}}
(\zeta^\dagger )\Big)+{\mathcal O}(h^\infty )\,,
  \end{split}
\end{equation}
where
$\widehat{\zeta }=(i/f)^{1/2}\zeta ^\dagger$.
The functions $G_\rho =G_\rho ^{l/r}$ are defined by line integrals in
(\ref{8.ft10}), (\ref{2.ft12}), and we take $G_{..}=G_{..}^l$ when 
$\mathrm{Im}(\e^{-\mathrm{i}\frac{\pi}{4}}\widehat{\zeta })>0$ and
$G_{..}=G_{..}^r$ when $\mathrm{Im}(\e^{-\mathrm{i}\frac{\pi}{4}}\widehat{\zeta })<0$.
\end{prop}

\section{The term \texorpdfstring{$\mathrm{II}(\zeta ,1-\delta)$}{} in (\ref{2.sto8})}\label{II}
\setcounter{equation}{0}

In this section we discuss the second term in the 
decomposition of Section \ref{Sto} for various values of the 
parameter $\varepsilon:=|\zeta|/\sqrt{h}$.

As in Proposition \ref{1sto7}, let $a$, $\phi $ be holomorphic on some
neighborhood of $(0,0)\in \mathbb{C}^2_{\omega ,\tilde{\omega }}$ and
let $\chi \in C_0^\infty (\mathbb{C}^2_{\omega ,\tilde{\omega}})$
have its support in that neighborhood. More assumptions will be made below.

Recall that in (\ref{2.sto8}),
\begin{equation}\label{II.7}
 \mathrm{II}(\zeta ,1-\delta)
=-\int_0^{1-\delta} \hskip-5pt a(\zeta,\widetilde \zeta;h)\chi(\zeta,\widetilde\zeta)\, \e^{\frac{\mathrm{i}}{h}\phi(\zeta,\widetilde
   \zeta)}\,\dfrac{2\mathrm{i}\,\overline{\widetilde{\zeta}}}{1-\theta^2}\, \mathrm{d}\theta\,,
\end{equation}
where $(\zeta,\widetilde\zeta)\in \Gamma_\theta$ so $\widetilde\zeta$ is given by 
\eqref{Sto.26}, namely:
\begin{equation*}
 \widetilde \zeta(\theta)=\frac{2\mathrm{i}\theta}{1-\theta^2} \,\zeta+\frac{1+\theta^2}{1-\theta^2} \,\overline{\zeta}\,.
\end{equation*}

\noindent
Let us make a first change of variable:
\begin{equation*}
 \e^{\mathrm{i}\alpha}:=\frac{1+\mathrm{i}\theta}{1-\mathrm{i}\theta}=\frac{1-\theta^2}{1+\theta^2}+\mathrm{i}\frac{2\theta}{1+\theta^2}\,,\qquad \mathrm{i.e.}
 \quad
 \begin{cases}
  \cos\alpha=\dfrac{1-\theta^2}{1+\theta^2}\\
  \sin\alpha=\dfrac{2\theta}{1+\theta^2}\\
 \end{cases}
\end{equation*}
with $\theta\in[0,1]$ and $\alpha\in[0,\frac{\pi}{2}]$. 
Since $\theta=\tan (\frac{\alpha}{2})$\,,
$$
\mathrm{d}\theta =\frac{1}{2}\frac{\mathrm{d}\alpha}{(\cos \frac{\alpha}{2})^2}=\frac{\mathrm{d}\alpha }{1+\cos \alpha }=\frac{\mathrm{d}\alpha}{1+\frac{1-\theta^2}{1+\theta^2}}\,,
$$

\noindent
$
(1+\theta ^2)^{-1}\mathrm{d}\theta =\dfrac{\mathrm{d}\alpha}{2}
$
 and we get
\begin{equation}\label{II.16}
\frac{2}{1-\theta ^2}\,\mathrm{d}\theta =\frac{1}{\cos \alpha}\,\mathrm{d}\alpha\,.
\end{equation}

\noindent
Then \eqref{Sto.26} becomes with some abuse of notation
\begin{align}\label{II.17}
 \widetilde \zeta(\theta)=\widetilde \zeta(\alpha)&=\mathrm{i}\frac{\sin\alpha}{\cos\alpha}\, \,\zeta+\frac{1}{\cos\alpha}\,\overline{\zeta}\nonumber\\
 &= \mathrm{i}\tan\alpha\,\zeta+\sqrt{1+(\tan\alpha)^2}\,\, \overline{\zeta}\,.
\end{align}

\noindent
We introduce the change of variable:
\begin{equation}\label{II.18}
 \begin{cases}
  t=\tan\alpha \\
  \mathrm{d}  t=(1+(\tan\alpha)^2)\, \mathrm{d} \alpha\\
  \alpha\in [0,\frac{\pi}{2}[
 \end{cases}
 \iff\quad
  \begin{cases}
   \alpha=\arctan(t)\\
   \mathrm{d} \alpha=\dfrac{1}{1+t^2}\, \mathrm{d}  t\\
   t\in[0,+\infty[
  \end{cases}\,.
\end{equation}
(\ref{II.16}) gives
\begin{equation}\label{II.18,5}
\frac{2}{1-\theta ^2}\mathrm{d}  \theta =\frac{1}{\sqrt{1+t^2}}\mathrm{d}  t\,.
\end{equation}

\par
\eqref{Sto.26}, \eqref{II.17} become with some abuse of notation 
\begin{equation}
\begin{split}\label{II.19}
 \widetilde \zeta(\theta)=\widetilde \zeta(\zeta
  ,\theta)=\widetilde \zeta (t)=\widetilde{\zeta }(t,\zeta )
  &=\mathrm{i} t \,\zeta+\sqrt{1+t^2}\,\, \overline{\zeta}\\
& =\mathrm{i} t\,\zeta+t\,\overline{\zeta}+\big(\sqrt{1+t^2}-t\big)\,\overline{\zeta}\\
& =t(\mathrm{i}\zeta+\overline{\zeta})+\dfrac{1}{\sqrt{1+t^2}+t}\,\overline{\zeta}\,.
\end{split}
\end{equation}
With $0\not=\zeta=|\zeta|\sigma$, $\sigma\in S^1$, we obtain
\begin{equation}\label{II.20}
\begin{split}
\widetilde\zeta(t) & =|\zeta|\left((\mathrm{i}\sigma+\overline{\sigma})\, t+\overline{\sigma}\, \frac{1}{\sqrt{1+t^2}+t}\right)\\
&= |\zeta|\left(2s_+\, \e^{\mathrm{i}\frac{\pi}{4}}t+\overline{\sigma }g_1(t) \right)\,,\\
\end{split}
\end{equation}
 where 
\begin{equation}\label{II.20,1}
 \sigma=s_+\, \e^{-\mathrm{i}\frac{\pi}{4}}+s_-\, \e^{\mathrm{i}\frac{\pi}{4}}\hbox{ with }
s_\pm=\mathrm{Re}  \left(\e^{\pm \mathrm{i}\frac{\pi}{4}}\sigma\right),\hbox{ satisfying}
 \,\,\,\, s_+^2+s_-^2=1\,,
\end{equation}
and 
\begin{equation}\label{II.20,2}
 g_1(t):=\dfrac{1}{\sqrt{1+t^2}+t}\,,\ \ t\geq 0\,,
\end{equation}
satisfies
\begin{equation}\label{II.20,4}
\partial_t^kg_1(t)={\mathcal  O}\left((1+t)^{-1-k}\right),
\end{equation}
for every $k\in \mathbb{N}$.

\begin{remark}
 We can change the variable directly between $\theta\in[0,1[$ and $t\in [0,+\infty[$ by:
 \begin{equation}\label{II.21}
\begin{cases}
  t=\dfrac{2\theta}{1-\theta^2},\\
  \mathrm{d}t=2\dfrac{1+\theta^2}{(1-\theta^2)^2}\,
  \mathrm{d}\theta=\dfrac{1+\theta^2}{1-\theta^2}\dfrac{2}{1-\theta^2}\,
  \mathrm{d}  \theta\,.
  \end{cases}
 \end{equation}
Since 
$$
1+t^2=1+\dfrac{4\theta^2}{(1-\theta^2)^2}=\dfrac{(1-\theta^2)^2+4\theta^2}{(1-\theta^2)^2}=\left(\dfrac{1+\theta^2}{1-\theta^2}\right)^2\,,
$$ 
we get \eqref{II.18,5}.
\end{remark}

\par We obtain the following result:

\begin{prop}\label{Inte.1}
Let $\zeta\in \mathbb{C}\setminus \e^{\mathrm{i}\frac{\pi}{4}}\mathbb{R}  $. Then
$s_+\ne 0$, so $\widetilde{\zeta }(t)\to \infty $ when $t\to +\infty $
and $\mathrm{II}(\zeta ,1)$ is well defined and equal to
\begin{equation}\label{II.22}\begin{split}
    -\int_0^1 a(\zeta,\widetilde
    \zeta(\theta);h)\,&\chi(\zeta,\widetilde\zeta(\theta))\,
    \e^{\frac{\mathrm{i}}{h}\phi(\zeta,\widetilde
      \zeta(\theta))}\,\dfrac{2\mathrm{i}\,\overline{\widetilde{\zeta}(\theta)}}{1-\theta^2}d\theta
    \\
    &= -\int_0^{+\infty} \hskip-5pt a(\zeta,\widetilde
    \zeta(t);h)\,\chi(\zeta,\widetilde\zeta(t))\,
    \e^{\frac{1}{h}\varphi
      (\zeta,t)}\,\dfrac{\mathrm{i}\overline{\widetilde{\zeta}(t)}}{\sqrt{1+t^2}}\,\mathrm{d}  t\,,
\end{split}
\end{equation}
where
\begin{equation}\label{II.23}
\varphi (\zeta ,t)=\mathrm{i}\phi (\zeta ,\widetilde{\zeta }(\zeta ,t))
\end{equation}
and
$\widetilde{\zeta}(\zeta ,
t)$ is given by \eqref{II.20}, \eqref{II.20,1} and \eqref{II.20,2}.
\end{prop}

\par
We now make the more precise assumptions (\ref{6.ft3}),
(\ref{1.ft4}), (\ref{1,5.ft4}), (\ref{2.ft4}), (\ref{3.ft4}),
(\ref{5.ft7}), (\ref{6.ft7}) and use (\ref{8.ft4}), (\ref{1.ft5}).

\par
Recall that $\xi=\frac{\zeta+\widetilde\zeta}{2}$,
$\eta=\frac{\zeta-\widetilde\zeta}{2\mathrm{i}}$ and that the quadratic part of the phase $\phi$ is given by (\ref{2.ft4}), (\ref{1.ft5}):
\begin{equation}\label{II.26}
\begin{split}
 \phi_2(\xi,\eta)&=\frac{\lambda}{2}\xi^2-\frac{\mu}{2}\eta^2+\rho \xi
 \eta \\
 &=\frac{\lambda }{8}(\zeta+\widetilde\zeta)^2+\frac{\mu
 }{8}(\zeta-\widetilde\zeta)^2
 +\mathrm{i}\frac{\rho}{4}(\zeta +\widetilde{\zeta })(\widetilde{\zeta }-\zeta )\,,
\end{split}
\end{equation}
where $\lambda ,\, \mu >0$, $\rho \in \mathbb{R}  $.
Rewrite $\zeta+\widetilde\zeta$ and $\widetilde\zeta-\zeta$, using \eqref{II.20}:
\begin{equation}\label{II.27}
 \zeta+\widetilde\zeta=\vert \zeta\vert\,\Big(2s_{+} \e^{\mathrm{i}\frac{\pi}{4}} t +\overline{\sigma}g_1(t) +\sigma\Big)\,,
 \end{equation}
\begin{equation}\label{II.28}
 \widetilde\zeta-\zeta=\vert \zeta\vert\Big(2s_{+} \e^{\mathrm{i}\frac{\pi}{4}} t
 + \overline{\sigma}g_1(t) -\sigma\Big)\,.
\end{equation}
Then, 
 \begin{align}\label{II.29}
\varphi_2(\zeta ,t):=\mathrm{i}&\phi_2(\zeta ,\widetilde{\zeta }(\zeta ,t))\nonumber\\
 =\mathrm{i}&|\zeta|^2\left(\frac{\lambda}{8} \left(2s_{+} \e^{\mathrm{i}\frac{\pi}{4}} t +\overline{\sigma}g_1(t)
   +\sigma\right)^2 +\frac{\mu }{8}  \left(2s_{+} \e^{\mathrm{i}\frac{\pi}{4}} t +
   \overline{\sigma}g_1(t) -\sigma\right)^2\right.\nonumber\\ &+\left.\mathrm{i}\frac{\rho}{4}\left(2s_{+} \e^{\mathrm{i}\frac{\pi}{4}} t +\overline{\sigma}g_1(t)
   +\sigma\right) \left(2s_{+} \e^{\mathrm{i}\frac{\pi}{4}} t +\overline{\sigma}g_1(t)
   -\sigma\right)\right)\nonumber\\
=\mathrm{i} &|\zeta |^2\left( \left(\frac{\lambda +\mu }{8}+\mathrm{i}\frac{\rho
      }{4}\right)\left(2s_+\e^{\mathrm{i}\frac{\pi}{4}}t+\overline{\sigma }g_1(t)\right)^2+\frac{\lambda -\mu }{4}\sigma   
\left(2s_+\e^{\mathrm{i}\frac{\pi}{4}}t+\overline{\sigma }g_1(t)
  \right)\right.\nonumber\\  &\left. +\left(\frac{\lambda +\mu }{8}-\mathrm{i}\frac{\rho}{4}\right)\sigma ^2 \right).
\end{align}

We complete the square from the first two terms inside the large parenthesis
in the last expression and get after some computations,
\begin{equation}\label{II.30}
  \begin{split}
    \varphi_2(\zeta ,t) =&|\zeta |^2\varphi_2(\sigma ,t)\\
    =&\mathrm{i}\frac{|\zeta |^2 }{8}
    \left(
     (\lambda +\mu +\mathrm{i} 2\rho)\left(2s_+\e^{\mathrm{i}\frac{\pi}{4}}t+\overline{\sigma}g_1(t)+\frac{(\lambda -\mu )\sigma }{\lambda +\mu +\mathrm{i} 2\rho} \right)^2
    \right. \\ & 
    \left. +(\lambda +\mu -\mathrm{i} 2\rho) \left(1-\frac{(\lambda -\mu )^2}{(\lambda +\mu)^2+4\rho^2} \right) \sigma^2\right)\, ,
  \end{split}
\end{equation}
$\zeta =|\zeta |\sigma $.

\begin{lemma}\label{Inte.2}
 For $\zeta=\vert \zeta\vert\,\sigma\in \mathbb{C}\setminus \e^{\mathrm{i}\frac{\pi}{4}}\mathbb{R}  $,
 \begin{equation}\label{II.35}
 \begin{aligned}
   \e^{\frac{1}{h}\varphi_2(\zeta,t)}=&\exp\left({\frac{\mathrm{i}(\lambda+\mu+\mathrm{i} 2\rho
         )}{8}\frac{|\zeta |^2}{h}\left((2s_{+}
       \e^{\mathrm{i}\frac{\pi}{4}})\,t+\frac{\lambda-\mu}{\lambda+\mu+\mathrm{i} 2\rho}\sigma\right)^2}\right)\\
   &\times 
   \exp\left({\frac{\mathrm{i}(\lambda+\mu+\mathrm{i}2\rho)}{8}\frac{|\zeta |^2}{h}\big(\upsilon_0(\sigma,t)+\upsilon_1(t)+\upsilon_2(\sigma,t)\big)}\right),
   \end{aligned}
\end{equation}
where $\upsilon_0, \upsilon_1, \upsilon_2$ are smooth bounded functions, given  by:
\begin{align}
 \upsilon_0(\sigma,t)&=\frac{4\lambda\mu+\rho^2}{(\lambda+\mu)^2+4\rho^2}\sigma^2+4s_{+} \e^{\mathrm{i}\frac{\pi}{4}}\overline{\sigma}g_1(t)t=\mathcal{O}(1)\label{II.36}\\
 \upsilon_1(t)&=2\frac{\lambda-\mu}{\lambda+\mu+\mathrm{i} 2\rho}g_1(t)\underset{+\infty}{\sim} \frac{\lambda-\mu}{\lambda+\mu+\mathrm{i} 2\rho}.\frac{1}{t}\label{II.37}\\
 \upsilon_2(\sigma,t)&=\overline{\sigma}^2g_1(t)^2\underset{+\infty}{\sim}
 \frac{\overline{\sigma}^2}{4t^2}\,,\label{II.38}
\end{align}
satisfying
\begin{equation}\label{II.38,5}
\partial ^k \upsilon_\nu ={\mathcal  O}((1+t)^{-\nu -k}) \hbox{ on }[0,+\infty
[\hbox{ for }\ \nu =0,1,2,\ k\in \mathbb{N}\,. 
\end{equation}
\end{lemma}

The functions $\upsilon_\nu $ belong to the symbol space $\mathbb{S}^{-\nu }=\mathbb{S}^{-\nu }([0,+\infty[)$ in the sense of the following definition,
\begin{dref}\label{1add1}
  For $T_0\in [0,+\infty[$, we put
\begin{multline}\label{4.add2}
 \mathbb{S}^m=\mathbb{S}^m([T_0,+\infty[)=(1+t)^m\mathbb{S}^0([T_0,+\infty [)\\=\{u\in C^\infty([T_0,+\infty[);\, \partial_t^k u=\mathcal{O}((1+t)^{m-k}),\, \forall k\in \mathbb{N}\}\,.
 \end{multline}
 Here $m\in \mathbb{R}$. If $F(t)>0$, $T_0\leq t<\infty,$ we put $F\mathbb{S}^m=\{Fu;\, u\in \mathbb{S}^m\}$. Clearly
\begin{equation}\label{5.add2}
   (1+t)^n \mathbb{S}^m=\mathbb{S}^{m+n}
\end{equation} 
\end{dref}

Summing up Proposition \ref{Inte.1} and Lemma \ref{Inte.2}, we obtain
the following
result when $\phi =\phi _2$:
\begin{prop}\label{Inte.3}
 For $\zeta=\vert \zeta\vert\,\sigma\in \mathbb{C}\setminus \e^{\mathrm{i}\frac{\pi}{4}}\mathbb{R}$, we have
 \begin{multline}\label{II.39}
   -\int_0^1 a(\zeta,\widetilde
   \zeta(\theta);h)\,\chi(\zeta,\widetilde\zeta(\theta))
   \e^{\frac{\mathrm{i}}{h}\phi_2(\zeta,\widetilde
     \zeta(\theta))}\,\dfrac{2\mathrm{i}\,\overline{\widetilde{\zeta}(\theta)}}{1-\theta^2}\,\,\mathrm{d}  \theta\\=
   \frac{1}{\mathrm{i}}\int_0^{+\infty}
   b(\zeta,t)
   \exp\left(\frac{|\zeta|^2}{h}\varphi_2(\sigma ,t)\right)\,\mathrm{d}  t\,,\end{multline}
 where
 \begin{multline}\label{II.39,2}
b(\zeta ,t)=a(\zeta,\vert \zeta\vert\big(2s_{+}\e^{\mathrm{i}\frac{\pi}{4}}
t+\overline{\sigma}g_1(t)\big);h)\\\times \chi\Big(\zeta,\vert\zeta\vert\big(2s_{+}\e^{\mathrm{i}\frac{\pi}{4}}
t+ \overline{\sigma}g_1(t)\big)\Big)
\vert \zeta\vert\Big[2s_{+} \e^{-\mathrm{i}\frac{\pi}{4}}\, t+\sigma
   g_1(t)\Big]\dfrac{1}{\sqrt{1+t^2}}\,,
 \end{multline}

 \begin{align}\label{II.39,4}
   \varphi_2 (\sigma ,t)=&\frac{(\lambda+\mu+\mathrm{i} 2\rho)}{8}\times\\
   &\Big(-\big(2s_+\,t+\e^{-\mathrm{i}\frac{\pi}{4}}
   \frac{\lambda-\mu}{\lambda+\mu+\mathrm{i} 2\rho
      }\sigma\big)^2+
 \mathrm{i}\big(\upsilon_0(\sigma,t)+\upsilon_1(t)+\upsilon_2(\sigma,t)\big)\Big)\,,\nonumber
\end{align}

\noindent
where $\upsilon_j\in \mathbb{S}^{-j}([0,+\infty [\,\ j=0,1,2$ are given by
\eqref{II.36}--\eqref{II.38} and $s_+, g_1$ are given by
\eqref{II.20,1} and \eqref{II.20,2} respectively.
\end{prop}

\subsection{Intermediate values of \texorpdfstring{$\varepsilon $}{}}

Set $\varepsilon :=\frac{|\zeta|}{\sqrt{h}}>0$. The case
$\varepsilon \ge h^{-\delta }$ for some small but fixed $\delta >0$
has been treated in Section \ref{le}, and we shall next study
the case when $\varepsilon \, $ is of the order of 1, or more
precisely when
\begin{equation}\label{II.39,6}
h^\delta \le \varepsilon \,    \le h^{-\delta },
\end{equation}
for some small fixed $\delta >0$. In this case we will not try to get the asymptotics of the integral in the left hand side of (\ref{II.39})
over the whole half line $[0,+\infty [$ but rather for the integral
over a half line $[T,+\infty [$ for a suitable not too large $T>0$ and leave the integral over $[0,T]$ for a numerical study. 

\par
We assume that $\zeta ={\mathcal O}(1)$ and $\zeta \in \mathbb{C}\setminus V$,
where $V$ is a conical open neighborhood in $\mathbb{C}\setminus \{0 \}$ of
$\e^{\mathrm{i}\frac{\pi}{4}}{\mathbb{R}}$. Then

\begin{equation}\label{II.39,8}\mathrm{supp}\big[b(\zeta,\cdot)\big]\subset \Big[0,\frac{C}{|\zeta|}\Big],\quad \hbox{ for some } C>0. \end{equation}

We have 
 $$
 \partial_t^k(a\chi)(\zeta,\vert \zeta\vert\big(2s_{+} \e^{\mathrm{i}\frac{\pi}{4}} t+ \overline{\sigma}g_1(t)\big);h)=\mathcal{O}(|\zeta|^k)=\mathcal{O}(t^{-k}),\quad \forall k\in\mathbb{N}\,.
 $$
 Since we also have $|\zeta|=\mathcal{O}(1)$, we get
 $$
 \partial_t^k(a\chi)(\zeta,\vert \zeta\vert\big(2s_{+} \e^{\mathrm{i}\frac{\pi}{4}} t+ \overline{\sigma}g_1(t)\big);h)=\mathcal{O}(|\zeta|^k)=\mathcal{O}(\frac{1}{(1+t)^{k}}),\,\,\, \forall k\in\mathbb{N}\,,
 $$
 i.e.,
 $$
(a\chi ) (\zeta,\vert \zeta\vert\big(2s_{+} \e^{\mathrm{i}\frac{\pi}{4}} \cdot +
\overline{\sigma}g_1(\cdot )\big);h)\in \mathbb{S}^0 ([0,+\infty [).
 $$
 It follows that 
 \begin{equation}\label{1.add1}
  \partial_t^kb(\zeta,t)=\mathcal{O}(1)\frac{|\zeta|}{(1+t)^{k}},\,\,\, \forall k\in\mathbb{N}\, ,
 \end{equation}
i.e., $b\in |\zeta |\mathbb{S}^0([0,+\infty [)$.

\par We next look at $\varphi_2$: We have $\upsilon_j\in \mathbb{S}^0([0,+\infty[)$, $j=0,1,2$,
 \begin{equation}\label{2.add1}
  \varphi_2(\sigma ,t)\equiv -\frac{\lambda+\mu+\mathrm{i} 2\rho }{8}\Big(2s_{+}
  \,t+\e^{-\mathrm{i}\frac{\pi}{4}}\frac{\lambda-\mu}{\lambda+\mu+\mathrm{i}2\rho}\sigma\Big)^2\ \mathrm{mod\,}\mathbb{S}^0([0,+\infty[),
 \end{equation}
 
\begin{equation}\label{3.add1}
   \partial_t\varphi_2(\sigma ,t)\equiv -\frac{\lambda+\mu+\mathrm{i} 2\rho}{2}s_+\Big( 2s_{+}
   \,t+\e^{-\mathrm{i}\frac{\pi}{4}}\frac{\lambda-\mu}{\lambda+\mu+\mathrm{i} 2\rho}\sigma\Big)
   \mathrm{mod\,}\mathbb{S}^{-1}([0,+\infty[)\, ,
\end{equation}

\begin{equation}\label{4.add1}
  \partial_t\varphi_2(\sigma ,t)\in \mathbb{S}^1([0,+\infty 
  [) \,.
\end{equation}

The assumption that $\zeta $ is not in a conical neighborhood of $\e^{\mathrm{i}\frac{\pi}{4}}]0,+\infty [$ tells us that $|s_+|\geq \frac{1}{\mathcal{O}(1)}$. We have
\begin{equation}\label{5.add1}
-\partial_t\mathrm{Re} \, \varphi_2\asymp\big\vert \partial_t\varphi_2 \big\vert 
\asymp t\qquad\mathrm{when}\quad t\geq T_0\gg1\,,
\end{equation}
and we get
\begin{equation}\label{1.add2}
  \frac{1}{\mathrm{Re} \, (\partial_t\varphi_2)}\, ,\  \frac{1}{\partial_t\varphi_2}\in \mathbb{S}^{-1}([T_0,+\infty [)\,.
\end{equation}

\par 
We can construct an asymptotic anti-derivative of
$b \exp\big(\frac{|\zeta|^2}{h}\varphi_2\big)$ in the limit
$t\rightarrow +\infty$. Look for a suitable $c(t,\zeta;h)$ such that
$$
\partial_t\Big(c e^{\frac{|\zeta|^2}{h}\varphi_2}\Big)\approx b e^{\frac{|\zeta|^2}{h}\varphi_2}\,,
$$
i.e.,
$$
\frac{|\zeta|^2}{h}(\partial_t\varphi_2)c+\partial_t c\approx b\,,
$$
$$
c+\frac{h}{|\zeta|^2\partial_t\varphi_2} \partial_t c\approx \frac{h}{|\zeta|^2\partial_t\varphi_2}b\,.
$$

\par
As an approximate solution we take
\begin{equation}\label{2.add2}
 c^{(N)}=\sum_{k=0}^{N}\Big(-\frac{h}{|\zeta|^2\partial_t\varphi_2}\partial_t\Big)^{k} \frac{h}{|\zeta|^2\partial_t\varphi_2} b =: \sum_{k=0}^{N} c_k\,. 
\end{equation}
 Then
 \begin{equation}\label{3.add2}
  c^{(N)}+\frac{h}{|\zeta|^2\partial_t\varphi_2} \partial
 _tc^{(N)}=\frac{h}{|\zeta |^2\partial_t\varphi_2 }b+\frac{h}{|\zeta
    |^2\partial_t\varphi_2 }\partial_t c_N \,.
\end{equation}
We have on $[T_0,+\infty [$,
 \begin{equation}\label{6.add2}
   \partial_t\varphi_2\in \mathbb{S}^1,\quad (\partial_t\varphi_2)^{-1}\in \mathbb{S}^{-1},\quad b\in |\zeta|\mathbb{S}^0\,.
\end{equation}
By \eqref{6.add2},  we get on $[T_0,+\infty [$, 
\begin{equation}\label{1.add3}
  c_k\in\Big(\frac{h}{(t|\zeta|)^2}\Big)^k\frac{h}{t|\zeta|}\mathbb{S}^0,
\end{equation}
  hence for $T\geq T_0$
\begin{equation}\label{2.add3}
  {{c_k}_\vert}_{[T,+\infty[}\in
\Big(\frac{h}{(T|\zeta|)^2}\Big)^k
  \frac{h}{t|\zeta|}\mathbb{S}^0([T,+\infty [)\,.
\end{equation}
Assume that $\displaystyle \frac{h}{(T|\zeta|)^2}\leq 1$, i.e.,
$T|\zeta|\geq \sqrt{h}$, $T\geq \frac{\sqrt{h}}{|\zeta|}=\frac{1}{\varepsilon }$.
If
\begin{equation}\label{3.add3}
 T\geq \max\big(T_0,\frac{\sqrt{h}}{|\zeta|}\big)\, ,
\end{equation}
we get
\begin{equation}\label{4.add3}
{{c_k}_\vert}_{[T,+\infty[} \in  \frac{h}{t|\zeta|}\mathbb{S}^0([T,+\infty [),
\hbox{ hence }
{{c^{(N)}}_\vert}_{[T,+\infty[}
\in \frac{h}{t|\zeta|}\mathbb{S}^0([T,+\infty [).
\end{equation}
In particular,
\begin{equation}\label{4,5.add3}
c^{(0)}=c_0=\frac{h}{|\zeta |^2\partial_t\varphi_2}b\in
\frac{h}{t|\zeta |}\mathbb{S}^0([T,+\infty [).
\end{equation}

Returning to the original problem of integration, we write \eqref{3.add2} as 
 $$
 \frac{|\zeta|^2\partial_t\varphi_2}{h}\, c^{(N)} + \partial_t c^{(N)}=b+\partial_t c_N\,,
 $$
 or
 \begin{equation}\label{5.add3}
 \partial_t\Big(\e^{\frac{|\zeta|^2}{h}\varphi_2} c^{(N)}\Big)= \e^{\frac{|\zeta|^2}{h}\varphi_2}\big(b+\partial_tc_N\big)\,,
\end{equation}
 where by (\ref{1.add3}),
 \begin{equation}\label{6.add3}
   \partial_t c_N\in \Big(\frac{h}{(t|\zeta|)^2}\Big)^{N+1} \frac{|\zeta|}{t}
   \mathbb{S}^0([T,+\infty [)\,.
   \end{equation}
   Recall that $b\in |\zeta|\mathbb{S}^0$ by \eqref{6.add2}.

 \par
 Integrate \eqref{5.add3} from $T$ to $+\infty$:
 \begin{equation}\label{1.add4}
 \int_T^{+\infty} b(t) \e^{\frac{|\zeta|^2}{h}\varphi_2(\sigma, t)}\, \mathrm{d}t  
 =-\e^{\frac{|\zeta|^2}{h}\varphi_2(\sigma ,T)} c^{(N)}(T)- 
 \int_T^{+\infty} \e^{\frac{|\zeta|^2}{h}\varphi_2(\sigma ,t)}  \big(\partial_t c_N\big)(t)\, \mathrm{d}t  \,.
\end{equation}

 Assume first that $b$ is an elliptic positive element of $|\zeta|\mathbb{S}^0$:
 \begin{equation}\label{2.add4}
 0<b\asymp |\zeta|\,.
\end{equation}
 By \eqref{6.add3} we get
\begin{equation}\label{3.add4}
 \partial_t c_N=\mathcal{O}(1)\Big(\frac{h}{(T|\zeta|)^2}\Big)^{N+1} b
\end{equation}
and \eqref{1.add4}, \eqref{3.add4} give
 \begin{equation}\label{4.add4}
\int_T^{+\infty}  \left(1+\mathcal{O}(1)\left(\frac{h}{(T|\zeta|)^2}\right)^{N+1}\right) b(t) \e^{\frac{|\zeta|^2}{h}\varphi_2(\sigma ,t)}\, \mathrm{d}t  =-\e^{\frac{|\zeta|^2}{h}\varphi_2(\sigma ,T)} c^{(N)}(T).
\end{equation}
With $N=0$ we get in view of (\ref{4.add3}),
\begin{align}\label{4,5.add4}
\int_T^\infty \left(1+{\mathcal  O}\left(\frac{h}{(T|\zeta|)^2} \right)
\right) b(t) \e^{\frac{|\zeta|^2}{h}\varphi_2 (\sigma,t)}\,\mathrm{d}t   &=-\e^{\frac{|\zeta|^2}{h}\varphi_2 (\sigma,T)}c_0(T)\nonumber\\
&={\mathcal  O}(1)\e^{\frac{|\zeta|^2 }{h}\varphi_2 (\sigma ,T)}\frac{h}{T|\zeta |^2}|\zeta|\,.
\end{align}

In this discussion we can replace $\varphi_2$ by $\mathrm{Re}(\varphi_2)$, in particular in (\ref{4,5.add4}).  Assume that $T\ge T_0$,
$T\gg \frac{\sqrt{h}}{|\zeta|}$ (cf. \eqref{3.add3}). With $N=0$, we first get
$$
\int_T^{+\infty }b(t)\,\e^{\frac{|\zeta |^2}{h}(\mathrm{Re} \,\varphi_2)
(\sigma ,t)}\,\mathrm{d}t  ={\mathcal  O}(1)\e^{\frac{|\zeta|^2}{h}(\mathrm{Re} \varphi_2)(\sigma ,t)}\frac{h}{T|\zeta|}\,,
$$
which implies the estimate
\begin{equation}\label{5.add4}
\int_T^{+\infty }b(t)\,\e^{\frac{|\zeta |^2}{h} \varphi_2
(\sigma ,t)}\,\mathrm{d}t  ={\mathcal  O}(1)\e^{\frac{|\zeta |^2}{h} \varphi_2(\sigma ,t)}\frac{h}{T|\zeta|}\,,
\end{equation}
where we can weaken the ellipticity, positivity and symbol class
assumptions on $b$ and merely assume that $b={\mathcal 
  O}(|\zeta |)$ on $[T,+\infty [$.

\par
Using this and (\ref{6.add3}) to estimate the last integral in
(\ref{1.add4}), we get
\begin{prop}\label{1add4}
For $T\ge T_0$, $T\gg \sqrt{h}/|\zeta|$, we have for every $b\in
|\zeta |\mathbb{S}^0([0,+\infty [)$ and $N\in \mathbb{N}$:
 \begin{equation}\label{6.add4}
   \int_T^{+\infty} b(t) \e^{\frac{|\zeta|^2}{h}\varphi_2(\sigma ,t)}\,\mathrm{d}t  =-\e^{\frac{|\zeta|^2}{h}\varphi_2(\sigma ,T)} \left(
     c^{(N)}(T)+\mathcal{O}(1)\left(\frac{h}{(T|\zeta|)^2}\right)^{N+1}\frac{h}{T|\zeta
       |}\right).
 \end{equation}
 Here $c^{(N)}$ is defined in (\ref{2.add2}) and we
 recall (\ref{4.add3}).
\end{prop}

  We shall extend this result to the case when the quadratic phase
  $\phi_2 $ is replaced by  the full  phase $\phi $, where 
  \begin{equation}\label{1.nonq1}
\phi (\zeta ,\widetilde{\zeta })=\phi_2(\zeta ,\widetilde{\zeta })+r(\zeta ,\widetilde{\zeta })\,,
  \end{equation}
and
  \begin{equation}\label{2.nonq1}
r(\zeta ,\widetilde{\zeta })={\mathcal  O}((\zeta ,\widetilde{\zeta })^3)
\end{equation}
is holomorphic in $\mathrm{neigh}((0,0);\mathbb{C}^2)$ and the cutoff has
its support in a sufficiently small neighborhood of $(0,0)\in {\mathbb{
  C}}^2_{\zeta ,\tilde{\zeta }}$.

\par 
We study the symbol properties of $r(\zeta ,\widetilde{\zeta
}(\zeta ,t))$, where $\widetilde{\zeta }(\zeta ,t)$ is given by
(\ref{II.20}) and $\zeta =|\zeta |\sigma $,
\begin{equation}\label{3.nonq1}
\widetilde{\zeta }(\zeta ,t)=\left((i\zeta +\overline{\zeta
  })t+\overline{\zeta }g_1(t) \right),
\end{equation}
$g_1(t)=\left(t+\sqrt{1+t^2} \right)^{-1}\in \mathbb{S}^{-1}([0,+\infty [)$ and
we keep the assumption that $|s_+(\zeta )|\ge \frac{1}{{\mathcal  O}(1)}$, so that
\begin{equation}\label{4.nonq1}
|i\zeta +\overline{\zeta }|\asymp |\zeta |\,.
\end{equation}

\par 
For any fixed $T>0$, we have
\begin{equation}\label{5.nonq1}
r(\zeta ,\widetilde{\zeta }(\zeta ,t))={\mathcal  O}((\zeta
,\overline{\zeta })^3),\hbox{ uniformly for }t\in [0,T]\,,
\end{equation}
and this function is real-analytic for $\zeta \in \mathrm{neigh}(0;\mathbb{C})$.

\par
By (\ref{3.nonq1}),
\begin{equation}\label{6.nonq1}
  \partial_t^j\partial_{\zeta }^k\partial_{\overline{\zeta }}^\ell
  (\widetilde{\zeta }(\zeta ,t))
  =\begin{cases}
    {\mathcal  O}(1)|\zeta |^{1-k-\ell}(1+t)^{1-j},\ \ 0\le k+\ell \le 1\\
    0,\ \ k+\ell \ge 2
  \end{cases}
\end{equation}
and this implies that
\begin{equation}\label{7.nonq1}
\partial_t^j\partial_{\zeta }^k\partial_{\overline{\zeta }}^\ell
(\widetilde{\zeta }(\zeta ,t))
={\mathcal  O}(1)|\zeta |^{(1-k-\ell )_+}(1+t)^{1-j}\,.
\end{equation}

\par For $T>0$ sufficiently large, we have by (\ref{3.nonq1}), (\ref{4.nonq1})
\begin{equation}\label{1.nonq2}
|\widetilde{\zeta }(\zeta ,t)|\asymp |\zeta |t,\ \ t\in [T,+\infty [\,.
\end{equation}

\par By (\ref{2.nonq1}) and Taylor expansion, we can write
\begin{equation}\label{4.nonq2}
r=\sum_{p+q=3}r_{p,q},
\end{equation}
where $r_{p,q}$ are holomorphic in a neighborhood of $(0,0)$ and
satisfy
\begin{equation}\label{5.nonq2}
r_{p,q}(\zeta ,\widetilde{\zeta })={\mathcal  O}(1)\zeta ^p\widetilde{\zeta
}^q,\ \ p,q\ge 0.
\end{equation}
(If necessary we shrink the support of $\chi $, to be contained in a
neighborhood of $(0,0)$, where (\ref{5.nonq2}) holds.)

\par By the chain rule,
\begin{equation}\label{6.nonq2}
  \begin{split}
\partial _t^j\partial_{\zeta }^m\partial_{\overline{\zeta }}^\ell
&\left(r_{p,q}(\zeta ,\widetilde{\zeta }(\zeta ,t))\right)= \hbox{ a
  finite linear combination of terms}\\
&\left(\partial_{\zeta }^M\partial_{\widetilde{\zeta }}^L\, r_{p,q}
\right)(\zeta ,\widetilde{\zeta }(\zeta ,t))\prod_{\nu =1}^L\left(
\partial_{\zeta }^{m_\nu }\partial_{\overline{\zeta }}^{\ell_\nu
}\partial_t^{j_\nu }\widetilde{\zeta }(\zeta ,t)
\right)
  \end{split}
\end{equation}
with
\begin{equation}\label{7.nonq2}
m_\nu +\ell_\nu +j_\nu\ne 0,\ \ m_\nu +\ell_\nu \le 1, 
\end{equation}
\begin{equation}\label{8.nonq2}
M+\sum_\nu m_\nu =m,\ \ \sum_\nu \ell_\nu =\ell,\ \ \sum_{\nu }j_\nu =j\,.
\end{equation}
When $\ell =j=0$, we also have the simple term $\partial_{\zeta
}^m(r_{p,q})(\zeta ,\widetilde{\zeta }(\zeta ,t))$ with $M=m$, $L=0$.

By (\ref{7.nonq1}), (\ref{1.nonq2}) the general term in
(\ref{6.nonq2}) is
\begin{equation}\label{9.nonq2}
  ={\mathcal  O}(1)\zeta ^{(p-M)_+}(\zeta t)^{(q-L)_+}\prod_{\nu
  =1}^L\left(\zeta ^{(1-m_\nu -\ell_\nu )_+}t^{1-j_\nu } \right)\,.
\end{equation}
Since $|\zeta |\le {\mathcal  O}(1)$, $|\zeta |t\le {\mathcal  O}(1)$,
we can drop the positive part subscripts
and see
that the general term in (\ref{6.nonq2}) is
$$
={\mathcal  O}(1)\zeta^{\big(p-M+q-L+\overset{L}{\underset{\nu=1}{\sum}}(1-m_\nu -\ell_\nu
  )\big)}\,\, t^{\big(q-\overset{L}{\underset{\nu=1}{\sum}}j_\nu\big)}
={\mathcal  O}(1)\zeta ^{p+q-m-\ell} t^{q-j}\,.
$$

\par 
By (\ref{4.nonq2}), (\ref{6.nonq2}), we get
\begin{equation}\label{1.nonq3}
\partial_{\zeta }^m\partial_{\overline{\zeta }}^\ell \partial
_t^j(r(\zeta ,\widetilde{\zeta }(\zeta ,t)))={\mathcal  O}(1)\zeta
^{3-m-\ell}t^{3-j}, \quad t\ge T\,.
\end{equation}
In particular,
\begin{equation}\label{2.nonq3}
\partial_t(r(\zeta ,\widetilde{\zeta }(\zeta ,t))={\mathcal  O}(1)\zeta
^3t^2, \ t\ge T
\end{equation}
and more generally,
\begin{equation}\label{3.nonq3}
r(\zeta ,\widetilde{\zeta }(\zeta ,\cdot ))\in \zeta ^3 \mathbb{S}^3([0,+\infty [),\
\hbox{ uniformly in }\zeta\,.
\end{equation}

Recall that $\varphi_2(\zeta ,t)=\mathrm{i}\phi_2(\zeta ,\widetilde{\zeta
}(\zeta ,t))$ satisfies (\ref{5.add1}),
\begin{equation}\label{4.nonq3}
|\partial_t\varphi_2 (\zeta ,t)|\asymp |\zeta |^2t,\ \ t\ge T_0\,.
\end{equation}

When
$\mathrm{supp}(\chi )$ is contained in a sufficiently small
neighborhood of $(0,0)$,
we can make $|\zeta |t$ as small as we like by
(\ref{1.nonq2}). Then we get $|\zeta |^3t^2\ll |\zeta |^2t$, so
\begin{equation}\label{4,5.nonq3}
|\partial_t (r(\zeta ,\widetilde{\zeta }(\zeta ,t)))|
\ll |\partial_t\varphi_2(\zeta ,t)|\,, 
\end{equation}  
implying
\begin{equation}\label{5.nonq3}
|\partial_t\varphi(\zeta ,t)|\asymp |\zeta |^2t,\hbox{ when }t\in
[T_0,+\infty [ \hbox{ and } (\zeta ,\widetilde{\zeta }(\zeta ,t))\in
\mathrm{supp}(\chi )\,.
\end{equation}

\par By the same argument,
\begin{equation}\label{6.nonq3}
\partial_t\varphi \in |\zeta|^2 \mathbb{S}^1([T_0,+\infty [),\hbox{ uniformly in
}\zeta 
\end{equation}
and this and (\ref{5.nonq3}) imply that (\ref{5.add1}), (\ref{1.add2})
hold with $\varphi_2$ replaced by $\varphi$. The proof of Proposition
\ref{1add4} then runs with $\varphi_2$ replaced by $\varphi$. Hence,
\begin{prop}\label{1nonq3} Assume that $\chi $ in (\ref{7.ft7})  has
  sufficiently small support.
For $T\ge T_0$, $T\gg \frac{\sqrt{h}}{|\zeta|}$, we have for every $b\in
|\zeta |\mathbb{S}^0([0,+\infty [)$ and $N\in \mathbb{N}$:
 \begin{equation}\label{7.nonq3}
   \int_T^{+\infty} b(t) \e^{\frac{|\zeta|^2}{h}\varphi (\sigma ,t)}\,
   \mathrm{d}t  =-\e^{\frac{|\zeta|^2}{h}\varphi (\sigma ,T)} \left(
     c^{(N)}(T)+\mathcal{O}(1)\left(\frac{h}{(T|\zeta|)^2}\right)^{N+1}\frac{h}{T|\zeta
       |}\right).
 \end{equation}
 Here $c^{(N)}$ is defined as in (\ref{2.add2}), but with $\varphi_2$
 replaced by $\varphi$,  (\ref{4.add3}) remains valid.
\end{prop}

\par\noindent 
The integral $\int _0^T b(t) \e^{\frac{1}{h}\varphi(\zeta ,t)}\,
\mathrm{d}t$ can be computed to wanted precision with standard numerical methods, no asymptotic expressions are needed for this expression.

\subsection{The small \texorpdfstring{$\varepsilon$ limit}{}}

We now turn to the case when
$\varepsilon \le h^\delta $ for some small fixed $\delta >0$.
Recall some quantities. We consider the integral in equation (\ref{II.22}), which, recalling (\ref{II.39}), can be written as
\begin{equation}\label{1.sep1}
  \frac{1}{\mathrm{i}}\int_0^{+\infty} b(\zeta ,t;h)\, \e^{\frac{1}{h}\varphi(\zeta ,t)}\, \mathrm{d} t\,,
\end{equation}

\begin{equation}\label{2.sep1}
b(\zeta ,t;h)=(a\chi )(\zeta ,\widetilde{\zeta }(\zeta
,t))\frac{\overline{\widetilde{\zeta }}(\zeta ,t)}{\sqrt{1+t^2}}\,,
\end{equation}
\begin{equation}\label{3.sep1}
\varphi(\zeta ,t)=\mathrm{i}\phi (\zeta ,\widetilde{\zeta}(\zeta ,t))\,.
\end{equation}
Here $\widetilde{\zeta }(\zeta ,t)$ is given by (\ref{II.20}),
(\ref{3.nonq1}),
\begin{equation}\label{4.sep1}
\widetilde{\zeta }(\zeta ,t)=|\zeta |(2s_+(\sigma )\e^{\mathrm{i}\frac{\pi}{4}}t+\overline{\sigma }g_1(t))=(\mathrm{i}\zeta +\overline{\zeta})t+\overline{\zeta }g_1(t)
\end{equation}
and $\sigma =\dfrac{\zeta}{|\zeta|}$, with the assumption that
$|s_+(\sigma )|\ge \frac{1}{{\mathcal  O}(1)}$ or equivalently that $\sigma $ is
not in a conic neighborhood of $\e^{\mathrm{i}\frac{\pi}{4}}\mathbb{R}$ in ${\mathbb{C}}\setminus \{0 \}$. $g_1$ is given by (\ref{II.20,2}),
\begin{equation}\label{5.sep1}
g_1(t)=\frac{1}{t+\sqrt{t^2+1}}=\sqrt{t^2+1}-t \,.
\end{equation}

\par 
Recall that $\varepsilon \,    =\dfrac{|\zeta |}{\sqrt{h}}$. The first identity in $(\ref{4.sep1})$ can be written
\begin{align}\label{6.sep1}
\widetilde{\zeta }(\zeta,t)&=\sqrt{h}\varepsilon \,    \left(2s_+(\sigma) \e^{\mathrm{i}\frac{\pi}{4}}t+\overline{\sigma }g_1(t)\right)\\
  &=\sqrt{h}\varepsilon \,    \widetilde\zeta(\sigma,t)\,.\nonumber
  \end{align}
  
Let $\phi_2(\zeta ,\widetilde{\zeta})$ be the leading quadratic form part in the Taylor expansion of $\phi $ at $(\zeta ,\widetilde{\zeta})=(0,0)$,
\begin{equation}\label{7.sep1}
\phi (\zeta,\widetilde{\zeta })=\phi_2(\zeta ,\widetilde{\zeta
})+r(\zeta ,\widetilde{\zeta })\,,
\end{equation}
\begin{equation}\label{8.sep1}
  r(\zeta ,\widetilde{\zeta })={\mathcal  O}((\zeta ,\widetilde{\zeta}
  )^3)\hbox{ in }\mathrm{neigh}((0,0);\mathbb{C}^2)\,.
  \end{equation}
  
Recall that 
\begin{equation}\label{9.sep1}
\varphi_2(\zeta ,t)=\mathrm{i}\phi_2(\zeta ,\widetilde{\zeta }(\zeta ,t))\,,
\end{equation}
where $\phi_2$ is given by (\ref{II.26}).

The function $\phi_2$ is positively homogeneous of degree 2 and $\widetilde{\zeta
}(\zeta ,t)$ is positively homogeneous of degree 1 with respect to
$\zeta $, hence $\varphi (\zeta ,t)$ is positively homogeneous of
degree 2 with respect to $\zeta $,
\begin{equation}\label{1.sep2}
  \varphi_2(\zeta ,t)=|\zeta |^2\varphi_2(\sigma ,t)\,,
\end{equation}
or equivalently,
\begin{equation}\label{2.sep2}
  \frac{1}{h}\varphi_2(\zeta ,t)=\varepsilon^2\varphi_2(\sigma ,t)\,.
\end{equation}
Recall (\ref{II.39,4}) in an abbreviated form,
 \begin{equation}\label{4.sep2}
  \varphi_2 (\sigma  ,t)=\\
   \frac{(\lambda+\mu+\mathrm{i} 2\rho)}{8}\Big(-\big(2s_+\,t+\e^{-\mathrm{i}\frac{\pi}{4}}
   \frac{\lambda-\mu}{\lambda+\mu+\mathrm{i} 2\rho
      }\sigma\big)^2+\mathrm{i} \upsilon(\sigma,t)\Big)\,,
  \end{equation}
where $\displaystyle \upsilon=\sum_{j=0}^{2}\upsilon_j\in \mathbb{S}^0([0,+\infty [)$ is described in
\eqref{II.36}--\eqref{II.38} and $s_+$ is given by
\eqref{II.20,1}. Write
\begin{equation}\label{4,5.sep2}\varphi=\mathrm{i}\phi (\zeta ,\widetilde{\zeta
  }(\zeta ,t))\,,
  \end{equation}
\begin{equation}\label{5.sep2}
\varphi =\varphi_2+\psi ,\ \psi (\zeta ,t)=\mathrm{i} r(\zeta ,\widetilde{\zeta
}(\zeta ,t))\,,
\end{equation}
\begin{equation}\label{6.sep2}
  r(\zeta ,\widetilde{\zeta })=\sum_{k=3}^\infty r_k(\zeta ,
  \widetilde{\zeta })=\sum_{k=3}^{K-1} r_k(\zeta ,
  \widetilde{\zeta })+r^K(\zeta ,\widetilde{\zeta })\,,
\end{equation}
in $\mathrm{neigh}((0,0);\mathbb{C}^2)$, where $r_k={\mathcal  O}((\zeta
,\widetilde{\zeta })^k)$ is homogeneous of degree $k$, $r^K={\cal
  O}((\zeta ,\widetilde{\zeta })^K)$.
Then since $\zeta =h^{\frac12}\varepsilon \sigma $, $\widetilde{\zeta }(\zeta
,t)=h^{\frac12}\varepsilon \widetilde{\zeta }(\sigma ,t)$\,,
\begin{equation}\label{7.sep2}
  \frac{1}{h}\psi = \mathrm{i}\sum_3^{K-1} h^{\frac{k}{2}-1}\varepsilon ^k
  \underbrace{r_k(
    \sigma ,\widetilde{\zeta }(\sigma
    ,t))}_{=:r_k(\sigma ,t)}+\mathrm{i}h^{\frac{K}{2}-1}\varepsilon ^K\widetilde{r}_K(\sigma
  ,t;h^{\frac{1}{2}}\varepsilon )\,,
\end{equation}
where \begin{equation}\label{8.sep2}
\widetilde{r}_K(\sigma ,t;h^{\frac{1}{2}}\varepsilon ):=(\varepsilon
h^{\frac{1}{2}})^{-K}r^K(\varepsilon h^{\frac{1}{2}}\sigma , \varepsilon
h^{\frac{1}{2}}\widetilde{\zeta }(\sigma ,t))\,.
\end{equation}
Here $r_k(\sigma ,t)$, $\widetilde{r}_K(\sigma,
t;h^{\frac{1}{2}}\varepsilon )$ are well defined for
\begin{equation}\label{2.sep3}
h^{\frac{1}{2}}\varepsilon \, t={\mathcal  O}(1)\,,
\end{equation}
and belong respectively to
$(1+t)^k\mathbb{S}^0([0,{\mathcal  O}(1)(h^{\frac{1}{2}}\varepsilon 
)^{-1}[)$ and $(1+t)^K\mathbb{S}^0([0,{\mathcal  O}(1)(h^{\frac{1}{2}}\varepsilon 
)^{-1}[)$. In (\ref{7.sep2}) we have 
\begin{equation}\label{3.sep3}\begin{split}
h^{\frac{k}{2}-1}\varepsilon \,    ^kr_k(\sigma ,t))&\in h^{-1}h^{\frac{k}{2}}(\varepsilon 
(1+t))^k\mathbb{S}^0\,,\\
h^{\frac{K}{2}-1}\varepsilon \,    ^K\widetilde{r}_K(\sigma ,t;\varepsilon
h^{\frac{1}{2}}))&\in h^{-1}h^{\frac{K}{2}}(\varepsilon \,    (1+t))^K\mathbb{S}^0\,.
\end{split}
\end{equation}
In particular for $K=3$
\begin{equation}\label{8.sep3}
\frac{1}{h}\psi \in h^{\frac{1}{2}}(\varepsilon (1+t))^3\mathbb{S}^0.
\end{equation}

\par We have $h^{-1}\varphi _2 \in (\varepsilon (1+t))^2\mathbb{S}^0$
and
\begin{equation}\label{4.sep3}
  \frac{1}{h}\mathrm{Re} \, \varphi_2 (\zeta ,t)=
\varepsilon^2\mathrm{Re} \, \varphi_2(\sigma ,t)\asymp -(\varepsilon \, t)^2,\hbox{ for }t\gg 1.
\end{equation}
(\ref{8.sep3}), (\ref{4.sep3}) imply that
\begin{equation}\label{5.sep3}
\frac{1}{h}|\psi |\ll -\frac{1}{h}\mathrm{Re} \, \varphi _2(\zeta ,t)
\end{equation}
when $t\gg 1$ and $h^{\frac12}\varepsilon t\ll 1$. The last bound is
equivalent to $|\zeta |t\ll 1$, or:
$$
|\zeta |\ll 1\hbox{ and } |\widetilde{\zeta }(\zeta ,t)|\ll 1
$$
which is satisfied for our integrals if $\mathrm{supp}\chi $ is
contained in a sufficiently small neighborhood of $(0,0)\in
\mathbb{C}^2$, which we assume from now on. In this region we have
$$
\e^{\frac{1}{h}\varphi (\zeta ,t)}={\mathcal O}(1)\e^{-(h^{-2\delta
  }+(\varepsilon t)^2)/{\mathcal
    O}(1)},\hbox{ when }\varepsilon t\ge h^{-\delta},
$$
for any fixed $\delta $ with $0<\delta \ll 1$. It follows that the
integral (\ref{1.sep1}) can be decomposed as
\begin{equation}\label{6.sep3}
\int_0^{+\infty }b(\zeta ,t)\e^{\frac{1}{h}\varphi(\zeta ,t)} \, \mathrm{d}t=
\int_{0\le \varepsilon t\le h^{-\delta }} \hskip-20pt b(\zeta
,t)\e^{\frac{1}{h}\varphi(\zeta ,t)} \, \mathrm{d}t
+{\mathcal O}\left(\frac{1}{\varepsilon} \right)\e^{-h^{-2\delta }/{\mathcal O}(1)}\,,
\end{equation}
when $b={\mathcal O}(1)$.

\par\smallskip\noindent 
On the interval $[0,(\varepsilon h^\delta )^{-1}]$ we have $\varepsilon
t\le h^{-\delta }$ and we choose $\delta \in
]0,1/6[$ and restrict the attention to this interval until further
notice.  (\ref{3.sep3}), (\ref{8.sep3}) give for $k,\, K\ge 3$,
\[
  \begin{split}
h^{\frac{k}{2}-1}\varepsilon^k r_k(\sigma ,t) &\in \, h^{k(\frac{1}{2}-\delta
  )-1} \, \mathbb{S}^0,\\
h^{\frac{K}{2}-1}\varepsilon ^K\widetilde{r}_K(\sigma ,t;\varepsilon h^{\frac12})
&\in \, h^{K(\frac{1}{2}-\delta
  )-1} \, \mathbb{S}^0,\\
\frac{1}{h}\psi &\in \, h^{\frac{1}{2}-3\delta } \, \mathbb{S}^0\,,
\end{split}
\]
where the exponents are positive,
$$
k(\frac{1}{2}-\delta )-1,\quad K(\frac{1}{2}-\delta )-1,\quad
\frac{1}{2}-3\delta >0\,. 
$$
We can therefore Taylor expand the exponential function at $0$ and get
on $[0,(\varepsilon h^\delta )^{-1}]$ (cf.\ (\ref{8.sep3})),

$$
\e^{\frac{\psi}{h}}=\sum_{n=0}^{N-1}\frac{1}{N!}\left(\frac{\psi}{h}\right)^n+{\mathcal
  O}(1)h^{\frac{N}{2}}(\varepsilon (1+t))^{3N}\,.
$$
Here we apply (\ref{7.sep2}) to each $\psi /h$ and see that on
$[0,(\varepsilon h^\delta )^{-1}]$
\begin{equation}\label{2.sep4}\begin{split}
\e^{\frac{\psi}{h}}=&{\mathcal 
  O}(1) h^{\frac{N}{2}}(\varepsilon (1+t))^{\widetilde{N}} +\hbox{a linear combination of products }\\
 & \prod_{j=1}^n \left(h^{\frac{k_j}{2}-1}\varepsilon \,    ^{k_j}r_{k_j}(\sigma
  ,t) \right), 
\quad\hbox{with }\sum (\frac{k_j}{2}-1)<\frac{N}{2},
\end{split}
\end{equation}
including the term 1 for $n=0$.
Here $\widetilde{N}=\widetilde{N}(N)$ is large enough.
The coefficients in the linear combination are constant and explicitly
computable.
The products in (\ref{2.sep4}) are
$$
=
{\mathcal  O}(1)h^{\overset{n}{\underset{j=1}{\sum}}(\frac{k_j}{2}-1)}
(\varepsilon (1+t))^{\underset{j=1}{\overset{n}{\sum}} k_j} \,.
$$

\par
We now return to the integral (\ref{1.sep1}). First we
shall only use that
\begin{equation}\label{rc.1}
b={\mathcal  O}(1)\,.
\end{equation}
From (\ref{4.sep3}), (\ref{6.sep3}) we get
\begin{equation}\label{rc.2}
  \begin{split}
\int_0^{+\infty}  b(\zeta
  ,t;h)\,\e^{\frac{1}{h}\varphi_2(\zeta ,t)}\, \mathrm{d} t&={\mathcal  O}\left(\frac{1}{\varepsilon \,    }
  \right),\\
  \int_{(\varepsilon h^\delta )^{-1}}^{+\infty}  b(\zeta
  ,t;h)\,\e^{\frac{1}{h}\varphi_2(\zeta ,t)}\, \mathrm{d} t&=
  {\mathcal O}\left(\frac{1}{\varepsilon } \right)\e^{-h^{-2\delta }/{\mathcal O}(1)}\,,
\end{split}
  \end{equation}
and similarly with $\phi $ instead of $\phi _2$.
  More generally, using that $\e^{-\varepsilon ^2 \frac{t^2}{C}} (\varepsilon \, t)^{\widetilde N}=\mathcal{O}(1)$, for all $C>0$, we have 
  \begin{equation}\label{rc.2,5}\begin{split}
  \int_0^{+\infty} b(\zeta,t;h)\, \e^{\frac{1}{h}\varphi_2(\zeta
    ,t)}(\varepsilon \,    t)^{\widetilde N}
  \, \mathrm{d}t&=\mathcal{O}\left(\frac{1}{\varepsilon }\right),\\
  \int_{(\varepsilon h^\delta )^{-1}}^{+\infty} b(\zeta,t;h)\, \e^{\frac{1}{h}\varphi_2(\zeta
    ,t)}(\varepsilon \,    t)^{\widetilde N}
  \, \mathrm{d}t
  &=
  {\mathcal O}\left(\frac{1}{\varepsilon } \right)\e^{-h^{-2\delta }/{\mathcal O}(1)}\,,
\end{split}
  \end{equation}
and similarly with $\phi $ instead of $\phi_2$.

\par 
Multiply (\ref{2.sep4}) with $b\e^{\varphi _2/h}$, integrate from
$0$ to $(\varepsilon h^\delta )^{-1}$ and apply the above estimates:
\[\begin{split}
    &\int_0^{(\varepsilon h^\delta )^{-1}} b(\zeta
    ,t;h)\,\e^{\frac{1}{h}\varphi (\zeta ,t)}\,\mathrm{d} t ={\mathcal
      O}\left(\frac1\varepsilon \,    \right)h^{\frac{N}{2}}\ + \hbox{ a linear
      combination of terms}\\
    &\int_0^{(\varepsilon h^\delta )^{-1}} b(\zeta ,t;h)\,
    \e^{\frac{1}{h}\varphi_2 (\zeta ,t)}\prod_{j=1}^n
    \left(h^{\frac{k_j}{2}-1}\varepsilon 
      ^{k_j}r_{k_j}(\sigma,\widetilde{\zeta }(\sigma ,t)) \right)
    \mathrm{d} t
\end{split}
\]
with $\sum (\frac{k_j}{2}-1)<\frac{N}{2}.$
Here the general integral on the right hand side is equal to
\begin{equation}\label{rc.3}\begin{split}
 & \int_0^{\infty }\mathrm{idem\,} \, \mathrm{d}t-
  \int_{(\varepsilon h^\delta )^{-1}}^{\infty }\mathrm{idem\,} \, \mathrm{d}t\\
  &=\int_0^{+\infty}\mathrm{idem\,} \, \mathrm{d}t+{\mathcal
    O}\left(\frac{1}{\varepsilon }\right)h^{\underset{j=1}{\overset{n}{\sum}}(\frac{k_j}{2}-1)}\e^{-h^{-2\delta }/{\mathcal O}(1)}\,,
\end{split}
\end{equation}
where "idem" denotes the integrand in the general integral in the
right hand side of (\ref{rc.3}).

Recall that $\varphi=\varphi_2+\psi$ and apply (\ref{2.sep4}) to get
\begin{equation}\label{rc.3,5}\begin{split}
 &\int_0^{+\infty} b(\zeta ,t;h)\,\e^{\frac{1}{h}\varphi (\zeta ,t)}\,\mathrm{d} t ={\mathcal 
  O}(\frac1\varepsilon \,    )h^{\frac{N}{2}}\ + \hbox{ a linear
  combination of terms}\\
&\int_0^{+\infty} b(\zeta ,t;h)\, \e^{\frac{1}{h}\varphi_2 (\zeta ,t)}\prod_{j=1}^n \left(h^{\frac{k_j}{2}-1}\varepsilon \,    ^{k_j}r_{k_j}(\sigma,\widetilde{\zeta }(\sigma ,t)) \right)
\mathrm{d} t 
\end{split}
\end{equation}
with $\sum (\frac{k_j}{2}-1)<\frac{N}{2}.$
As in \eqref{rc.2,5}, each integral in the right hand side of (\ref{rc.3,5}) is
 ${\mathcal  O}(\frac1\varepsilon)h^{\underset{j=1}{\overset{n}{\sum}} (\frac{k_j}{2}-1)}$.

\par
Let $\varphi_2^0$ denote the left hand side of (\ref{4.sep2})
with $\upsilon$ replaced by 0, so that
\begin{equation}\label{rc.4}
  \varphi_2(\sigma ,t)=\varphi_2^0(\sigma ,t)+\widetilde{\upsilon}(\sigma
  ,t),\ \
\widetilde{\upsilon} (\sigma ,t):=\frac{\lambda +\mu +\mathrm{i} 2\rho }{8}\,\mathrm{i} \upsilon(\sigma ,t)
\end{equation}
and
\begin{equation}\label{rc.5}
\frac{1}{h}\varphi_2(\zeta ,t)=\varepsilon^2\varphi_2^0(\sigma
,t)+\varepsilon^2\widetilde{\upsilon}(\sigma ,t)\,.
\end{equation}
Write
$$
\e^{\varepsilon^2\widetilde{\upsilon}}=\sum_{\nu=0}^{L-1}\frac{1}{\nu !}(\varepsilon^2\widetilde{\upsilon})^\nu +{\mathcal  O}(\varepsilon^{2L})\,,
$$
$$
\e^{\frac{1}{h}\varphi_2(\zeta ,t)}=\sum_{\nu=0}^{L-1}\frac{1}{\nu !}(\varepsilon^2\widetilde{\upsilon})^\nu \e^{\varepsilon^2\varphi_2^0(\sigma ,t)} +{\mathcal 
  O}(\varepsilon \,    ^{2L})\, \e^{\varepsilon^2\varphi_2^0(\sigma ,t)}\,.
$$
Using this in (\ref{rc.3,5}), we get
\begin{equation}\label{rc.6}
  \begin{split}
 &\int_0^{+\infty} b(\zeta ,t;h)\, \e^{\frac{1}{h}\varphi (\zeta ,t)}\mathrm{d} t ={\mathcal 
  O}(\frac{1}{\varepsilon })\left(h^{\frac{N}{2}}+\varepsilon 
  ^{2L}\right) + \hbox{ a linear
  combination}\\
&\hbox{of terms }\int_0^{+\infty} b(\zeta ,t;h)\,\e^{\varepsilon ^2\varphi_2^0 (\sigma,t)}
(\widetilde{\upsilon}(\sigma ,t)\varepsilon^2)^\nu 
\prod_{j=1}^n \left(h^{\frac{k_j}{2}-1}\varepsilon^{k_j}r_{k_j}(\sigma
  ,\widetilde{\zeta }(\sigma ,t)) \right)
\mathrm{d} t\\
&\hbox{with }\sum (\frac{k_j}{2}-1)<\frac{N}{2},\ \nu <L.
\end{split}
\end{equation}

\noindent 
{\bf Study of the general term in (\ref{rc.6}).}

For $k\in [3,+\infty [\cap \mathbb{N}$, we have 
\begin{equation}\label{1.gt1}
  r_k(\sigma,\widetilde\zeta (\sigma,t))=r_k\big(\sigma ,
  2s_+(\sigma)\e^{\mathrm{i}\frac{\pi}{4}}t+\overline{\sigma}g_1(t)\big)\,,
\end{equation}
where $r_k(\zeta,\widetilde\zeta)$ is holomorphic near $(0,0)$ and positively homogeneous of degree $k$. Here $g_1\in \mathbb{S}_{cl}^{-1}([0,+\infty[)$, $g_1\sim \frac{1}{2t}+\frac{c_{-2}}{t^2}+\ldots$, $t\longrightarrow +\infty$ and we get 
$$
 r_k(\sigma,\widetilde\zeta(\sigma,t))=\sum_{\nu+\mu=k} r_k^{\nu,\mu} \sigma^{\nu}\big(2s_+(\sigma) \e^{\mathrm{i}\frac{\pi}{4}}t+\overline{\sigma} g_1(t)\big)^{\mu} \in \mathbb{S}^k_{cl}([0,+\infty[)\,,
$$

\begin{align}\label{2.gt1}
 r_k(\sigma,\widetilde\zeta(\sigma,t)) &=\sum_{\nu+\mu+\rho=k} r_k^{\nu,\mu,\rho} t^{\mu} g_1(t)^{\rho}\\
 &\sim \sum_{\underset{\nu,\mu,\rho,\alpha\in\mathbb{N}}{\nu+\mu+\rho=k}}
   r_k^{\nu,\mu,\rho,\alpha} t^{\mu-\rho-2\alpha}\,,\qquad
   t\to +\infty\,,\nonumber
\end{align}
where we keep in mind that  $r_k^{\nu ,\mu ,..}$ depend on $\sigma$.  
This can be expressed as 
\begin{equation}\label{3.gt1}
 r_k(\sigma,\widetilde\zeta(\sigma,t))\sim \sum_{\mathbb{Z}\ni\mu\leq k}
 r_{k,\mu} t^\mu,\qquad t\longrightarrow +\infty\,.
\end{equation}
We apply this to \eqref{rc.6} and write
\begin{equation}\label{4.gt1}
 \prod_{j=1}^{n}\Big(h^{\frac{k_j}{2}-1}\varepsilon ^{k_j}r_{k_j}(\sigma,\widetilde\zeta(\sigma,t))\Big)=h^{\frac{\widetilde K}{2}-n}\varepsilon ^{\widetilde K} \widetilde r_{\widetilde K}(\sigma,t)\,,
\end{equation}
\begin{equation}\label{5.gt1}
 \widetilde {K}=\sum_{j=1}^n k_j,\qquad \widetilde r_{\widetilde K}(\sigma,t)\in \mathbb{S}^{\widetilde K}_{cl}([0,+\infty[)\,.
\end{equation}
In the non-trivial case $n\geq 1$, we have $1\leq n\leq [\frac{\widetilde K}{3}]$, where $[\frac{\widetilde K}{3}]$ denotes the integer part; the largest integer $\leq \frac{\widetilde K}{3}$.

\par
By \eqref{rc.2}, (\ref{rc.2,5}) the contribution from \eqref{4.gt1} to
the corresponding term in \eqref{rc.6} is $\mathcal{O}(1/\varepsilon \,    )\varepsilon ^{2\nu}h^{\frac{\widetilde K}{2}-[\frac{\widetilde K}{3}]}$ so we can discard the terms with $\frac{\widetilde K}{2}-[\frac{\widetilde K}{3}]\geq \frac{N}{2}-1$ and only keep the ones with $\frac{\widetilde K}{2}-[\frac{\widetilde K}{3}]< \frac{N}{2}-1$.

\par
Recall that $b$ is given in \eqref{2.sep1}, and assume for simplicity
that $a(\zeta,\widetilde\zeta)$ is independent of $h$, $|\zeta|$ is
assumed to be small and $\varepsilon \,    t=\mathcal{O}(h^{-\frac{1}{2}})$
in the region where $\chi=1$ so the exponential factor becomes
$\mathcal{O}\big(\e^{-1/\mathcal{O}(h)}\big)$ in the cutoff region
where $\varepsilon \,    t \asymp h^{-\frac{1}{2}}$. For this reason we
simply ignore the cutoff $\chi$.

\par
We have
$$
a(\zeta,\widetilde \zeta)=\sum_{j+k<N} a_{j,k} \zeta^j\widetilde\zeta^k+\mathcal{O}(1)\big(|\zeta|+|\widetilde\zeta|\big)^N,
$$
so
\begin{equation}\label{1.gt2}
\hskip-5pt
 a(\zeta,\widetilde \zeta(\zeta,t))=\sum_{j+k<N} a_{j,k} \sigma^j
 (h^{\frac{1}{2}}\varepsilon )^{j+k}\Big(2s_+\e^{\mathrm{i}\frac{\pi}{4}} t+\overline{\sigma}g_1(t)\Big)^k
 +\mathcal{O}(1) h^{\frac{N}{2}}\Big(\varepsilon (1+t)\Big)^N.
\end{equation}
The general term belongs to $(h^{\frac{1}{2}}\varepsilon )^{j+k}\mathbb{S}^k_{cl}([0,+\infty[)$ and is of the form $(h^{\frac{1}{2}}\varepsilon )^{j+k}c_{j,k}(t)$, where $c_{j,k}\in \mathbb{S}^k_{cl}([0,+\infty[)$.

\par
Regrouping the terms with the same $j+k$, we can write the double sum
in \eqref{1.gt2} as
\begin{equation}\label{2.gt2}
 \sum_{j=0}^{N-1} (h^{\frac{1}{2}}\varepsilon)^{j} c_j(t),\qquad c_j\in \mathbb{S}^j_{cl}([0,+\infty[)\,.
\end{equation}
\par
From \eqref{2.sep1}, (\ref{4.sep1}) it follows that
$(h^{\frac12}\varepsilon )^{-1}b(\zeta,t;h)$ has the same structure as in \eqref{1.gt2}, \eqref{2.gt2},
\begin{equation}\label{3.gt2}
\frac{b(\zeta,t)}{h^{\frac12}\varepsilon \,    }=\mathcal{O}(1)h^{\frac{N}{2}} \Big(\varepsilon (t+1)\Big)^N+\sum_{j=0}^{N-1}
(h^{\frac{1}{2}}\varepsilon )^j b_j(t)\,,
\end{equation}
\begin{equation}\label{4.gt2}
 b_j\in \mathbb{S}^j_{cl}([0,+\infty[)\,.
\end{equation}
Recall \eqref{rc.4}, where $\upsilon={\overset{2}{\underset{j=0}{\sum}}} \upsilon_j$ and
\eqref{II.36}-\eqref{II.38}, to see that $\widetilde \upsilon\in
\mathbb{S}^0_{cl}([0,+\infty[)$. Combine this with \eqref{3.gt2},
\eqref{4.gt2}, \eqref{4.gt1}. \eqref{rc.6}, where now
$(h^{\frac12}\varepsilon )^{-1}b={\mathcal O}(1)$, gives
\begin{equation}\label{1.gt3}\begin{split}
\frac{1}{h^{\frac12}\varepsilon} &\int_0^{+\infty} b(\zeta ,t;h)\e^{\frac{1}{h}\varphi (\zeta ,t)}\mathrm{d} t ={\mathcal O}(\frac{1}{\varepsilon })\left(h^{\frac{N}{2}}+\varepsilon ^{2L}\right) +\, \hbox{a linear
  combination }\\
&\hbox{ of terms }\int_0^{+\infty} \e^{\varepsilon ^2\varphi_2^0 (\sigma,t)}
\Big(\sum_{j=0}^N (h^{\frac{1}{2}}\varepsilon )^j b_j(t)\Big)        \big(\widetilde{\upsilon}(\sigma ,t)\varepsilon^2\big)^\nu h^{\frac{\widetilde K}{2}-n}
\varepsilon ^{\widetilde K}\widetilde r_{\widetilde K}(\sigma,t) \, \mathrm{d} t
\end{split}
\end{equation}
with $\widetilde K, N, \nu\in \mathbb{N}, \nu<L, \frac{\widetilde K}{2}-[\frac{\widetilde K}{3}]<\frac{N}{2}, n\leq [\frac{\widetilde K}{3}]$ and
\begin{equation}\label{2.gt3}
 b_j\in \mathbb{S}^j_{cl}([0,+\infty[),\quad \widetilde \upsilon\in \mathbb{S}^0_{cl}([0,+\infty[),\quad \widetilde r_{\widetilde K}(\sigma,t) \in \mathbb{S}^{\widetilde K}_{cl}([0,+\infty[)\,.
\end{equation}
\medskip\noindent
{\bf Analytic singularity at $\varepsilon \,    =0$. }
Rather than computing the full asymptotics of (\ref{1.gt3}) when
$\varepsilon \to 0$, we settle for terms that are non-analytic there.
Recall that $\displaystyle \int_0^{+\infty} \e^{-\frac{t^2}{2}}\, \mathrm{d} t=\sqrt{\frac{\pi}{2}}$. By scaling 
\begin{equation}\label{1.gt4}
 \int_0^{+\infty} \e^{-\varepsilon ^2 \frac{t^2}{2}}\, \mathrm{d} t=\frac{1}{\varepsilon }\sqrt{\frac{\pi}{2}},\qquad \varepsilon >0.
\end{equation}
Differentiate with respect to $\varepsilon $ :
\begin{align}\label{2.gt4}
 \int_0^{+\infty} \e^{-\varepsilon ^2 \frac{t^2}{2}} \varepsilon \,    t^2\, \mathrm{d} t & =
 \frac{1}{\varepsilon ^2}\sqrt{\frac{\pi}{2}}\,\nonumber\\
\int_0^{+\infty} \e^{-\varepsilon ^2 \frac{t^2}{2}} t^2\, \mathrm{d} t & =
\frac{c_2}{\varepsilon ^3}\,,\quad c_2=\sqrt{\frac{\pi}{2}}\not=0\,.
 \end{align}
Repeating the argument gives
\begin{equation}\label{3.gt4}
 \int_0^{+\infty} \e^{-\varepsilon ^2 \frac{t^2}{2}} t^{2k}\, \mathrm{d} t  =
\frac{c_{2k}}{\varepsilon ^{2k+1}}\,,\quad c_{2k}:=\frac{(2k)!}{2^k k!}\sqrt{\frac{\pi}{2}}\not=0 \,.
\end{equation}
Let $\chi\in C^\infty([0,+\infty[;[0,1])$, $\chi(t)=1$ for $t\gg 1$ and $\chi(t)=0$ for $0\leq t\ll 1$. Then 
\begin{equation}\label{4.gt4}
 \int_0^{+\infty} \e^{-\varepsilon ^2 \frac{t^2}{2}} \chi(t) t^{2k}\, \mathrm{d} t  =
\frac{c_{2k}}{\varepsilon ^{2k+1}} + \mathrm{an}(\varepsilon )\,,
\end{equation}
where $\mathrm{an}(\varepsilon )$ denotes a function, analytic near $\varepsilon =0$.

\par\medskip
In order to reach $k\leq 0$, we take $k=0$ in \eqref{4.gt4}:
$$
\int_0^{+\infty} \e^{-\varepsilon ^2 \frac{t^2}{2}} \chi(t)\, \mathrm{d} t=\frac{c_0}{\varepsilon }+\mathrm{an}(\varepsilon)\,, \quad c_0=\sqrt{\frac{\pi}{2}}\not=0\,,
$$
$$
\int_0^{+\infty} \e^{-\varepsilon ^2\frac{t^2}{2}} \varepsilon \,    t^2 \frac{1}{t^2}\chi(t)\, \mathrm{d} t={c_0}+\mathrm{an}(\varepsilon ).
$$
Take primitives with respect to $\varepsilon $ :
$$
\int_0^{+\infty} \e^{-\varepsilon ^2 \frac{t^2}{2}} \frac{1}{t^2}\chi(t)\, \mathrm{d} t=c_{-2}-{c_0}\varepsilon +\mathrm{an}(\varepsilon )=\mathrm{an}(\varepsilon)  
$$
with $c_{-2}=\int_0^{+\infty} \frac{\chi(t)}{t^2}\, \mathrm{d}t<+\infty$. Recall that $\chi^{(p)}(0)=0$ for all $p\in \mathbb{N}$. 

\noindent
Iterate:
\begin{equation}
\int_0^{+\infty} \e^{-\varepsilon ^2 \frac{t^2}{2}} \frac{1}{t^{2k}}\chi(t)\, \mathrm{d} t=\mathrm{an}(\varepsilon ),\qquad k=1,2,\ldots\,.
\end{equation}

\noindent
Next, start from 
\begin{equation}
 \int_0^{+\infty} \e^{-\frac{t^2}{2}} t\,\mathrm{d} t=1
\end{equation}
and scale 
\begin{equation}\label{1.gt5}
 \int_0^{+\infty} \e^{-\varepsilon ^2 \frac{t^2}{2}} t\,\mathrm{d} t=\frac{1}{\varepsilon ^2} \,.
\end{equation}
Differentiate with respect to $\varepsilon $ :
$$
\int_0^{+\infty} \e^{-\varepsilon ^2\frac{t^2}{2}} t^3\,\mathrm{d} t=\frac{c_3}{\varepsilon ^4}\quad \hbox{ with } \, c_3=2
$$
and more generally,
\begin{equation}\label{2.gt5}
 \int_0^{+\infty} \e^{-\varepsilon ^2 \frac{t^2}{2}} t^{2k+1}\,\mathrm{d} t=\frac{c_{2k+1}}{\varepsilon ^{2k+2}}\quad \hbox{ with } \, c_{2k+1}=2^k k!\,,\,\,    
 \hbox{ for all } \, k\in \mathbb{N}\,.
\end{equation}
To go in the opposite direction, start from
$$
\int_0^{+\infty} \e^{-\varepsilon ^2 \frac{t^2}{2}} \chi(t) t\,\mathrm{d} t=\frac{1}{\varepsilon ^2}+\mathrm{an}(\varepsilon )\,.
$$
$$
\int_0^{+\infty} \e^{-\varepsilon ^2 \frac{t^2}{2}}\varepsilon \,    t^2 \frac{1}{t} \chi(t) \,\mathrm{d} t=\frac{1}{\varepsilon }+\mathrm{an}(\varepsilon )
$$
and integrate with respect to $\varepsilon $ :
$$
\int_0^{+\infty} \e^{-\varepsilon ^2 \frac{t^2}{2}}\frac{1}{t} \chi(t) \,\mathrm{d} t=\ln\frac{1}{\varepsilon }+\mathrm{an}(\varepsilon )\,.
$$
Write this as 
$$
\int_0^{+\infty} \e^{-\varepsilon ^2 \frac{t^2}{2}}\varepsilon \,    t^2\frac{1}{t^3} \chi(t) \,\mathrm{d} t=\varepsilon \ln\frac{1}{\varepsilon }+\mathrm{an}(\varepsilon )
$$
and integrate with respect to $\varepsilon $ 
$$
\int_0^{+\infty} \e^{-\varepsilon ^2 \frac{t^2}{2}}\frac{1}{t^3} \chi(t) \,\mathrm{d} t=c_{-3}\varepsilon ^2\ln\frac{1}{\varepsilon }+\mathrm{an}(\varepsilon )\,.
$$
Iterate:
\begin{equation}\label{3.gt5}
 \int_0^{+\infty} \e^{-\varepsilon ^2 \frac{t^2}{2}}\frac{1}{t^k} \chi(t) \,\mathrm{d} t=c_{-k}\varepsilon ^{k-1}\ln\frac{1}{\varepsilon }+\mathrm{an}(\varepsilon )\,,
\end{equation}
for $1\leq k\in \mathbb{N}$ odd. Summing up, we have for all $k\in \mathbb{Z}$,
\begin{equation}\label{1.gt6}
\hskip-10pt\int_0^{+\infty} \e^{-\varepsilon ^2 \frac{t^2}{2}} \, {t^k}\chi(t)\,\mathrm{d} t=
 \begin{cases}
  c_k \varepsilon ^{-1-k}+\mathrm{an}(\varepsilon ),& \hbox { if }
  k\geq 0\,,\\
  \mathrm{an}(\varepsilon ), & \hbox { if }  0\not=k\in -2\mathbb{N} ,\\
  c_k\varepsilon ^{-1-k}\ln(\frac{1}{\varepsilon})+\mathrm{an}(\varepsilon )\,, &
  \hbox{ if }  k\in -(2\mathbb{N}+1)\,.\\
 \end{cases}
\end{equation}

\noindent
Using \eqref{1.gt6}, we compute for any fixed $\tau\in \mathbb{C}$ the following integral 
\begin{equation}\label{1,05.gt6}
\int_0^{+\infty} \e^{-\varepsilon ^2 \frac{t^2}{2}} \, (t+\tau)^k\chi(t)\,\mathrm{d}t\,.
\end{equation}
 In fact, we get 

\begin{itemize}
    
    \item For $k\in \mathbb{N}$,
    \begin{equation}\label{1,1.gt6}
       \begin{split}
\hskip-20pt \int_0^{+\infty} \e^{-\varepsilon ^2 \frac{t^2}{2}} \, (t+\tau)^k\chi(t)\,\mathrm{d} t&=\sum_{j=0}^k d_j^{(k)}\varepsilon^{-1-j}+ \mathrm{an}(\varepsilon )\\
       &=\frac{d_k^{(k)}}{\varepsilon^{k+1}}+\frac{d_{k-1}^{(k)}}{\varepsilon^{k}}+\cdots+\frac{d_0^{(k)}}{\varepsilon}+ \mathrm{an}(\varepsilon)\,,
       \end{split}
    \end{equation}
where for all $j\in\{0,\ldots, k\}\,,\quad\displaystyle d_j^{(k)}=\binom{k}{j} \tau^j\times  \begin{cases}
\frac{(2p)!}{2^p p!}\sqrt{\frac{\pi}{2}}  &\mathrm{if} \,\,\,  j=2p\\
{} & {}\\
2^p p!  &\mathrm{if}  \,\,\,  j=2p+1
\end{cases}$\,.

    \item  For $k\in \mathbb{Z}_-$, 
    
    \begin{align}\label{1,2.gt6}
       \int_0^{+\infty} \e^{-\varepsilon ^2 \frac{t^2}{2}} \, (t+\tau)^k\chi(t)\,\mathrm{d} t=\underset{-k+l \, \mathrm{ is \, odd}}{\sum_{l=0}^N} d_l^{(k)}&\varepsilon^{-1-k+l} \ln (\frac{1}{\varepsilon})+\\
       & + \mathcal{O}(\varepsilon^{-k+N-\delta})+ \mathrm{an}(\varepsilon)\,,\nonumber 
    \end{align}
    for any fixed integer $N$ and a positive real number $\delta$ small enough ($0<\delta\ll 1$). 
    Note that the coefficients $d_l^{(k)}$'s are computable.        
\end{itemize}

\noindent
\eqref{1,1.gt6} and \eqref{1,2.gt6} remain valid with modified coefficients, 
if we replace $\e^{-\varepsilon ^2 \frac{t^2}{2}}$ by $\e^{-\widetilde{\varepsilon}^2 \frac{t^2}{2}}$ where $\widetilde{\varepsilon}$ is a constant with $\mathrm{Re}( \widetilde{\varepsilon})>0$.

\smallskip
Notice that if $f=\mathcal{O}\Big((t+1)^{-N-2}\Big)$, then
\begin{equation}\label{2.gt6}
 F(\varepsilon ):=\int_0^{+\infty} \e^{-\varepsilon ^2 t^2} \, f(t)\, \mathrm{d} t \in C^N\big(\mathrm{neigh}(0);\mathbb{R}\big)
\end{equation}
and by Taylor expansion 
\begin{equation}\label{3.gt6}
 F(\varepsilon )=\sum_{\nu=0}^{N-1} (\partial_\varepsilon ^\nu F)(0) \,\,\frac{\varepsilon ^\nu}{\nu !}+\mathcal{O}(\varepsilon ^N)\,.
\end{equation}

We return to \eqref{1.gt3} and write the general term on the RHS as 
\begin{equation}\label{4.gt6}
\sum_{j=0}^N h^{\frac{j}{2}+\frac{\widetilde K}{2}-n}
\varepsilon ^{2\nu}\int_0^{+\infty}
\e^{\varepsilon ^2\varphi^0_2(\sigma,t)} \, \widetilde \upsilon(\sigma, t)^\nu
\varepsilon ^{ j+\widetilde K} b_j(t) \widetilde r_{\widetilde K}(\sigma,t)\,dt
\end{equation}
with $\widetilde K, N, \nu, \widetilde \upsilon, b_j, \widetilde
r_{\widetilde K}$ as
after (\ref{1.gt3}) and in \eqref{2.gt3}.

We split the integral in two by means of $1=(1-\chi)+\chi$ with $\chi$
as above. $(1-\chi)$ is supported in a bounded interval $[0,R[$ and we see that 
\begin{equation}\label{5.gt6}
 \int_0^{+\infty} \e^{\varepsilon ^2\varphi^0_2(\sigma,t)} \,
 \widetilde \upsilon(\sigma, t)^\nu \varepsilon ^{2\nu +j+\widetilde K} b_j(t) \widetilde r_{\widetilde K}(\sigma,t)(1-\chi(t))\,\mathrm{d} t=\mathrm{an}(\varepsilon )\,.
\end{equation}

\noindent
Recall that $\varphi^0_2(\sigma,t)$ is given by \eqref{4.sep2} with $\upsilon=0$. The results \eqref{1,1.gt6} \eqref{1,2.gt6} remain valid for 
\begin{equation}\label{6.gt6}
\int_0^{+\infty} \e^{\varepsilon ^2\varphi^0_2(\sigma,t)} t^k \chi(t)\, \mathrm{d} t.
\end{equation}
In fact, using a contour deformation to reach the critical point 
$$
\tau_c:=- \e^{-\mathrm{i}\frac{\pi}{4}}\frac{\lambda-\mu}{2s_+(\lambda+\mu+\mathrm{i} 2\rho)}\sigma
$$
of $t\mapsto\varphi^0_2(\sigma,t)$ and the fact that the function $\chi$ is equal to 1 for large $t$, we show that modulo an analytic function of $\varepsilon$ the integral in \eqref{6.gt6} is equal to that given by \eqref{1,05.gt6} with $\tau=\tau_c$, modified as mentioned after \eqref{1,2.gt6}. 

We get the same symbol space $\mathbb{S}^m_{cl}([T,+\infty[)$ for large $T,$ if we replace all powers of $t$ by powers of $\widetilde t:= 2 s_+ t+ \e^{-\mathrm{i}\frac{\pi}{4}}\dfrac{\lambda-\mu}{\lambda+\mu+\mathrm{i} 2\rho}\,\sigma$ and it easy to go from one representation to the other.

Next look at the complementary term in the splitting of the integral
in \eqref{4.gt6}, namely $h^{\frac{j}{2}+\frac{\widetilde{K}}{2}-n}$ times
\begin{equation}\label{1.gt7}
 \int_0^{+\infty} \e^{\varepsilon ^2\varphi^0_2(\sigma,t)} \, \widetilde
 \upsilon(\sigma, t)^\nu \varepsilon ^{2\nu +j+\widetilde K} b_j(t) \widetilde r_{\widetilde K}(\sigma,t)\chi(t)\,\mathrm{d} t\,.
\end{equation}
In view of the above remark we may shift to the variable $\widetilde t$, where $\varphi^0_2(\sigma,t)=-C\varepsilon ^2 \widetilde t^2$, $\mathrm{Re}( C)>0$,
then drop the tildes and reduce us (essentially) to the case 
$\e^{\varepsilon ^2\varphi^0_2(\sigma,t)}=\e^{-\varepsilon ^2 t^2}$.

\noindent
We have $\widetilde \upsilon^\nu b_j \widetilde r_{\widetilde K}\in \mathbb{S}_{cl}^{j+\widetilde K}([T,+\infty[)$, $\mathrm{supp}\chi\subset [T,+\infty[$ and for $M\in \mathbb{N}$, $M\gg 1$,
\begin{equation}\label{2.gt7}
 \widetilde \upsilon^\nu b_j \widetilde r_{\widetilde K}(t)=\sum_{\mu=-M-1}^{j+\widetilde K} d_\mu t^\mu +\mathcal{O}\Big((t+1)^{-M-2}\Big)\,.
\end{equation}
Insert the different terms of RHS\eqref{2.gt7} into \eqref{1.gt7}:
\begin{itemize}
 
 \item When $0\leq \mu\leq j+\widetilde K$, we get in view of \eqref{1,1.gt6}
\begin{equation}\label{3.gt7}
  \varepsilon ^{2\nu +j+\widetilde K-\mu}\int_0^{+\infty}
  \e^{\varepsilon ^2\varphi_2^0(\sigma,t)} \varepsilon ^\mu d_\mu
  t^{\mu}\chi(t)\, \mathrm{d} t=\varepsilon ^{2\nu +j+\widetilde K-\mu} \, \Big(\frac{\widetilde d_\mu}{\varepsilon }+\mathrm{an}(\varepsilon )\Big)\,.
\end{equation}

 \item When $\mu\leq -1$, we get 
 
 \begin{align}\label{4.gt7}
   \varepsilon ^{2\nu +j+\widetilde K}\int_0^{+\infty}
   &\e^{\varepsilon ^2\varphi_2^0(\sigma,t)}
   d_\mu t^{\mu}\chi(t)\, \mathrm{d} t =
   \varepsilon ^{2\nu +j+\widetilde K} \Big(\widetilde d_\mu \varepsilon ^{-\mu-1}\ln(\frac1\varepsilon ) \\
&  + \underset{-\mu+l \, \mathrm{ is \, odd }}{\sum_{l=0}^{N}} \widetilde d_l^{(\mu)}
    \varepsilon ^{-1-\mu +l} \ln(\frac1\varepsilon )+\mathcal{O}(\varepsilon^{-\mu+N-\delta})
    +\mathrm{an}(\varepsilon )\Big)\,.\nonumber
    \end{align}
 
 \end{itemize}

Finally from \eqref{2.gt6}, \eqref{3.gt6}, we get 
\begin{equation}\label{5.gt7}
  \varepsilon ^{2\nu +j+\widetilde K}\hskip-5pt\int_0^{+\infty}
  \hskip-5pt \e^{\varepsilon ^2\varphi_2^0(\sigma,t)}
  \mathcal{O}\left((t+1)^{-M-2}\right)\chi(t)\, \mathrm{d} t=\varepsilon ^{2\nu +j+\widetilde K}\left(\mathrm{an}(\varepsilon )+\mathcal{O}(\varepsilon ^M)\right)\,.
\end{equation}

\medskip 
Combining \eqref{1.gt7}-\eqref{5.gt7}, we get
\begin{align}\label{6.gt7}
  &\int_0^{+\infty} \e^{\varepsilon ^2\varphi^0_2(\sigma,t)} \widetilde
    \upsilon(\sigma, t)^\nu
    \varepsilon ^{2\nu +j+\widetilde K}  b_j(t) \widetilde r_{\widetilde K}(\sigma,t)\chi(t)\mathrm{d} t\\
  &= \hskip-10pt\sum_{0\leq \mu\leq j+\widetilde K} \hskip-10pt
    \varepsilon ^{2\nu +j+\widetilde K-\mu} \frac{\widetilde d_\mu}{\varepsilon }+ \hskip-10pt\sum_{-M-1\leq \mu\leq -1} \hskip-10pt \varepsilon ^{2\nu +j+\widetilde K} \widetilde d_\mu \varepsilon ^{-\mu-1}\ln(\frac1\varepsilon )
    +\varepsilon ^{2\nu +j+\widetilde K}\mathcal{O}(\varepsilon ^M)+\mathrm{an}(\varepsilon )\,.\nonumber
\end{align}
Use this and \eqref{5.gt6} in \eqref{1.gt3} with $M\geq 2L$,
\begin{multline}\label{1.gt8}
  \frac{1}{h^{\frac12}\varepsilon}\int_0^{+\infty} b(\zeta ,t;h)\, \e^{\frac{1}{h}\varphi (\zeta
    ,t)}\,\mathrm{d} t\\ ={\mathcal
    O}(\frac1\varepsilon )\left(h^{\frac{N}{2}-1}+\varepsilon ^{2L}\right)+
  \mathrm{an}(\varepsilon ) + \hbox{a linear
    combination of terms}\\
  \sum_{j=0}^{N} h^{\frac{j}{2}+\frac{\widetilde
      K}{2}-n}\varepsilon ^{2\nu}\left(\sum_{\mu =0}^{ j+\widetilde
      K} \varepsilon ^{j+\widetilde K-\mu} \frac{\widetilde
      d_\mu}{\varepsilon }+ \sum_{\mu=-M-1}^{ -1}
    \varepsilon ^{j+\widetilde K} \widetilde d_\mu
    \varepsilon ^{-\mu-1}\ln(\frac1\varepsilon )\right)
\end{multline}
with 
\begin{equation}\label{2.gt8}
 \widetilde K, N, \nu\in \mathbb{N},\ \nu<L, \ \frac{\widetilde K}{2}-\left[\widetilde K/3\right]<\frac{N}{2}-1,\ n\leq \left[\widetilde K/3\right]\,.
\end{equation}
Here $\widetilde d_\mu$ also depends on $j,\widetilde K,\nu$ and
$\sigma $. The sum in RHS\eqref{1.gt8} can be written
$$
\varepsilon ^{2\nu}\sum_{j=0}^N\left(h^{\frac{j}{2}+\frac{\widetilde K}{2}-n}p_j(\varepsilon )\frac{1}{\varepsilon }+\varepsilon ^{j+\widetilde K}h^{\frac{j}{2}+\frac{\widetilde K}{2}-n}q_j(\varepsilon )\ln(\frac1\varepsilon )\right)\,,
$$
where $p_j,q_j$, $j=1,\ldots,N$ are polynomials in $\varepsilon $  and
after extracting some powers of $h$, we get
\begin{align}\label{3.gt8}
  &\frac{1}{h^{\frac12}\varepsilon \,    }\int_0^{+\infty} b(\zeta ,t;h)\,\e^{\frac{1}{h}\varphi (\zeta
    ,t)}\,\mathrm{d} t ={\mathcal
    O}(\frac1\varepsilon )\left(h^{\frac{N}{2}-1}+\varepsilon ^{2L}\right)
  +\mathrm{an}(\varepsilon ) + \hbox{a linear
   } \nonumber\\
  &\hbox{ combination of terms } \varepsilon ^{2\nu} h^{\frac{\widetilde K}{2}-n} \sum_{j=0}^N
  h^{\frac{j}{2}}\left(p_j(\varepsilon )\frac{1}{\varepsilon }+\varepsilon ^{j+\widetilde
      K}q_j(\varepsilon )\ln(\frac1\varepsilon )\right)
\end{align}
still with \eqref{2.gt8} valid.

\par\medskip\noindent 
{\bf The leading term in (\ref{rc.6})}.\label{lt}
\par\noindent 
We insert (\ref{1.gt2}) into the general terms in RHS(\ref{rc.6})
and use (\ref{2.sep1}) where
$$
\frac{\overline{\widetilde{\zeta }}(\zeta
  ,t)}{\sqrt{1+t^2}}=
h^{\frac12}\varepsilon \frac{2s_+(\sigma)\e^{-\mathrm{i}\frac{\pi}{4}}t+\sigma g_1(t)}{\sqrt{1+t^2}}
=h^{\frac12}\varepsilon \Big(2s_+(\sigma )\e^{-\mathrm{i}\frac{\pi}{4}}+{\mathcal O}((1+t)^{-2})\Big).
$$
The leading contribution to (\ref{rc.6}) is then
\begin{equation}\label{1.lt1}
  \int_0^{+\infty} h^{\frac12}\varepsilon   2s_+(\sigma )\e^{-\mathrm{i}\frac{\pi }{4}}\e^{\varepsilon^2\varphi_2^0(\sigma ,t)}a(0,0) \, \mathrm{d}t={\mathcal O}(h^{\frac12}) \,,
\end{equation}
obtained by taking the leading term $a(0,0)$ in (\ref{1.gt2}), $\nu
=0$ and the trivial product $``\overset{0}{\underset{j=1}{\prod}}(...)"=1$. For all other
terms we gain some power of $\varepsilon $  or $h^{\frac12}$. Notice that the
remainder ${\mathcal O}(\frac1\varepsilon )(h^{N/2}+\varepsilon \,    ^{2L})$ is to be
replaced by ${\mathcal O}(h^{\frac12})(h^{N/2}+\varepsilon \,    ^{2L})$ since our
$b$ is now ${\mathcal O}(h^{\frac12}\varepsilon \,    )$ rather than ${\mathcal O}(1)$.

\medskip\noindent 
{\bf In conclusion}, for $b$ given in
(\ref{2.sep1}), we have after comparing the asymptotic decomposition
\eqref{1.gt2}, \eqref{3.gt2}, \eqref{1.gt3},
\begin{align}\label{2.lt1}
 \frac{1}{\mathrm{i}} \int_0^{+\infty} b(\zeta ,t;h)\e^{\frac{1}{h}\varphi (\zeta ,t)} \, \mathrm{d}t
  =\frac{1}{\mathrm{i}} \times & h^{\frac12}\varepsilon 2s_+(\sigma ) \e^{-\mathrm{i}\frac{\pi}{4}}a(0,0)\times \\
  &\int_0^{+\infty} \e^{\varepsilon^2\varphi_2^0(\sigma ,t)}\, \mathrm{d}t+{\mathcal O}(h)+{\mathcal O}(h^{\frac12}\varepsilon)+{\mathcal O}(\varepsilon^2)\,,\nonumber
  \end{align}
where the first term on the right hand side is ${\mathcal O}(h^{\frac12})$.

\noindent
Recall that $\varphi_2^0(\zeta ,t)$ is given by RHS(\ref{4.sep2})
with "$\upsilon$" there replaced by $0$. 
Remembering the treatment of \eqref{6.gt6} with $k=0$ we see that
\begin{equation}\label{3.lt1}
  \begin{split}
\int_0^{+\infty} \e^{\varepsilon^2\varphi_2^0(\sigma
  ,t)} \, \mathrm{d}t&=\mathrm{an}(\varepsilon)+\int_0^{+\infty} \exp
\left(-\varepsilon^2\frac{\lambda +\mu +\mathrm{i} 2\rho }{8}(2s_+t)^2
\right) \, \mathrm{d}t\\
&=\mathrm{an}(\varepsilon)+\int_0^{+\infty} \exp
\left(-\varepsilon^2s_+^2(\lambda +\mu +\mathrm{i} 2\rho
  )\frac{t^2}{2}\right) \, \mathrm{d}t\\
  &=\mathrm{an(\varepsilon )}+\frac{\sqrt{2\pi}}{2\varepsilon 
    |s_+|(\lambda +\mu +\mathrm{i} 2\rho )^{\frac12}} \,.
  \end{split}
\end{equation}

\noindent
Inserting this in (\ref{2.lt1}) gives
\begin{equation}\label{4.lt1}\begin{split}
\frac{1}{\mathrm{i}}  \int_0^{+\infty} & b(\zeta ,t;h) \, \e^{\frac{1}{h}\varphi (\zeta ,t)} \, \mathrm{d}t\\
  &=\frac{h^{\frac12}}{\mathrm{i}}
\frac{s_+(\sigma )}{|s_+(\sigma )|}
\frac{\sqrt{2\pi }\, a(0,0) \e^{-\mathrm{i}\frac{\pi}{4}}}{(\lambda +\mu +\mathrm{i} 2\rho)^{\frac12}}
+{\mathcal O}(h^{\frac12}\varepsilon)+{\mathcal O}(h)+{\mathcal O}(\varepsilon^2)\\
&=-h^{\frac12}
\frac{s_+(\sigma )}{|s_+(\sigma )|}
\frac{\sqrt{2\pi }\, a(0,0) \e^{\mathrm{i}\frac{\pi}{4}}}{(\lambda +\mu +\mathrm{i} 2\rho)^{\frac12}}
+{\mathcal O}(h^{\frac12}\varepsilon)+{\mathcal O}(h)+{\mathcal O}(\varepsilon^2)\,,
  \end{split}
\end{equation}
where $s_+(\sigma)=\mathrm{Re}\big(\e^{\mathrm{i}\frac{\pi}{4}}\sigma\big)$ with $\sigma=\frac{\zeta}{|\zeta|}\not=0$.
Since $\mathrm{Re}\big(\e^{\mathrm{i}\frac{\pi}{4}}\zeta\big)=-
\mathrm{Im}\big(\e^{-\mathrm{i}\frac{\pi}{4}}\zeta\big)$ then 
\begin{equation*}
\begin{split}
 \frac{s_+(\sigma )}{|s_+(\sigma )|}=\mathrm{sgn}\big(\mathrm{Re}(\e^{\mathrm{i}\frac{\pi}{4}}\zeta)\big)&=
\begin{cases}
                                                                                                -1 \, &\mathrm{if} \quad \mathrm{Re}\big(\e^{\mathrm{i}\frac{\pi}{4}}\zeta\big)<0\\
1 \,  &\mathrm{if} \quad \mathrm{Re}\big(\e^{\mathrm{i}\frac{\pi}{4}}\zeta\big)>0                                                                                                                \end{cases}\\
&=\begin{cases}
                                                                                                -1 \, &\mathrm{if} \quad \mathrm{Im}\big(\e^{-\mathrm{i}\frac{\pi}{4}}\zeta\big)>0\\
1 \,  &\mathrm{if} \quad \mathrm{Im}\big(\e^{-\mathrm{i}\frac{\pi}{4}}\zeta\big)<0                                                                                                                \end{cases}\,,
\end{split}
\end{equation*}
hence 
\begin{equation*}
-\frac{s_+(\sigma )}{2|s_+(\sigma )|}
=\begin{cases}
                                                                                    \frac{1}{2} \, &\mathrm{if} \quad \mathrm{Im}\big(\e^{-\mathrm{i}\frac{\pi}{4}}\zeta\big)>0\\
-\frac{1}{2} \,  &\mathrm{if} \quad \mathrm{Im}\big(\e^{-\mathrm{i}\frac{\pi}{4}}\zeta\big)<0                                                                                                                \end{cases}\,.
\end{equation*}

Therefore \eqref{4.lt1} becomes 
\begin{equation}\label{5.lt1}\begin{split}
 \frac{1}{\mathrm{i}}&  \int_0^{+\infty}  b(\zeta ,t;h) \, \e^{\frac{1}{h}\varphi (\zeta ,t)} \, \mathrm{d}t\\
 &=\Big(\hskip-5pt-\frac{s_+(\sigma )}{2|s_+(\sigma )|}\Big)\times h^{\frac12}
 \frac{2\sqrt{2\pi }\, a(0,0) \e^{\mathrm{i}\frac{\pi}{4}}}{(\lambda +\mu +\mathrm{i} 2\rho)^{\frac12}}
 +{\mathcal O}(h^{\frac12}\varepsilon)+{\mathcal O}(h)+{\mathcal O}(\varepsilon^2)\\
 &=\left[h^{\frac12}
 \frac{2\sqrt{2\pi }\, a(0,0) \e^{\mathrm{i}\frac{\pi}{4}}}{(\lambda +\mu +\mathrm{i} 2\rho)^{\frac12}}\times
 \begin{cases}
  \frac{1}{2} \, &\mathrm{if} \quad \mathrm{Im}\big(\e^{-\mathrm{i}\frac{\pi}{4}}\zeta\big)>0\\
-\frac{1}{2} \,  &\mathrm{if} \quad \mathrm{Im}\big(\e^{-\mathrm{i}\frac{\pi}{4}}\zeta\big)<0  
 \end{cases}\right]+R(h,\varepsilon)\,,
   \end{split}
\end{equation}
where the remainder term $R(h,\varepsilon)$ is ${\mathcal O}(h^{\frac12}\varepsilon)+{\mathcal O}(h)+{\mathcal O}(\varepsilon^2)$.

\subsection{Summary of this section}\label{sum}
Assume $\zeta \ne 0$. After a simultaneous rotation in ${\mathbb C}_\zeta $ and ${\mathbb{C}}_{\xi +i\eta }$, (which amounts to changing the
family of contours $\Gamma_\theta$,)  we may assume that $\zeta $
does not belong to a conic neighborhood of $\e^{\mathrm{i}\frac{\pi}{4}}{\mathbb{R}}$ in
${\mathbb{C}}\setminus \{0 \}$. In Proposition \ref{Inte.1} we have seen
that $\mathrm{II}(\zeta ,1)=\underset{\delta \to 0}{\lim}\mathrm{II}(\zeta
,1-\delta )$ exists and is given by the second integral in
(\ref{II.22}), rewritten in Proposition \ref{Inte.3} as
$\frac{1}{\mathrm{i}}\int_0^{+\infty}  b(\zeta ,t)\e^{\frac{|\zeta |^2}{h}\varphi(\sigma
 ,t)} \, \mathrm{d}t$. Here $\widetilde{\zeta }(t,\zeta )$ is given by
(\ref{II.20}), where $|s_+|\ge 1/{\mathcal O}(1)$.

\par 
Set $\varepsilon :=\frac{|\zeta |}{\sqrt{h}}$. When $\varepsilon \ge h^{-\delta
}$ for some fixed $\delta >0$, we have a two term asymptotic expansion of the original integral in Section \ref{le} and the Stokes reduction
to $\mathrm{I}+\mathrm{II}+\mathrm{III}$ is not needed. The
intermediate case is when $h^{\delta }\le \varepsilon \le h^{-\delta }$
for some small and fixed $\delta >0$. In this case Proposition
\ref{1nonq3} gives an asymptotic formula for $\int_T^{+\infty} b(\zeta
,t)\, \e^{\frac{|\zeta |^2}{h}\varphi (\sigma ,t)} \, \mathrm{d}t$ when 
$T\gg \max (1,\varepsilon^{-1})$, assuming that $\chi $ in (\ref{7.ft7})  has
  sufficiently small support:
  For $T\ge T_0\gg 1$, $T\gg \varepsilon^{-1}$, we have for every $b\in
|\zeta |\mathbb{S}^0([0,+\infty [)$ and $N\in \mathbb{N}$:
 \begin{align*}
   \int_T^{+\infty} b(t) \e^{\frac{|\zeta|^2}{h}\varphi (\sigma ,t)}\,
   \mathrm{d}t=-\e^{\frac{|\zeta|^2}{h}\varphi (\sigma ,T)}& \,\,\times\\
   &\left(c^{(N)}(T)+\mathcal{O}(1)\left(\frac{h}{(T|\zeta|)^2}\right)^{N+1}\frac{h}{T|\zeta|}\right)\,.\nonumber
 \end{align*}
 Here
 $$
 c^{(N)}=\sum_{k=0}^{N}\Big(-\frac{h}{|\zeta|^2\partial_t\varphi}\partial_t\Big)^{k}
 \frac{h}{|\zeta|^2\partial_t\varphi } b =: \sum_{k=0}^{N} c_k\,.
 $$
 We know that $\partial _t\varphi \asymp t$ for large $t$ and we can
 view (\ref{7.nonq3}) as an asymptotic formula in increasing powers of
 $(T\varepsilon)^{-1}$.

 \par 
 We now turn to the case when $\varepsilon $ is small in the sense
 that $\varepsilon \le h^{\delta }$ for fixed small $\delta >0$. It
 turns out that $\varepsilon $ can be used as a small asymptotic
 parameter and $\frac{1}{\varepsilon h^{1/2}}\int_0^{+\infty}
 b \, \e^{\frac{1}{h}\varphi (\zeta ,t)}\, \mathrm{d} t$  can first be developed as in (\ref{rc.6}), where "$b$" corresponds to our "$\frac{1}{\varepsilon
 h^{\frac12}} \, b$" and further computations lead to (\ref{1.gt8}),
 \begin{multline*}
  \frac{1}{h^{\frac12}\varepsilon}\int_0^{+\infty} b(\zeta  ,t;h)\, \e^{\frac{1}{h}\varphi (\zeta
    ,t)}\,\mathrm{d}t\\ ={\mathcal O}(\frac1\varepsilon )\left(h^{\frac{N}{2}-1}+\varepsilon ^{2L}\right)+
  \mathrm{an}(\varepsilon ) + \hbox{a linear
    combination of terms}\\
  \sum_{j=0}^{K} h^{\frac{j}{2}+\frac{\widetilde
      K}{2}-n}\varepsilon ^{2\nu}\left(\sum_{\mu =0}^{ j+\widetilde
      K} \varepsilon ^{j+\widetilde K-\mu} \frac{\widetilde
      d_\mu}{\varepsilon }+ \sum_{\mu=-M-1}^{ -1}
    \varepsilon ^{j+\widetilde K} \widetilde d_\mu
    \varepsilon ^{-\mu-1}\ln(\frac1\varepsilon )\right)
\end{multline*}
with (\ref{2.gt8}):
$$
 K,\widetilde K, N, \nu\in \mathbb{N},\ \nu<L, \ \frac{\widetilde K}{2}-\left[\widetilde K/3\right]<\frac{N}{2}-1,\ n\leq \left[\widetilde K/3\right]\,.
$$
Here $\widetilde d_\mu$ also depends on $j,\widetilde K,\nu$ and
$\sigma $.

\par This can also be written as
\begin{align*}
  \frac{1}{h^{\frac12}\varepsilon \,    }\int_0^{+\infty} b(\zeta ,t;h)\,&\e^{\frac{1}{h}\varphi (\zeta,t)}\,\mathrm{d}t ={\mathcal
    O}(\frac1\varepsilon )\left(h^{\frac{N}{2}-1}+\varepsilon ^{2L}\right)
  +\mathrm{an}(\varepsilon ) + \hbox{a linear
    combination} \\
  &\hbox{of terms } \varepsilon ^{2\nu} h^{\frac{\widetilde K}{2}-n} \sum_{j=0}^K
  h^{\frac{j}{2}}\left(p_j(\varepsilon )\frac{1}{\varepsilon }+\varepsilon ^{j+\widetilde
      K}q_j(\varepsilon )\ln(\frac1\varepsilon )\right)\,,
\end{align*}
where $p_j$, $q_j$ are polynomials in $\varepsilon $.
still with (\ref{2.gt8}) valid.

From these considerations together with Proposition \ref{Inte.3}, one obtains the leading-order contribution to the second term in the decomposition of Section \ref{Sto}, namely $\mathrm{II}(\zeta,1)$, given by (\ref{5.lt1}):
\begin{equation}\begin{split}\label{6.lt1}
\mathrm{II}(\zeta,1)&=\frac{1}{\mathrm{i}} \int_0^{+\infty}  b(\zeta ,t;h) \, \e^{\frac{1}{h}\varphi (\zeta ,t)} \, \mathrm{d}t\\
 &=\left[h^{\frac12}
 \frac{2\sqrt{2\pi }\, a(0,0) \e^{\mathrm{i}\frac{\pi}{4}}}{(\lambda +\mu +\mathrm{i} 2\rho)^{\frac12}}\times
 \begin{cases}
  \frac{1}{2} \, &\mathrm{if} \quad \mathrm{Im}\big(\e^{-\mathrm{i}\frac{\pi}{4}}\zeta\big)>0\\
-\frac{1}{2} \,  &\mathrm{if} \quad \mathrm{Im}\big(\e^{-\mathrm{i}\frac{\pi}{4}}\zeta\big)<0  
 \end{cases}\right]+R(h,\varepsilon)\,,
   \end{split}
\end{equation}
where the remainder term $R(h,\varepsilon)$ is ${\mathcal O}(h^{\frac12}\varepsilon)+{\mathcal O}(h)+{\mathcal O}(\varepsilon^2)$.

\section{The term \texorpdfstring{$\mathrm{III}(\zeta,1-\delta)$}{} in (\ref{2.sto8})}\label{tt}
\setcounter{equation}{0}

A discussion of the third term in the decomposition 
of Section \ref{Sto} is presented in this section. 

As in Proposition \ref{1sto7} we assume that $a$, $\phi $ are
holomorphic in $\mathrm{neigh}((0,0);\mathbb{C}^2)$ and that $\chi \in
C_0^\infty (\mathbb{C}^2)$ has its support in such a neighborhood. 

\noindent
We start by looking at the deformation vector field associated to the
coordinates $w_+$, $w_-$, associated to the $\theta $-dependent basis
$e_+$, $e_-$ in Lemma \ref{1ft1} and in (\ref{2.ft1}). Write
(\ref{1.ft1}) as
\begin{equation}\label{tt.1}
  \begin{split}
&e_+=\left(\check{f}(\theta )\e^{-\mathrm{i}\frac{\pi}{4}},f(\theta )\e^{\mathrm{i}\frac{\pi}{4}}\right)\in \Gamma_\theta ,\\
&e_-=\left(f(\theta )\e^{\mathrm{i}\frac{\pi}{4}},\check{f}(\theta )\e^{-\mathrm{i}\frac{\pi}{4}}\right)\in \Gamma_\theta ,
  \end{split}
\end{equation}
where 
\begin{equation}\label{tt.2}
f(\theta )=\frac{1+\theta }{\sqrt{2(1+\theta ^2)}},\quad \check{f}(\theta
):=f(-\theta)\,.
\end{equation}
Then a calculation shows that
\begin{equation}\label{tt.3}
\partial_\theta f(\theta )=\frac{1}{1+\theta  ^2}\check{f}(\theta
),\ \ \partial_\theta\check{f}(\theta )=-\frac{1}{1+\theta
  ^2}f(\theta )\,.
\end{equation}
(Note that $f(\theta )^2+\check{f}(\theta )^2=1$, so $f(\theta )=\sin
\alpha $, $\check{f}(\theta )=\cos \alpha $, for a suitable $\alpha
=\alpha (\theta )$.) A simple computation, using this and the fact that $\e^{\mathrm{i}\frac{\pi}
{4}}=\mathrm{i} \e^{-\mathrm{i}\frac{\pi}{4}}$ show that
\begin{equation}\label{tt.4}\begin{cases}
    \dot{e}_+(\theta )=\dfrac{\mathrm{i}}{1+\theta ^2}e_-(\theta )\\
\dot{e}_-(\theta )=\dfrac{\mathrm{i}}{1+\theta ^2}e_+(\theta )
  \end{cases}\,,
\end{equation}
where $\dot{e}_\pm=\partial_\theta e_\pm$\,.

\par\medskip
The deformation vector field at a point
$w_+e_+(\theta)+w_-e_-(\theta)\in \Gamma_\theta$ is given by the
vector
\begin{align*}
 \widetilde{\nu }:&=w_{+}\dot e_{+}+w_{-}\dot e_{-}\\
 & = \frac{1}{1+\theta^2}\big[w_{-}\,\mathrm{i} e_{+}+w_{+}\,\mathrm{i} e_{-}\big]\,.
\end{align*}
 If we identify the vector $\mathrm{i} e_\pm$ with the real vector field 
 $$
 2\mathrm{Re}(\mathrm{i}\partial_{w_\pm})=
 \mathrm{i}\partial_{w_\pm}+\overline{\mathrm{i}\partial_{w_\pm}}=\mathrm{i}\partial_{w_\pm}-\mathrm{i}\partial_{\overline{w}_\pm}\,,
 $$
 we obtain the deformation field
 \begin{equation}\label{tt.5}
   \widetilde{\nu
   }\simeq\frac{1}{1+\theta^2}\sum_{\bullet\in\{+,-\}}2\mathrm{Re}
   (\mathrm{i}\,w_{\widehat{\bullet}}\partial_{\bullet}),\hbox{ where }  
   \widehat{+}:=-,\ \ \widehat{-}:=+\,.
 \end{equation}
 Recall that the coordinates $w_{+}, w_{-}$ are real on $\Gamma_\theta.$
 

\medskip\par
We now turn to $\mathrm{III}(\zeta ,1-\delta)$ in (\ref{2.sto8})
which can be written
\begin{equation}\label{tt.7}
  \mathrm{III}(\zeta ,1-\delta)=\int_0^{1-\delta}k(\theta
  )\,\mathrm{d}\theta\,,
\end{equation}
where 
\begin{align}\label{tt.8}
k(\theta)=\frac{1}{\pi}\iintJ_{\Gamma_{\theta}}
\frac{1}{\omega-\zeta}\,   a(\omega,\widetilde \omega;h)\, \e^{\frac{\mathrm{i}}{h}\phi(\omega,\widetilde \omega)}&\times
\\
\nu(\theta;\omega,\widetilde\omega)&\mathord{\rfloor}
\mathrm{d}_{\omega,\tilde\omega}\left(
\chi(\omega,\widetilde\omega)\, \frac{\mathrm{d}\widetilde\omega\wedge \mathrm{d}\omega}{2\mathrm{i}}\right)\,,\nonumber
\end{align}
where $\nu$ is given in (\ref{Sto.7}). Now $\nu $ and $\widetilde{\nu
}$ in (\ref{tt.5}) are both deformation fields associated to the
family $\Gamma_\theta $ which however preserve different coordinates
on $\Gamma_\theta $. Thus $\nu -\widetilde{\nu }$ is
tangential to each $\Gamma_\theta $ and as explained before and in
Lemma \ref{Err} we can replace $\nu $ in the expression for $J(\theta
)$ with $\widetilde{\nu }$ and get
\begin{align}\label{tt.9}
 k(\theta)=\frac{1}{\pi}\iintJ_{\Gamma_{\theta}}
  \frac{1}{\omega-\zeta}  a(\omega,\widetilde \omega;h)\, &\e^{\frac{\mathrm{i}}{h}\phi(\omega,\widetilde \omega) }  \times\\
 &\widetilde{\nu}(\theta;\omega,\widetilde\omega)\mathord{\rfloor}
 \mathrm{d}_{\omega,\tilde\omega}\left(\chi(\omega,\widetilde\omega)\frac{\mathrm{d}\widetilde\omega\wedge \mathrm{d}\omega}{2\mathrm{i}}\right)\,.\nonumber
\end{align}
Here we recall that $w_\pm$ are defined as the coordinates on $\Gamma
_\theta$ in (\ref{2.ft1}) and to each $\zeta $ we associate the unique
$\widetilde{\zeta }$ with $(\zeta ,\widetilde{\zeta })\in \Gamma
_\theta $ so that

\begin{equation}\label{tt.10}
(\zeta ,\widetilde{\zeta })=z_+e_++z_-e_-, \quad z_\pm\in \mathbb{R}\,.
\end{equation}

\medskip\noindent
(This is different from (\ref{7.ft1})). Then as in (\ref{4.ft1})
\begin{equation}\label{tt.10,5}
z_+=\frac{\sqrt{2(1+\theta ^2)}}{1-\theta }\mathrm{Re} (\zeta \e^{\mathrm{i}\frac{\pi}{4}}),\quad
z_-=\frac{\sqrt{2(1+\theta^2)}}{1+\theta}\mathrm{Re} (\zeta \e^{-\mathrm{i}\frac{\pi}{4}})\,.
\end{equation}

\medskip\noindent
Using (\ref{3.ft1}), (\ref{tt.10,5}) and (\ref{6.ft1}), we obtain
\[\begin{split}
    k(\theta)=\frac{1}{\pi}&\iintJ\frac{a\big(\omega(w_{+}, w_{-}),\widetilde\omega(w_{+}, w_{-}) ;h\big) \e^{\frac{\mathrm{i}}{h}\phi\big(\omega(w_{+}, w_{-}),\widetilde\omega(w_{+}, w_{-})\big)}}{\dfrac{1-\theta}{\sqrt{2(1+\theta^2)}}\, \e^{-\mathrm{i}\frac{\pi}{4}} (w_{+}-z_{+}) +\dfrac{1+\theta}{\sqrt{2(1+\theta^2)}} \e^{\mathrm{i}\frac{\pi}{4}}(w_{-}-z_{-})}\times\\
    \frac{1}{2(1+\theta^2)}&\left(\sum_{\bullet\in\{+,-\}}2\mathrm{Re}(w_{\widehat{\bullet}}\,
      \mathrm{i}\partial_{w_\bullet})\right)\mathord{\rfloor}\left(\left(\partial_{\overline{w}_{+}}\chi )
    \mathrm{d}\overline{w}_{+}+ (\partial_{\overline{w}_{-}}\chi )
    \mathrm{d}\overline{w}_{-}\right){\wedge}\mathrm{d}w_{+}\wedge \mathrm{d}w_{-}\right),
\end{split}
\]
which can be rewritten as
 \begin{equation}\begin{split}\label{tt.11}
     k(\theta)&=\frac{1}{\pi\sqrt{2(1+\theta^2)}}\iintJ_{\mathbb{R}^2_{w_+,w_-}}
     \hskip-25pt\frac{\e^{\mathrm{i}\frac{\pi}{4}}a\big(\omega(w_{+},
       w_{-}),\widetilde\omega(w_{+}, w_{-}) ;h\big)
       \e^{\frac{\mathrm{i}}{h}\phi\big(\omega(w_{+},
         w_{-}),\widetilde\omega(w_{+}, w_{-})\big)}}
     {\big[(1-\theta) (w_{+}-z_{+})
       +\mathrm{i}(1+\theta)(w_{-}-z_{-})\big]}\times\\
       {} & {}\\
&\big(\sum_{\bullet\in\{+,-\}}2\mathrm{Re}(w_{\widehat{\bullet}}\,
\mathrm{i}\partial_{w_\bullet}) 
\big)
\mathord{\rfloor}\left((\partial_{\overline{w}_+}\chi)  \mathrm{d}\overline{w}_{+}\wedge \mathrm{d}w_{+}\wedge \mathrm{d}w_{-}
 + (\partial_{\overline{w}_{-}}\chi)  \mathrm{d}\overline{w}_{-}
 {\wedge} \mathrm{d}w_{+}\wedge \mathrm{d}w_{-}\right) .
 \end{split}
 \end{equation}
\par 
This can be expanded explicitly, and for our estimates we retain that
 \begin{equation}\label{tt.12}\begin{split}
  k(\theta)= \mathcal{O}(1)\iintJ_{\mathbb{R}^2_{w_+,w_-}}\hskip-20pt&\frac{\Big|a\big(\omega(w_{+}, w_{-}),\widetilde\omega(w_{+}, w_{-}) ;h\big) \e^{\frac{\mathrm{i}}{h}\phi\big(\omega(w_{+}, w_{-}),\widetilde\omega(w_{+}, w_{-})\big)}\Big| }{\big|(1-\theta) (w_{+}-z_{+}) +\mathrm{i}(1+\theta)(w_{-}-z_{-})\big|}\times\\
  &\times |\overline{\partial_{\omega,\tilde\omega}}\chi|
  \,\mathrm{d}w_{+}\mathrm{d}w_{-}\,.
 \end{split}
 \end{equation}

 \noindent
Here we recall that $\chi$ denotes an almost holomorphic extension of
 the original $\chi$ on $\Gamma_0$. We can choose $\chi$ so that
 $\overline{\partial}\chi$ vanishes on
 $\mathrm{neigh}\big((0,0);\mathbb{C}^2\big)$. More precisely, we can
 assume that $\overline{\partial}\chi$, viewed as a function on
 $\Gamma_\theta $, has its support in a set
\begin{equation}\label{tt.13}
 K_\theta=([-b,b]\times [-b,b])\setminus ([-a,a]\times [-a,a])
\subset\mathbb{R}^2_{w_+,w_-}\simeq \Gamma_\theta  
\end{equation}
where $0<a<b<\infty$ are independent of $\theta $. The almost 
holomorphicity tells us that 
\begin{equation}\label{tt.14}
 \overline{\partial}\chi=\mathcal{O}(\theta^\infty)\quad\mathrm{on}\,\,\, K_\theta\,.
\end{equation}

\par 
We assume that $\phi $ satisfies (\ref{6.ft3}), (\ref{1.ft4},
(\ref{2.ft4})) and (\ref{3.ft4}). Then by (\ref{6.ft4}), we have
$$
\mathrm{Im}\, \phi_2\asymp \theta (|\xi |^2+|\eta |^2)\hbox{ on }\Gamma_\theta\,.
$$
Assuming that $\phi$ is real valued on $\mathbb{R}^2_{\xi ,\eta }$ and
that we work in a sufficiently small neighborhood of $(0,0)\in {\mathbb{C}}^2$ (where $\chi $ has its support), we have
\begin{equation}\label{tt.15}
\mathrm{Im}\, \phi \asymp \theta (|\xi |^2+|\eta |^2)
\hbox{ on } \Gamma_\theta \cap\mathrm{supp\,}(\chi )\,.
\end{equation}
In particular,
\begin{equation}\label{tt.16}
 \mathrm{Im}\,\phi\asymp \theta\,\,\hbox{ on }\,\, K_\theta\,.
\end{equation}

\par 
Decompose $K_\theta$ into a disjoint union of rectangles $R$ of
horizontal size (i.e. in the $w_+$-direction) $\asymp 1$ and of
vertical size (i.e. in the $w_-$-direction) $\asymp 1-\theta$.
 Let $(w_+(R),w_-(R))$ be the center of $R$. Restrict the integral in
 the right hand side of \eqref{tt.12} to such a rectangle
 $R$:

 \begin{equation}\label{tt.17}
 k_R(\theta):=\iintJ_{R}\frac{ \mathcal{O}(\theta^\infty)\, \e^{-\frac{\theta}{Ch}}}{\big|(1-\theta) (w_{+}-z_{+}) +\mathrm{i}(1+\theta)(w_{-}-z_{-})\big|}\,\mathrm{d}w_{+}\mathrm{d}w_{-}\,.
 \end{equation}
 Rescale in the vertical direction:
 \begin{equation}\label{tt.18}
  (1-\theta)\,\widetilde{z}_{-}:=(1+\theta)\,z_{-},\quad \mathrm{and}\quad
  (1-\theta)\,\widetilde{w}_{-}:=(1+\theta)\,w_{-}\,.
 \end{equation}
 The corresponding points $(w_{+},\widetilde{w}_{-})$ form a rectangle $\widetilde R$ of horizontal and vertical size $\asymp 1$ and with center 
$ (w_{+}(R),\widetilde{w}_-(R))=(w_{+}(R),\frac{1+\theta}{1-\theta}{w}_-(R))$. We get 
\begin{equation}\label{tt.19}\begin{split}
 k_R(\theta) & =\iintJ_{\widetilde R(R)}\frac{\mathcal{O}(\theta^\infty)\, \e^{-\frac{\theta}{Ch}}}{(1-\theta)\big|(w_{+}-z_{+}) +\mathrm{i}(\widetilde{w}_{-}-\widetilde{z}_{-})\big|}\cdot \frac{1-\theta}{1+\theta}\,\mathrm{d}w_{+}\mathrm{d}\widetilde{w}_{-}\\
 &= \frac{\mathcal{O}(\theta^\infty)\,
   \e^{-\frac{\theta}{Ch}}}{1+|z_{+}-w_{+}(R)|+|\widetilde{z}_{-}-\widetilde{w}_{-}(\widetilde
   R)|}\,.
 \end{split}
\end{equation}
Here $w_{+}(R),$ $\widetilde \omega_-(R)=\mathcal{O}\Big(\dfrac{1}{1-\theta}\Big)$.


\medskip 
Very roughly, the numbers
$|{z}_{+}-{w}_{+}(R)|$ and $|\widetilde{z}_{-}-\widetilde{w}_{-}(\widetilde
R)|$ are distributed like the integers in an interval $[0,\simeq \mathcal{O}(1)]$,
$[0,\simeq (1-\theta)^{-1}]$ respectively, so when summing the $k_R(\theta)$ over
the rectangles $R$ in our decomposition, we get from \eqref{tt.19}:
\begin{equation}\label{tt.21}
k(\theta)  =\mathcal{O}(\theta^\infty)\, \e^{-\frac{\theta}{Ch}}\sum_{k=0}^{(1-\theta)^{-1}}\frac{1}{1+k}\,.
\end{equation}
The last sum in \eqref{tt.21} is $\mathcal{O}(1)(1+\ln\frac{1}{1-\theta})$ and we get from \eqref{tt.21}:
\begin{equation}\label{tt.22}
 k(\theta)=\mathcal{O}(\theta^\infty)\, \e^{-\frac{\theta}{Ch}} \left(1+\ln\frac{1}{1-\theta}\right)\,.
\end{equation}
Using this in (\ref{tt.7}) gives
\begin{equation}\label{tt.23}
  \mathrm{III}(\zeta ,1-\delta)={\mathcal  O}(h^\infty)\,.
\end{equation}
We also see that 
\begin{equation}\label{tt.24}
\mathrm{III}(\zeta ,1-\delta )\longrightarrow
\mathrm{III}(\zeta ,1):=\int_0^1 k(\theta )\, \mathrm{d}\theta \in \mathbb{C}\,,
\end{equation}

when $\delta \rightarrow 0$ and $\mathrm{III}(\zeta ,1)={\mathcal  O}(h^\infty)$.

\section{Conclusion}
In this paper, we have presented a detailed discussion of the 
integral (\ref{int.1}) for values of $0<h\ll 1$ with a stationary 
point of the phase $\phi(\omega)$ at $\omega=0$, and this for a wide 
range of values of the singularity $\zeta$. If $|\zeta|\gg \sqrt{h}$, 
standard steepest descent techniques can be applied, but for 
$|\zeta|\leq \sqrt{h}$, this is no longer possible. We applied a 
polarization approach to treat such cases, i.e. the problem is 
considered on the anti-diagonal of $\mathbb{C}^{2}$ with variables 
$\omega$, $\widetilde{\omega}$ where $\widetilde{\omega}=\overline{\omega}$. 
Here steepest descent contours were considered.
An application of Stokes' theorem allowed for a decomposition of the to 
be studied integral into three terms. For each of these terms the 
dominant contributions for small $h$ for various domains of $|\zeta|$ 
were presented.\\
The current study was for the signature $(+,-)$ of the Hessian matrix at 
the stationary point. A numerical study of interesting cases of these 
integrals appearing in the inverse scattering approach to the 
Davey-Stewartson II equations shows that other possible signatures and 
even degenerate cases are important in this context as well. The 
application of the present approach to these cases will be the 
subject of further research.

\begin{appendices}

\section{Review of wavefront sets for Remark \ref{1sto2}}\label{wf}
\setcounter{equation}{0} 

We here complete Remark \ref{1sto2} by
reviewing briefly some elements of the theory of wavefront sets. For
more details see e.g. \cite[Chapter 8, Section 8.2]{Ho83_01} and also \cite[Chapter 7]{GrSj94_01}. If
$X \subset {\mathbb{R}}^n$ is open, $u\in {\mathcal D}'(X )$,
$(x_0,\xi _0)\in X \times \big({\mathbb{R}}^n\setminus\{0\}\big)$, we say that
$(x_0,\xi _0)\not\in \mathrm{WF}(u)$ if there exist
$\phi \in C_0^\infty (X )$ with $\phi (x_0)\ne 0$ and a conic
neighborhood $V$ of $\xi _0$ in ${\mathbb{R}}^n\setminus\{0\}$ such that
$\widehat{\phi u}(\xi )={\mathcal O}(\langle \xi \rangle^{-\infty })$ in
$V$.  Here $\widehat{\phi u}={\mathcal F}(\phi u)$ denotes the Fourier
transform of $\phi u$, and we use the notation
$\langle \xi \rangle = (1+|\xi |^2)^{1/2}$. This defines a closed
conic subset $\mathrm{WF}(u)\subset X \times 
\big({\mathbb{R}}^n\setminus\{0\}\big)$, and the
definition extends to the case when $X $ is a smooth manifold where
now $\mathrm{WF}(u)\subset T^*X \setminus 0$ with $0$ here denoting
the zero section. If $\Omega $ is a differential form, we define
$\mathrm{WF}(\Omega )$ as the intersection of the wavefront sets of
the coefficients, when expressing $\Omega $ in local coordinates.

\medskip\par 
If $Y$ is a closed submanifold of the manifold $X$ with
$\mathrm{dim}Y=k<n=\mathrm{dim}X$ and $u\in {\mathcal D}'(X)$, then the
restriction ${{u}_\vert}_{Y}$ is well defined in ${\mathcal D}'(Y)$
provided that
\begin{equation}\label{1.wf1}
\mathrm{WF}(u)\cap N^*Y=\emptyset\,,
\end{equation}
where $N^*Y$ denotes the conormal bundle of $Y\subset X$. Moreover,
under the assumption (\ref{1.wf1}) we have
\begin{equation}\label{2.wf1}
\mathrm{WF}\left({u}_{{\vert}_{Y}} \right)\subset \iota
^*\Big(\mathrm{WF}(u)\cap T_Y^*X \Big)\,,
\end{equation}
where $\iota ^*$ denotes the natural projection $T^*_{Y}X\to
T^*Y$.

\smallskip\par 
The term "well defined" can be given a precise sense using
approximation of distributions with smooth functions in semi-normed spaces adapted to wavefront sets.

\smallskip\par 
We recall a simple proof of why ${u}_{{\vert}_{Y}}$ is well defined 
under the assumption (\ref{1.wf1}). This is a local statement, and we are free to choose suitable local coordinates on $X$ and on $Y$. We choose them near a point of $Y$ of the form $x=(x',x'')$, $x'\in {\mathbb{R}}^k$, $x''\in {\mathbb{R}}^{n-k}$, so that $Y$ (locally near the given
point) becomes $x''=0$ and the given point $x=0$.

\smallskip\par 
Let $\phi \in C_0^\infty \big(\mathrm{neigh}(0,{\mathbb{R}}^n)\big)$ satisfy
$\phi (0)\ne 0$. Then $N^*Y$ becomes $x''=0$, $\xi '=0$, $\xi ''\ne 0$ and (\ref{1.wf1}) tells us that if $\mathrm{supp}\phi $ is small enough, then
\begin{equation}\label{3.wf1}
{\mathcal F}(\phi u)(\xi )={\mathcal O}(\langle \xi \rangle^{-\infty })\hbox{
  when }|\xi '|\le |\xi |/C\,,
\end{equation}
if $C>0$ is large enough. By Fourier inversion
we have (formally to start with)
$$
{\phi u}{{_\vert}_{Y}}(x')=\phi u(x',0)=\frac{1}{(2\pi )^n}\int
{\mathcal F}(\phi u)(\xi )\, \e^{\mathrm{i} x'\cdot \xi'} \, \mathrm{d}\xi\,,
$$ 
and hence if ${\mathcal F}'$ denotes the Fourier transform on 
${\mathbb{R}}^k$:
\begin{equation}\label{1.wf2}
{\mathcal F}'\left({\phi u}_{\vert_{Y}} \right)(\xi ')=\frac{1}{(2\pi)^{n-k}}\int{\mathcal F}(\phi u)(\xi) \, \mathrm{d}\xi ''\,
\end{equation}
the integral converges for all $\xi '$ and gives a function of
temperate growth in $\xi '$. This leads to the definition ${\phi u}_{\vert_{Y}}\in {\mathcal E}'$ and to ${u}_{\vert_{Y}}\in {\mathcal D}'$. Further easy estimates give (\ref{2.wf1}). 

\medskip\par 
Consider a smooth map $f:X\longrightarrow Y$, where $X$, $Y$ are two open sets in ${\mathbb{R}}^k$, ${\mathbb{R}}^n$ respectively or more generally manifolds of dimension $k$, $n$. 
Let $u$ be a distribution on $Y$. To define
$u\circ f$ is equivalent to defining 
$({1_X\otimes u)}_{\vert_{\mathrm{graph}(f)}}$, 
where $(1_X\otimes u )(x,y):=u(y)$
and $\mathrm{graph}(f)=\big\{(x,f(x));\, x\in X \big\}$. 
We check that
\begin{equation}\label{2.wf2}
N^*_{(x,f(x))}\big(\mathrm{graph} (f)\big)=\big\{ (\xi ,\eta );\, \xi=-{}^t\hskip-2pt{f'}(x)\eta  \big\}\,,
\end{equation}
\begin{equation}\label{3.wf2}
\mathrm{WF}(1_X\otimes u)\subset \big\{ (x,0;y,\eta );\, (y,\eta )\in\mathrm{WF}(u) \big\}\,,
\end{equation}
hence
\begin{equation}\label{4.wf2}\begin{split}
  \mathrm{WF}&(1_X\otimes u)\cap N^*(\mathrm{graph}(f))\\
 =&\Big\{ (x,0;f(x),\eta );\, x\in X,\ (f(x),\eta )\in \mathrm{WF}(u), \ 0={}^t\hskip-2pt{f'}(x)\eta  \Big\}\,.
  \end{split}
\end{equation}
In particular,
\begin{equation}\label{5.wf2}
\mathrm{WF}\Big(1_X\otimes u\Big)\cap N^*\Big(\mathrm{graph}(f)\Big)=\emptyset\,,
\end{equation}
if
\begin{equation}\label{6.wf2}
{\mathcal N}\big({}^t\hskip-2pt{f'(x)}\big)\cap \Big(\mathrm{WF}(u)\cap T^*_{f(y)}Y\Big)=\emptyset\,,\quad \forall x\in X\,.
\end{equation}
Here ${\mathcal N}\big({}^t\hskip-2pt{f'(x)}\big)$ is the kernel of the linear transformation
$\big[{}^t\hskip-2pt{f'(x)}\big]$. See \cite[Theorem 8.2.4, page 263]{Ho83_01}.

\medskip\par 
The map $\iota^*$ can be identified with that of restriction of
differential 1-forms to the submanifold. When the submanifold is 
$\mathrm{graph}(f)$, and we identify 1-forms on that submanifold with 1-forms in $\mathrm{d}x$ we get at a point $(x,f(x))$:
\begin{equation}\label{7.wf2}
\iota^*(\xi \cdot \mathrm{d}x+\eta \cdot \mathrm{d}y)=\xi \cdot \mathrm{d}x+\eta \cdot  f'(x)\mathrm{d}x=\Big(\xi +{}^t\hskip-2pt{f'}(x)\eta \Big)\cdot \mathrm{d}x\,.
\end{equation}

\par\noindent
We combine this with (\ref{3.wf2}) and get
\begin{equation}\label{8.wf2}
  \mathrm{WF}\left( {(1_X\otimes u)}_{\vert_{\mathrm{graph}(f)}} \right)
  \subset
  \Big\{ (x,{}^t\hskip-2pt{f'}(x)\eta );\, (f(x),\eta )\in \mathrm{WF}(u) \Big\}\,,
\end{equation}
when (\ref{6.wf2}) holds.

\medskip\par 
From the identification 
${1_X\otimes u}_{\vert_{\mathrm{graph}(f)}}\simeq u\circ f$, 
we get the equivalent inclusion
\begin{equation}\label{1.wf3}
\mathrm{WF}(u\circ f)\subset \Big\{ (x,{}^t\hskip-2pt{f'}(x)\eta );\, (f(x),\eta )\in \mathrm{WF}(u) \Big\}\,,
\end{equation}
when (\ref{6.wf2}) holds.

\medskip\par 
Now consider a smooth map $f=f(t,x):\, I\times W\to Y\subset {\mathbb{R}}^n$ where $f(t,\cdot )$ is injective with injective differential for every
$t\in I$. Let $u \in {\mathcal D}'(Y) $ and assume that
(\ref{5.sto2}) holds with $\Gamma _t={\mathcal R}(f(t,\cdot ))$. This 
means that (\ref{6.wf2}) holds for $f(t,\cdot )$, and
a fortiori it also holds for $f$. Then $u \circ f(t,\cdot )\in
{\mathcal D}'(W)$ is well defined for every $t\in I$, and $u \circ f$
is well defined in ${\mathcal D}'(I\times W)$.

\medskip\par 
We conclude from (\ref{1.wf3}) that with $f=f(t,x)$ as above,
\begin{equation}\label{2.wf3}
  \mathrm{WF}(u \circ f)\subset
  \Big\{ (t,x;{}^t\hskip-2pt {(\partial _tf})\eta ,{}^t\hskip-2pt {(\partial _x}f)\eta )
;\ (f(t,x),\eta )\in \mathrm{WF}(u ) \Big\}\,.
\end{equation}
Here ${}^t\hskip-2pt{\big(\partial_x f\big)}\eta \ne 0$, and we conclude that
\begin{equation}\label{3.wf3}
(t,x;\tau ,\xi )\in \mathrm{WF}\big(u \circ f\big) \, \Longrightarrow \, \xi \ne 0\,.
\end{equation}

\medskip\par 
Restrict $t$ to a compact sub-interval. We may  decompose $u
\circ f$ by means of a partition of unity in the $x$ variable. Let $\chi \in C_0^\infty (I\times W)$. We claim that
\begin{equation}\label{4.wf3}
I(t):=\int \chi (t,x)(u \circ f)(t,x)\, \mathrm{d}x
\end{equation}
is a smooth function of $t$. To prove this, let $\widetilde{\chi
}(x)\in C_0^\infty $ be equal to 1 on $\mathrm{supp}(\chi )$. Then
\begin{equation}\label{5.wf3}
  \begin{split}
    I(t)=&\int \widetilde{\chi }(x)(\chi u \circ f)(t,x)\, \mathrm{d}x\\
    =&\int \widetilde{\chi }(x)\iintJ {\mathcal F}(\chi u \circ
    f)(\tau ,\xi )\, \e^{\mathrm{i}(t\tau +x\cdot \xi )}\frac{\, \mathrm{d}\tau \, \mathrm{d}\xi }{(2\pi )^{k+1}} \, \mathrm{d}x\\
=&\iintJ ({\mathcal F}'\widetilde{\chi })(-\xi ) {\mathcal F}(\chi u \circ
    f)(\tau ,\xi ) \, \e^{\mathrm{i} t\tau}\frac{\, \mathrm{d}\tau \, \mathrm{d}\xi }{(2\pi )^{k+1}}\,,
  \end{split}
\end{equation}
where ${\mathcal F}$ denotes the Fourier transform on ${\mathbb
  R}^{k+1}_{t,x}$ and ${\mathcal F}'$ that on ${\mathbb{R}}^k_x$.

\medskip\par 
By (\ref{3.wf3}) we have
\begin{equation}\label{1.wf4}
{\mathcal F}(\chi u \circ f)(\tau ,\xi )={\mathcal O}\Big(\langle (\tau ,\xi)
\rangle^{-\infty }\Big)\,, \ \hbox{ when }|\xi |\le |(\tau ,\xi )|/C\,,
\end{equation}
for some large enough $C>0$. We split the last integral  in
(\ref{5.wf3}) into two terms by restricting the integration to $|\xi
|\le |(\tau ,\xi )|/C$ and $|\xi |> |(\tau ,\xi )|/C$ respectively. In
the domain of integration for the first term we have

$$
({\mathcal F}'\widetilde{\chi })(-\xi ){\mathcal F}(\chi u \circ f)(\tau
,\xi )={\mathcal O}\Big(\langle (\tau ,\xi )\rangle ^{-\infty }\Big)\,.
$$
In that of the second term we have $\langle \xi \rangle\asymp \langle
(\tau ,\xi )\rangle$, hence for some $N_0\ge 0$,
$$
({\mathcal F}'\widetilde{\chi })(-\xi ){\mathcal F}(\chi u \circ f)(\tau
,\xi )={\mathcal O}(1)\langle \xi \rangle ^{-\infty }\langle (\tau ,\xi
)\rangle ^{N_0 }
={\mathcal O}(1)\langle (\tau ,\xi )\rangle ^{N_0-\infty }\,.
$$
It follows that the last integral in (\ref{5.wf3}) converges and is a
smooth function of $t$. Thus,
\begin{equation}\label{2.wf4}
I(t)\in C^\infty\,.
\end{equation}
This gives the smoothness of $\int_{\Gamma _t}\phi \Omega $ in Remark \ref{1sto2}.

\section{The special functions \texorpdfstring{$G$}{}, some complements.}
\label{spe}
\setcounter{equation}{0}
In this appendix we collect some basic facts on the 
special function $G$ being the Hilbert transform of the Dawson 
integral.

Let $Q(\omega)=q \,\dfrac{\omega^2}{2}$, $q\in \mathbb{C}\setminus\{0\}$. Assume
\begin{equation}\label{spe.1}
\mathrm{Im} Q(\omega)\asymp |\omega|^2,\quad \omega\in \e^{\mathrm{i}\frac{\pi}{4}}\mathbb{R}\,. 
\end{equation}
Equivalently, $\mathrm{Re} \,q >0$. Recall that 
\begin{equation}\label{spe.2}
 G_{iQ}^{l/r}(\zeta  )=\frac{1}{2\pi \mathrm{i}}\int_{\Gamma} \frac{1}{\zeta-\omega} \, \e^{\mathrm{i} Q(\omega)} \, \mathrm{d}\omega\,,
\end{equation}
where $\Gamma$ is a deformation of $\e^{\mathrm{i}\frac{\pi}{4}}\mathbb{R}$
oriented North East, avoiding $\zeta$ with this point to the left/right.

\par\medskip
Change of variable: $\omega=k\widehat \omega$, to obtain
\begin{equation}\label{spe.3}
 \mathrm{i} Q(\omega)=-\frac{\widehat\omega^2}{2} \Leftrightarrow
 \mathrm{i} q k^2 \widehat\omega^2=-\widehat\omega^2\Longleftrightarrow
 \mathrm{i} q k^2=-1\Leftrightarrow q k^2=\mathrm{i}\,. 
\end{equation}
Can take
\begin{equation}\label{spe.4}
 k=\frac{\e^{\mathrm{i} \frac{\pi}{4}}}{\sqrt{q}}\, 
\end{equation}
with the branch of the square root satisfying $\arg({\sqrt{q}})\in]-\frac{\pi}{4},\frac{\pi}{4}[$, so that $\arg({\frac{1}{\sqrt{q}}})\in]-\frac{\pi}{4},\frac{\pi}{4}[$. Then $\arg(k)\in]0,\frac{\pi}{2}[$.

\par\medskip\noindent
Let $G^{l,r}=G^{l,r}_{-\frac{\omega^2}{2}}$ be as in \cite{KlSjSt23_02}, modulo the multiplicative factor $(2\pi \mathrm{i})^{-1},$ and with $\e^{\mathrm{i} \frac{\pi}{4}}\mathbb{R}$ replaced by $\mathbb{R}$. We get 
\begin{align}\label{spe.5}
G^{l,r}_{\mathrm{i} Q}(\zeta)&=\frac{1}{2\pi \mathrm{i}}\int_{k^{-1}\Gamma}\frac{1}{\zeta-k\widehat\omega} \e^{-\frac{\widehat\omega^2}{2}}\, \mathrm{d}(k\widehat\omega)\nonumber\\
& =\frac{1}{2\pi \mathrm{i}}\int_{k^{-1}\Gamma}\frac{1}{\frac{1}{k}\zeta-\widehat\omega} \e^{-\frac{\widehat\omega^2}{2}}\, \mathrm{d}(\widehat\omega)\nonumber\\
&=G^{l,r}\Big(\frac{\zeta}{k}\Big)\,.
\end{align}
For the last equality, we see that $k^{-1}\Gamma$ can be deformed to a small deformation of $\mathbb{R}$.
Let $\zeta\in \mathbb{R}$:
\begin{align*}
 G^l(\zeta)&=\lim_{\varepsilon\searrow 0}G^l(\zeta+\mathrm{i} \varepsilon)=\lim_{\varepsilon\searrow 0}\frac{1}{2\pi \mathrm{i}}\int_{\mathbb{R}} \frac{1}{\zeta+\mathrm{i}\varepsilon-\omega} \e^{-\frac{\omega^2}{2}}\, \mathrm{d}\omega\,,\\
 G^r(\zeta)&=\lim_{\varepsilon\searrow 0}G^r(\zeta-\mathrm{i} \varepsilon)=\lim_{\varepsilon\searrow 0}\frac{1}{2\pi \mathrm{i}}\int_{\mathbb{R}} \frac{1}{\zeta-\mathrm{i}\varepsilon-\omega} \e^{-\frac{\omega^2}{2}}\, \mathrm{d}\omega\,.
\end{align*}
Here we encounter $\dfrac{1}{x\pm \mathrm{i} 0}\in \mathcal{S}'(\mathbb{R})$ and $\dfrac{1}{x+\mathrm{i} 0}+\dfrac{1}{x-\mathrm{i} 0}=2\mathrm{v.p}\Big(\frac{1}{x}\Big)$, 
so 
\begin{equation}\label{spe.6}
\frac{1}{2}\Big(G^l(\zeta)+G^r(\zeta)\Big)=\frac{1}{2\pi \mathrm{i}}\Big[(\mathrm{v.p}\frac{1}{x})*(\e^{-\frac{x^2}{2}})\Big](\zeta)\,. 
\end{equation}
The distribution $\mathrm{v.p}\Big(\frac{1}{x}\Big)$ is odd so $\frac{1}{2}\big(G^l+G^r\big)$
is odd on $\mathbb{R}$ and in particular $(G^l+G^r)(0)=0$. Let us check the last fact directly:
$$
G^l(0)+G^r(0)=\lim_{\varepsilon\searrow 0} \frac{1}{2\pi \mathrm{i}}\int_{\mathbb{R}} \, \underbrace{\Big(\frac{1}{\mathrm{i}\varepsilon-\omega}+\frac{1}{-\mathrm{i}\varepsilon-\omega}\Big)}_{=:g(\omega,\varepsilon)} \, \e^{-\frac{\omega^2}{2}}\, \mathrm{d}\omega\,.
$$
\begin{equation}\label{spe.7}\begin{cases}
    g(-\omega,\varepsilon)=\frac{1}{\mathrm{i}\varepsilon+\omega}+\frac{1}{-\mathrm{i}\varepsilon+\omega}\\
    {}\\
-g(\omega,\varepsilon)=\frac{1}{-\mathrm{i}\varepsilon+\omega}+\frac{1}{\mathrm{i}\varepsilon+\omega}
  \end{cases}
  \Longrightarrow \quad g(-\omega,\varepsilon)=-g(\omega,\varepsilon)\,.
\end{equation}

\medskip\noindent
Thus $g(\cdot,\varepsilon)$ is odd and hence $G^l(0)+G^r(0)=0$. This extends to $G^l_{-\mathrm{i} Q}(0)+G^r_{-\mathrm{i} Q}(0)=0$ by \eqref{spe.5}. Moreover, 
\begin{equation}\label{spe.8}
\begin{split}
   G^r(0)-G^l(0) & =\lim_{\varepsilon\searrow 0} \, \frac{1}{2\pi \mathrm{i}}\int_{\mathbb{R}}\Big( \frac{1}{-\mathrm{i}\varepsilon-\omega} -\frac{1}{\mathrm{i}\varepsilon-\omega}\Big)\, \e^{-\frac{\omega^2}{2}}\, \mathrm{d}\omega\\
   &= \, \lim_{\varepsilon\searrow 0} \frac{\varepsilon}{\pi}\int_{\mathbb{R}}
   \frac{\e^{-\frac{\omega^2}{2}}}{\omega^2+\varepsilon^2} \, \mathrm{d}\omega\\
   &=\, \lim_{\varepsilon\searrow 0} \frac{1}{\pi}\int_{\mathbb{R}}
   \frac{\e^{-\varepsilon^2\frac{t^2}{2}}}{t^2+1} \, \mathrm{d}t\\
   &=\frac{1}{\pi}\int_{\mathbb{R}}
   \frac{1}{t^2+1} \, \mathrm{d}t=1\,.
\end{split} 
\end{equation}
Then,  
\begin{equation}\label{spe.9}
G^r(0)=-G^l(0)=\frac{1}{2}\,.
\end{equation}
For more details concerning the properties of the special function $G$ we can consult \cite[Chapter 7, page 295]{AbSt84_01}\footnote{We thank Jan Dereziński for drawing our attention to the inhomogeneous differential equation satisfied by the special function $G$.}.

\section{Invariance of \texorpdfstring{$\mathrm{I } (\zeta ,1)$}{} under rotations}
\label{inv}
\setcounter{equation}{0}
In this appendix we study the invariance of the 
studied integral under rotations of the steepest descent contours in 
the polarization approach.

Let $\phi \in C^\infty \big(\mathrm{neigh}(0;{\mathbb R}^2);{\mathbb R}\big)$
with
\begin{equation}\label{1.inv1}
\phi (0)=0\,,\,\, \phi '(0)=0\,,\,\,\, \phi ''(0) \hbox{ non-degenerate of
  signature } (+,-)\, .
\end{equation}
Let $\phi _2$ be the quadratic Taylor polynomial of $\phi $ at $0$,
\begin{equation}\label{2.inv1}
  \phi (\xi ,\eta )=\phi _2(\xi ,\eta )+{\mathcal O}((\xi ,\eta )^3)\,,\quad
  (\xi ,\eta )\longrightarrow  0\,.
\end{equation}
Write \begin{equation}\label{3.inv1}
\phi _2(\xi ,\eta )=\frac{\lambda }{2}\xi ^2-\frac{\mu }{2}\eta
^2+\rho \xi \eta\,,\quad \lambda ,\mu ,\rho \in {\mathbb R}.
\end{equation}
Observe that
\begin{lemma}\label{1inv1}
  We have the following equivalences
  \begin{equation}\label{4.inv1}
\lambda >0 \Longleftrightarrow \phi _2(\xi ,0)\asymp \xi ^2
\end{equation}
and
\begin{equation}\label{5.inv1}
\mu >0 \Longleftrightarrow \phi _2(0,\eta) \asymp -\eta ^2\,.
\end{equation}
\end{lemma}

\par 
The set $S^+_{\phi _2}:=\phi_2^{-1}\big(]0,+\infty [\big)$ is the
disjoint union of two convex mutually antipodal sectors and $S^-_{\phi
_2}:=\phi_2^{-1}\big(]-\infty ,0[\big)$ has the same property. 
The conditions (\ref{4.inv1}), (\ref{5.inv1}) amount to $\dot{\mathbb R}u_1\subset
  S^+_{\phi _2}$, $\dot{\mathbb R}u_2\subset
  S^-_{\phi _2}$, where $u_1=(1,0)$, $u_2=(0,1)$.

  \par\medskip
  Let $\mathbf{U}:{\mathbb R}^2\longrightarrow {\mathbb R}^2$ be a rotation of the form
  $\mathbf{U}={\mathcal R}_\alpha$ as in (\ref{2.inv3}) below. ($\mathbf{U}$ is orthogonal with
  $\det \mathbf{U}=1$.) Let $\widetilde{\xi }$, $\widetilde{\eta }$ be the
  coordinates on ${\mathbb R}^2 $ for the orthonormal basis $(\mathbf{U}u_1, \mathbf{U}u_2)$,
  so that
 $$
\xi\, u_1+\eta\, u_2=\widetilde{\xi } \, \mathbf{U}u_1+\widetilde{\eta} \, \mathbf{U}u_2\,.
$$
Then
$$
\phi _2=\frac{\widetilde{\lambda }}{2} \, {\widetilde{\xi
}\,}^2-{\,}\frac{\widetilde{\mu }}{2} \, {\widetilde{\eta }\,}^2+ \, \widetilde{\rho
} \, \widetilde{\xi } \, \widetilde{\eta } \quad \hbox{ with } \, \widetilde{\lambda },\,
\widetilde{\mu }>0\,,
$$
if and only if 
\begin{equation}\label{6.inv1}
\mathbf{U}u_1\in S^+_{\phi _2}\,,\quad  \mathbf{U}u_2\in S^-_{\phi_2}\,.
\end{equation}

\par\smallskip
Assume that $\lambda ,\,\mu >0 $ in(\ref{3.inv1}) and let $\mathcal{C}_+$ denote the component of $S^+_{\phi _2}$ which contains ${\mathbb R}_+\,u_1$. Identifying ${\mathbb R}^2_{\xi ,\eta }\simeq {\mathbb C}_{\xi+i\eta }$, we have
\begin{equation}\label{1.inv2}
\mathcal{C}_+=\big\{\omega \in {\mathbb C}\setminus\{0\}\,;\, \gamma _1<\mathrm{arg}(\omega)
<\gamma _2 \big\}\,,
\end{equation}
where
\begin{equation}\label{2.inv2}
-\frac{\pi }{2}<\gamma _1<0<\gamma _2 <\frac{\pi }{2}\,.
\end{equation}
If $\mathbf{U}u_1\in \mathcal{C}_+$, then $\mathbf{U}u_2=\e^{\mathrm{i}\frac{\pi}{2}}\mathbf{U}u_1$ satisfies
\begin{equation}\label{3.inv2}
\mathrm{arg}(\mathbf{U}u_2)\in ]\gamma _1+\frac{\pi }{2},\gamma _2+\frac{\pi }{2}[\,,
\end{equation}
and $\mathbf{U}u_2\in S^-_{\phi _2}$ if and only if  $\mathrm{arg}(\mathbf{U}u_2)\in ]\gamma _2,\gamma_1+\pi [$, i.e. if and only if
\begin{equation}\label{4.inv2}
\mathrm{arg}(\mathbf{U}u_1)\in ]\gamma _2-\frac{\pi
  }{2},\gamma _1+\frac{\pi }{2}[\,.
\end{equation}
We conclude that
\begin{equation}\label{5.inv2}
\mathbf{U}u_1\in S^+_{\phi_2}\quad \hbox{ and }\quad \mathbf{U}u_2\in S^-_{\phi_2}
\end{equation}
if and only if 
\begin{equation}\label{6.inv2}
  \mathrm{arg}(\mathbf{U}u_1)\in \left]\max (\gamma_1,\gamma_2-\frac{\pi }{2}),\min
(\gamma _2,\gamma _1+\frac{\pi }{2})\right[=:J\,.
\end{equation}
The set of rotations $\mathbf{U}=\mathcal{R}_\alpha $ satisfying (\ref{6.inv2}) is
connected; the corresponding set of $\alpha $ is the interval in
(\ref{6.inv2}).

\par\medskip 
Under the assumption that
\begin{equation}\label{7.inv2}
\zeta \in \mathbb{C}\setminus\e^{\mathrm{i}\frac{\pi}{4}}{\mathbb R}\,,
\end{equation}
we have that $I(\zeta ):=\mathrm{I}(\zeta ,1)$ is given by
\begin{equation}\label{1.inv3}
 I(\zeta )=-\frac{1}{\pi}\iintJ_{\Gamma_{1}}\frac{1}{\omega-\zeta} \, 
a(\omega ,\widetilde{\omega };h)\chi(\omega,\widetilde{\omega}) \, \e^{\frac{\mathrm{i}}{h}\phi(\omega ,\widetilde{\omega })} \, \frac{\mathrm{d}\widetilde{\omega }\wedge \mathrm{d}\omega}{2\mathrm{i}} \,,
\end{equation}
where we recall that $\omega =\xi +i\eta $, $\widetilde{\omega }=\xi
-i\eta $ and $\Gamma_1=\e^{\mathrm{i}\frac{\pi}{4}}{\mathbb R}_\omega \times
\e^{\mathrm{i}\frac{\pi}{4}}{\mathbb R}_{\widetilde{\omega }}$ in ${\mathbb C}^2_{\omega
  ,\widetilde{\omega }}$. If $\zeta \ne 0$ belongs to a small conic
neighborhood of $\e^{\mathrm{i}\frac{\pi}{4}}\dot{{\mathbb R}}$, we shall use $\mathbf{U}=\mathcal{R}_\alpha$
for a suitable $\alpha \in J$ to achieve that a natural effective $\check{\zeta}$ avoids
a conic neighborhood of $\e^{\mathrm{i}\frac{\pi}{4}}\dot{{\mathbb R}}$.

\par\medskip 
In ${\mathbb R}^2$ we can identify $\mathcal{R}_\alpha $ with
\begin{equation}\label{2.inv3}
  \mathcal{R}_\alpha =\begin{pmatrix}\cos \alpha &-\sin \alpha \\
  \sin \alpha & \cos \alpha \end{pmatrix}\,.
\end{equation}
For $\mathbf{U}=\mathcal{R}_\alpha $ with $\alpha $ in the interval $J$ in
(\ref{6.inv2}), let $(\check{\xi },\check{\eta })$ denote the linear
coordinates  of a point $\rho =(\xi ,\eta )\in {\mathbb R}^2$ with respect
to the basis $(\mathbf{U}u_1, \mathbf{U}u_2)$. Then
\begin{equation}\label{3.inv3}
\begin{pmatrix}\xi \\
  \eta \end{pmatrix}=\mathcal{R}_\alpha \begin{pmatrix}\check{\xi }\\
  \check{\eta }\end{pmatrix}\,.
\end{equation}
In analogy with the definition of $\omega ,\, \widetilde{\omega }$,
we define $\check{\omega },\, \check{\widetilde{\omega }}$ by
\begin{equation}\label{4.inv3}
\check{\omega }= \check{\xi }+i  \check{\eta }\,,\qquad
\check{\widetilde{\omega }}=  \check{\xi }-i  \check{\eta }\,.
\end{equation}
Then (\ref{3.inv3}) gives
\begin{equation}\label{5.inv3}
\omega =\e^{\mathrm{i}\alpha}\,\check{\omega }\,,\qquad \widetilde{\omega}=\e^{-\mathrm{i}\alpha}\,\check{\widetilde{\omega }}\,.
\end{equation}
If $\check{a}$, $\check{\chi }$, $\check{\phi }$ denote the functions
$a$, $\chi $, $\phi $ as functions of $(\check{\omega
},\check{\widetilde{\omega }})$, so that 
$$
\check{ a}(\check{\omega } ,
\check{\widetilde{\omega }};h)=a(\omega ,\widetilde{\omega };h)
$$ 
and similarly for $\chi $, $\phi $, and $\zeta \not\in \e^{\mathrm{i}\frac{\pi}{
  4}}\dot{{\mathbb R}}$, then 
\begin{equation}\label{1.inv4}\begin{split}
  I(\zeta )&=
-\frac{1}{\pi}\iintJ_{\check{\Gamma}} \, \frac{1}{\e^{\mathrm{i}\alpha}\check{\omega}-\zeta} \, 
\check{a}(\check{\omega},\check{\widetilde{\omega}};h) \check{\chi}(\check{\omega } ,
\check{\widetilde{\omega }}) \, \e^{\frac{\mathrm{i}}{h}\check{\phi}(\check{\omega } ,
\check{\widetilde{\omega }})} \,
\frac{\mathrm{d}\check{\widetilde{\omega }}\wedge
  \mathrm{d}\check{\omega }}{2\mathrm{i}}\\
&=-
\e^{-\mathrm{i}\alpha }\frac{1}{\pi}\iintJ_{\check{\Gamma}} \, \frac{1}{\check{\omega}-\check{\zeta}} \, 
\check{a}(\check{\omega},\check{\widetilde{\omega}};h) \check{\chi}(\check{\omega } ,
\check{\widetilde{\omega }}) \, \e^{\frac{\mathrm{i}}{h}\check{\phi}(\check{\omega } ,
\check{\widetilde{\omega }})} \,
\frac{\mathrm{d}\check{\widetilde{\omega }}\wedge
  \mathrm{d}\check{\omega }}{2\mathrm{i}}\,,
\end{split}
\end{equation}
where $\check{\zeta }=\e^{-\mathrm{i}\alpha }\zeta $ and
\begin{equation}\label{2.inv4}
\check{\Gamma }=\Big\{(\check{\omega },\check{\widetilde{\omega }})\in \mathbb{C}^2;\,\,\, 
(\e^{\mathrm{i}\alpha }\check{\omega },\e^{-\mathrm{i}\alpha}\check{\widetilde{\omega
  }})\in \Gamma_1 \Big\}\,.
\end{equation}

\par 
Now drop the assumption that $\zeta \not\in \e^{\mathrm{i}\frac{\pi}{4}}\dot{{\mathbb 
    R}}$ and assume instead that
\begin{equation}\label{3.inv4}
\check{\zeta } \mbox{ is outside a conic neighborhood of } \e^{\mathrm{i}\frac{\pi}{4}}\dot{\mathbb
    R}\,\, \mathrm{in}\,\, \dot{\mathbb C}\, ,
\end{equation}
i.e.\
\begin{equation}\label{4.inv4}
\mathrm{arg}(\zeta) \not\in \mathrm{neigh}(\alpha +\frac{\pi }{4};{J})
\end{equation}
and keep the assumption (\ref{6.inv2}):
\begin{equation}\label{5.inv4}
\alpha \in J\, .
\end{equation}
We can then make the change of variables (\ref{3.inv3}) in the
original integral in (\ref{1.sto3}), (\ref{3.sto3}) and get first as in \eqref{1.inv4}
\begin{align}\label{6.inv4}
I(a\chi ,\phi ,\zeta ;h)=&-\frac{1}{\pi}\iintJ_{\mathbb R^2}\,\frac{1}{
\e^{\mathrm{i}\alpha }(\check\xi +\mathrm{i} \check\eta )-\zeta}\, \check a(\check\xi , \check\eta ;h) \check\chi (\check\xi , \check\eta )\,\e^{\frac{\mathrm{i}}{h}\check\phi (\check\xi
, \check\eta)}\,\mathrm{d}\check \xi \mathrm{d}\check\eta\nonumber\\
=&-\frac{\e^{-\mathrm{i}\alpha}}{\pi}\iintJ_{\mathbb R^2}\,\frac{1}{
(\check\xi +\mathrm{i} \check\eta )-\check{\zeta}}\, \check a(\check\xi , \check\eta ;h) \check\chi (\check\xi , \check\eta )\,\e^{\frac{\mathrm{i}}{h}\check\phi (\check\xi
, \check\eta)}\,\mathrm{d}\check \xi \mathrm{d}\check\eta\nonumber\\
=&: \e^{-\mathrm{i}\alpha}\check I(\check a \check\chi ,\check\phi ,\check \zeta ;h)\,.
\end{align}

\noindent
Then apply the results in Section \ref{Sto} and the following ones to
$\check I(\check a \check\chi ,\check\phi ,\check\zeta ;h)$. We get the
analogue of (\ref{2.sto8}), (\ref{3.sto8}):
\begin{equation}\label{7.inv4}
\e^{\mathrm{i}\alpha}I(a\chi ,\phi ,\zeta ;h)= \check I(\check a \check \chi ,\check \phi, \check\zeta ;h)
=\check {\mathrm{I}}(\check\zeta ,1) +\check {\mathrm{II}}(\check\zeta ,1) +\check {\mathrm{III}}(\check\zeta ,1)
\end{equation}
with the natural definitions of $\check {\mathrm{I}}$, $\check
{\mathrm{II}}$, $\check {\mathrm{III}}$, analogous to those of
$\mathrm{I}$, $\mathrm{II}$, $\mathrm{III}$. In particular,
\begin{equation}\label{1.inv5}
\hskip-5pt\check {\mathrm{I}}(\check\zeta ,1)=-\frac{1}{\pi }\iintJ_{(\check \omega , \check {\widetilde\omega })\in \Gamma_1} \, \frac{1}{\check \omega -\check\zeta}
  \, \check{a}(\check{\omega},\check{\widetilde{\omega}};h) \check{\chi}(\check{\omega } ,
\check{\widetilde{\omega }}) \, \e^{\frac{\mathrm{i}}{h}\check{\phi}(\check{\omega } ,
\check{\widetilde{\omega }})} \,
\frac{\mathrm{d}\check{\widetilde{\omega }}\wedge
  \mathrm{d}\check{\omega }}{2\mathrm{i}}\,, 
\end{equation}
which differs from the first integral in (\ref{1.inv4}) only by the
choice of integration contour, cf.\ (\ref{2.inv4}).

\par\medskip 
We can express (\ref{1.inv5}) as an integral in the $\omega $,
$\widetilde{\omega }$ variables (cf. \eqref{5.inv3}, $\zeta=\mathrm{e}^{\mathrm{i}\alpha}\check\zeta$, 
\begin{equation}\label{2.inv5}
\mathrm{e}^{-\mathrm{i}\alpha}  \check{\mathrm{I}}(\zeta ,1)=-\frac{1}{\pi }\iintJ_{\Gamma_\alpha (1)}
  \frac{\e^{\frac{\mathrm{i}}{h}\phi(\omega,\widetilde \omega)}}{\omega-\zeta } \, a(\omega,\widetilde \omega;h)\,\chi(\omega,\widetilde \omega) \,  \, \frac{\mathrm{d}\widetilde{\omega }\wedge  \mathrm{d}\omega}{2\mathrm{i}}=:\mathrm{I}_\alpha(\zeta,1)\,,
\end{equation}
where
\begin{equation}\label{3.inv5}
\Gamma _\alpha (1)=\big\{(\omega ,\widetilde{\omega })\in \mathbb{C}^2\, ; \,\, (\e^{-\mathrm{i}\alpha
}\omega , \e^{\mathrm{i}\alpha }\widetilde{\omega })\in \Gamma_1 \big\}\,,
\end{equation}
with $\Gamma_0(1)=\Gamma_1$.

\noindent
Keeping the assumptions (\ref{4.inv4}), (\ref{5.inv4}), we
reinstate the assumption that $\zeta \not\in \e^{\mathrm{i}\frac{\pi}{4}}{\mathbb R}$ in
the strengthened variant
$$
\mathrm{arg}(\zeta) \not\in \mathrm{neigh}(\frac{\pi }{4}; J)\,,
$$
and further strengthen the assumptions to
\begin{equation}\label{4.inv5}
\mathrm{arg}(\zeta) \not\in \mathrm{neigh}(K_\alpha ;J)\,,
\end{equation}
where we put
\begin{equation}\label{5.inv5}
K_\alpha =[\frac{\pi }{4},\frac{\pi }{4}+\alpha ]\hbox{ when }\alpha
>0, \, \hbox{ and } K_\alpha =[\frac{\pi }{4}+\alpha ,\frac{\pi }{4}]\hbox{ when
}\alpha \le 0\,.
\end{equation}
Then $\mathrm{I}_{t\alpha} (\zeta ,1)$ in (\ref{2.inv5}) is
well defined for $t\in [0,1] $, since $\omega - \zeta \ne 0$ for
$(\omega ,\widetilde{\omega })\in \Gamma _{t\alpha}(1)$.

\noindent
Along $\Gamma _{t\alpha} (1)$ we have $|\omega -\zeta |\ge |\zeta
|/{\mathcal O}(1)$ uniformly for $t\in [0,1] $. Moreover, $\mathrm{Im}
\phi (\omega ,\widetilde{\omega })\asymp |\omega
|^2+|\widetilde{\omega }|^2$ here, by (\ref{3.ft8}), provided that we
work in a sufficiently small neighborhood of $(0,0)$. Assuming that
$a(\omega ,\widetilde{\omega };h)$ is holomorphic in
$\mathrm{neigh}((0,0);\mathbb{C}^2)$ and that $\chi $ has its support in such a
small neighborhood and equals to 1 in another similar neighborhood, we
can apply Stokes' formula and get when $|\zeta |\ll 1$
\begin{equation}\label{2.inv6}
\mathrm{I}_\alpha (\zeta ,1)-
  \mathrm{I}_0 (\zeta ,1)
  =-\frac{1}{\pi }\int_0^1 \hskip-5pt\mathrm{d}t \iintJ_{\Gamma _{t\alpha}
    (1)}\frac{1}{\omega-\zeta} \, a\, \e^{\frac{\mathrm{i}}{h}\phi} \, \nu \rfloor \frac{\mathrm{d}\chi \wedge
  \mathrm{d}\widetilde{\omega }\wedge \mathrm{d}\omega}{2\mathrm{i}}\,,
\end{equation}
where $\nu $ is a deformation field associated to our family of
contours. As before, $\mathrm{d}\chi \wedge \mathrm{d}\widetilde{\omega }\wedge \mathrm{d}\omega=
\overline{\partial }\chi \wedge \mathrm{d}\widetilde{\omega }\wedge \mathrm{d}\omega
$ is supported away from $(0,0)$, hence
\begin{equation}\label{3.inv6}
a(\omega,\widetilde\omega;h)\, \e^{\frac{\mathrm{i}}{h}\phi(\omega,\widetilde\omega)} \, \nu \rfloor (\mathrm{d}\chi \wedge
  \mathrm{d}\widetilde{\omega }\wedge \mathrm{d}\omega )={\mathcal O}(1)\exp \Big(-\frac{1}{{\mathcal O}(h)}\Big)\,.
\end{equation}

\par 
We use $(t\alpha)$-dependent coordinates as in
\eqref{2.ft1}, to write
\begin{equation}\label{4.inv6}
\Gamma_1\ni (\check\omega ,\check{\widetilde{\omega }})=w_+e_++w_-e_- \,.
\end{equation}
(\ref{3.ft1}) gives $\check \omega =w_-\e^{\mathrm{i}\frac{\pi}{4}}$, where $\omega
=\e^{\mathrm{i}{t\alpha}}\check \omega $, hence $1/(\omega -\zeta)$ in
(\ref{2.inv6}) is
$$
=\frac{1}{\e^{\mathrm{i}({t\alpha} + \frac{\pi}{4})}w_- -\zeta}=\frac{{\mathcal O}(1)}{|w_-|+|\zeta
  |}\,,
$$
which produces a log contribution to the integral in (\ref{2.inv6})
and we get
\begin{equation}\label{1.inv7}
\mathrm{I}_\alpha (\zeta ,1)-\mathrm{I}_0 (\zeta ,1)={\mathcal
  O}(1)\, \e^{-\frac{1}{{\mathcal O}(h)}}\ln (\frac{1}{|\zeta |})\,.
\end{equation}
The factor $\ln (1/|\zeta |)$ in (\ref{1.inv7}) can be dropped if
$\ln (1/|\zeta |)\le \e^{1/{\mathcal O}(h)}$ i.e. if
\begin{equation}\label{2.inv7}
|\zeta |\ge \exp(-\exp (1/{\mathcal O}(h)))\,.
\end{equation}
Then, we obtain the following result
\begin{prop}
Let $\alpha\in J$ where $J$ is the interval in \eqref{6.inv2}. 
Assume that $\zeta$ satisfies \eqref{4.inv5} and \eqref{2.inv7}. Then 
$$
\mathrm{I}_\alpha (\zeta ,1)=\mathrm{I} (\zeta ,1)+{\mathcal
  O}(1)\, \e^{-\frac{1}{{\mathcal O}(h)}}\,,
$$
we recall that $\mathrm{I} (\zeta ,1)$ is the first term in the decomposition \eqref{3.sto8} (see also Proposition \ref{1ft3}), and $I_\alpha(\zeta,1)=\mathrm{e}^{-\mathrm{i}\alpha}\check{\mathrm{I}}(\check\zeta,1)$ the corresponding term in \eqref{7.inv4}, \eqref{1.inv5}, \eqref{2.inv5}.
\end{prop}

\begin{remark}
We have shown above that the asymptotics of the first term obtained via Stokes' formula for the decomposition of integral \eqref{6.inv4} is invariant under rotations. The same arguments prove that the other two integrals, $\mathrm{II}$ and $\mathrm{III},$ obtained via the Stokes formula are also invariant under rotations.

\end{remark}

\end{appendices}



\printbibliography

\end{document}